\documentclass{amsart}
\usepackage{comment}
\usepackage{version}
\usepackage{color}
\newenvironment{NB}{
\color{red}{\bf NB}. \footnotesize 
}{}
\excludeversion{NB}

\usepackage{amsmath,amssymb}
\usepackage[arrow, matrix]{xy}
\usepackage[dvips]{graphicx}
\usepackage{amsthm}

\newtheorem{thm}{Theorem}[section]

\newtheorem{cor}[thm]{Corollary}
\newtheorem{prop}[thm]{Proposition}
\newtheorem{defn}[thm]{Definition}
\newtheorem{lem}[thm]{Lemma}
\newtheorem{rem}[thm]{Remark}

\numberwithin{thm}{section}

\def\mbi#1{\boldsymbol{#1}}

\def\Ext{\text{Ext}}

\def\C{\mathbb{C}}
\def\R{\mathbb{R}}

\def\e{\varepsilon}
\def\Q{\mathbb{Q}}

\def\min{\mathop{\mathrm{min}}\nolimits}

\def\Hom{\mathop{\mathrm{Hom}}\nolimits}

\def\im{\mathop{\mathrm{im}}\nolimits}

\def\ker{\mathop{\mathrm{ker}}\nolimits}
\def\coker{\mathop{\mathrm{coker}}\nolimits}
\def\dim{\mathop{\mathrm{dim}}\nolimits}

\def\deg{\mathop{\mathrm{deg}}\nolimits}
\def\rank{\mathop{\mathrm{rank}}\nolimits}
\def\rk{\mathop{\mathrm{rk}}\nolimits}

\def\ch{\mathop{\mathrm{ch}}\nolimits}

\def\id{\mathop{\mathrm{id}}\nolimits}
\def\Ext{\mathop{\mathrm{Ext}}\nolimits}
\def\exp{\mathop{\mathrm{exp}}\nolimits}

\def\Coh{\mathop{\mathrm{Coh}}\nolimits}
\def\Spec{\mathop{\mathrm{Spec}}\nolimits}

\newcommand{\Res}{\operatornamewithlimits{Res}}
\newcommand{\Dec}{\operatornamewithlimits{Dec}}

\newcommand{\Hilb}{\mathop{\text{-}\mathrm{Hilb}}\nolimits}

\newcommand{\ba}{\Bar}
\newcommand{\wt}{\widetilde}
\newcommand{\wh}{\widehat}
\newcommand{\mk}{\mathfrak}
\newcommand{\mc}{\mathcal}
\newcommand{\mb}{\mathbb}

\newcommand{\SL}{\text{SL}}

\newcommand{\mo}{\mathcal{O}}
\newcommand{\E}{\mathcal{E}}
\newcommand{\F}{\mathcal{F}}
\newcommand{\PP}{\mathbb{P}}

\newcommand{\Z}{\mathbb{Z}}

\newcommand{\GL}{\operatorname{GL}}

\begin{document}

\title{Functional equations of Nekrasov functions proposed by Ito-Maruyoshi-Okuda}

\author{Ryo Ohkawa}

%\institute{
%\at Department of Mathematics, School of Fundamental Science and Engineering, Waseda University, 3--4--1 Okubo, Shinjuku-ku, Tokyo 169--8555, Japan \\
%\email{ohkawa.ryo@aoni.waseda.jp}\\
%Tel.: +81-3-5286-3195
%}

\date{}

\begin{abstract}
We prove functional equations of Nekrasov partition functions for $A_{1}$-singularity, suggested by Ito-Maruyoshi-Okuda \cite{IMO}.
Our proof is given by the computation similar to \cite{O}.
This is the method by Nakajima-Yoshioka \cite{NY3} based on the theory of wall-crossing formula developed by Mochizuki \cite{Mo}.
\end{abstract}

\keywords{Nekrasov partition functions \and Framed sheaves \and Instantons}

\subjclass[2000]{14D21 \and 57R57}

\maketitle

%%%%%%%%%%%%%%%%%%%%%%%%%%%%%%%%%%%%%%%%%%%%%%%%%%%%%%%%%%%%%%%%%%%%%%%%%%%
%%%%%%%%%%%%%%%%%%%%%%%%%%%%%%%%%%%%%%%%%%%%%%%%%%%%%%%%%%%%%%%%%%%%%%%%%%%
%%%%%%%%%%%%%%%%%%%%%%%%%%%%%%%%%%%%%%%%%%%%%%%%%%%%%%%%%%%%%%%%%%%%%%%%%%%

\section{Introduction}
\label{sec:intr}

Nekrasov partition functions are introduced by Nekrasov \cite{Nek}.
They are defined by integrations
$$
Z_{\PP^{2}}(\mbi \e, \mbi a, q)=\sum_{n=0}^{\infty} q^{n} \int_{M(r,n)} 1
$$ 
on moduli spaces $M(r,n)$ of framed sheaves on the plane $\PP^{2}$ with the rank $r$ and the second Chern class $n$, where the integrand $1$ can be replaced by various equivariant cohomology classes on $M(r,n)$ corresponding to physical theories.
These integrations are defined by localization for torus actions on moduli spaces, and variables $\mbi \e=(\e_{1} , \e_{1})$ and $\mbi a=(a_{1}, \ldots, a_{r})$ correspond to the diagonal $T^{2}$-actions on $\C^2 \subset \PP^2$ and $T^{r}$-actions by scale change of framings. 

Nekrasov's conjecture states that these partition functions give deformations of the Seiberg-Witten prepotentials for $N=2$ SUSY Yang-Mills theory.
This conjecture is proven in Braverman-Etingof \cite{BE}, Nekrasov-Okounkov \cite{NO} and Nakajima-Yoshioka \cite{NY1} independently.
In \cite{NY1}, they study relationships with similar partition functions defined for blow-up $\hat{\PP}^{2}$ of $\PP^{2}$ along the origin, and get {\it blow-up formula}. 
This is the formula for bilinear relations of $Z_{\PP^{2}}(\e_{1}, \e_{2}-\e_{1}, \mbi a, q)$ and $Z_{\PP^{2}}(\e_{1}- \e_{1}, \e_{2}, \mbi a, q)$, which correspond to $T^{2}$-fixed points on $(-1)$ curve of $\hat{\PP}^{2}$.
Furthermore these arguments are extended in \cite{NY3} to various cohomology classes other than $1$ using the theory of perverse coherent sheaves.

On the other hand, when $r=2$ and $\mbi a=(a,-a)$, Alday-Gaiotto-Tachikawa \cite{AGT} proposed the relation
$$
Z_{\PP^{2}}(\mbi \e, \mbi a, q)= \left( \frac{q}{\e_{1}^{2} \e_{2}^{2}} \right)^{- \Delta} \mc F_{c} \left( \Delta \Big | \frac{q}{\e_{1}^{2} \e_{2}^{2}} \right)
$$
where $\mc F_{c} \left( \Delta \Big | \frac{q}{\e_{1}^{2} \e_{2}^{2}} \right)$ is the conformal block with central charge $c= 1 + 6 (\e_{1} + \e_{2})^{2} /( \e_{1}\e_{2})^{2}$ and conformal weight $\Delta = (4a^{2} - (\e_{1} + \e_{2})^{2} ) / ( 4\e_{1} \e_{2} )$.
These are simplified $c=1, \Delta= a^{2} /(\e_{1} \e_{2})$ when we assume
\begin{eqnarray*}
\e_{1} + \e_{2} = 0.
\end{eqnarray*}
Gamayun-Iorgov-Lisovyy \cite{GIL} suggested that series of conformal blocks forms $\tau$ functions of Painlev\'e equations.
This conjecture are proved in Bershtein-Shchechkin \cite{BS1} and Iorgov-Lisovyy-Teschner \cite{ILT} by different methods.

In \cite{BS1}, bilinear relations for $Z_{\PP^{2}}(2 \e_{1}, \e_{2}-\e_{1}, \mbi a, q)$ and $Z_{\PP^{2}}(\e_{1}- \e_{2}, 2\e_{2}, \mbi a, q)$ are shown by representation theoretic method.
This is called $(-2)$ {\it blow-up formula} since the variables $(2 \e_{1}, \e_{2}-\e_{1})$ and $(\e_{1}- \e_{2}, 2\e_{2})$ correspond to $T^{2}$-fixed points on $(-2)$ curve.
Here we note that the condition $\e_{1}+ \e_{2}=0$ keeps to hold under these change of variables, hence this formula is applied to study conformal blocks with the central charge $c=1$.
In \cite{BS3}, they also constructed Painlev\'e $\tau$ function from the original Nakajima-Yoshioka blow-up formula where $c=-2$.
So it is plausible to say that study of $(-2)$ blow-up formula from geometrical point of view has interesting application to study of Painlev\'e $\tau$ function.
In particular, $K$-theoretic integration rather than cohomological one directly corresponds to discrete Painlev\'e $\tau$ functions proposed by \cite{BS2}, and conjectural functional equations are given there.  

With those applications in mind, we study similar partition functions for $(-2)$ curve, that is, ALE space of type $A_{1}$.
This space is derived equivalent to the quotient stack $[\C^2/ H]$, where $H = \lbrace \pm \id_{\C^{2}} \rbrace$ naturally acts on $\C^{2}$.
Moduli space of instantons on ALE space is singular, and two resolutions are given by moduli spaces of framed sheaves on the ALE space and the quotient stack $[\C^2/ H]$.
These two resolutions are also considered as moduli spaces of stable ADHM data ( or, quiver varieties defined by Nakajima \cite{N1} in more general setting) corresponding to different stability parameters.
We study relation between two partition functions defined by integrations over these moduli spaces of framed sheaves.
As a main result, we show formulas in Theorem \ref{main} among these two partition functions conjectured by Ito-Maruyoshi-Okuda \cite[(4.1), (4.2)]{IMO}.
We also note that the similar phenomenon to this paper is studied by Belavin-Bershtein-Feigin-Litvinov-Tarnopolsky \cite{BBFLT} from representation theoretic point of view.

Partition functions for the ALE space in this paper can be regarded as degree zero Hirota differential, that is, multiplication for $Z_{\PP^{2}}(2 \e_{1}, \e_{2}-\e_{1}, \mbi a, q)$ and $Z_{\PP^{2}}(\e_{1}- \e_{2}, 2\e_{2}, \mbi a, q)$ by \eqref{mult} in Appendix \ref{sec:comb}. 
To get higher order Hirota differential, we need to multiply integrands with power series $\exp ( \sum_{d=1}^{\infty} \mu(C)^{d} t^{d} )$, where $\mu(C)$ is a slant product of Chern class of universal framed sheaf with the fundamental cycle of $(-2)$ curves $C$ on ALE space.
In this paper, we consider the case where $t=0$.
We also expect that arguments in this paper can be extended to $K$-theoretic integrations in the near future.

In the rest of the paper, we give a proof of our main result, Theorem \ref{main}, as follows.
In \S 2, we give outline of the paper.
In \S 3, we explain the conjecture by Ito-Maruyoshi-Okuda, that is, Theorem \ref{main}.
In \S 4, we recall ADHM description of framed sheaves, and give a precise definition of integrations over moduli spaces of framed sheaves.
In \S 5, we recall Mochizuki method in the manner similar to our previous results \cite{O}.
In \S 6, we compute wall-crossing formulas, and complete a proof of Theorem \ref{main}.
In Appendix A, we construct moduli spaces of framed sheaves in terms of ADHM data.
In Appendix B, we recall combinatorial description of partition functions for the ALE space and the quotient stack $[\C^{2} / H]$ following \cite{NY1}.

The author thanks Hiraku Nakajima for many advices, Isamu Iwanari for telling him patching of stacks and \cite{I}, and Mikhail Bershtein for discussion about bilinear relations of Painlev\'e $\tau$ functions which gives him a motivation and helps him to write introduction in this paper. 
The author is grateful to the referee for her or his careful reading of the paper.
This paper is a part of the outcome of research performed under a Waseda University Grant for Special Research Projects (Project number: 2017S-077).

%%%%%%%%%%%%%%%%%%%%%%%%%%%%%%%%%%%%%%%%%%%%%%%%%%%%%%%%%%%%%%%%%%%%%%%%%%%
%%%%%%%%%%%%%%%%%%%%%%%%%%%%%%%%%%%%%%%%%%%%%%%%%%%%%%%%%%%%%%%%%%%%%%%%%%%
%%%%%%%%%%%%%%%%%%%%%%%%%%%%%%%%%%%%%%%%%%%%%%%%%%%%%%%%%%%%%%%%%%%%%%%%%%%

\section{Outline of the paper}

We do not include self-contained arguments in this paper, since many parts of the argument is similar as in the previous research \cite{NY3} and \cite{O}.
Instead, we give outline here, and in the main body of the paper, we indicate proofs of corresponding statements in loc. cit.
In addition, we try to explain in a unified manner as possible.

We reduce a proof of functional equations in Theorem \ref{main} to wall-crossing formula of moduli of ADHM data, or quiver variety associated to $A_{1}^{(1)}$-quiver.
To analize wall-crossing formula, we introduce enhanced master spaces in terms of ADHM data.
In our argument, vector spaces in \cite{O} are replaced to $\Z_{2}=\Z/2\Z$-graded vector spaces.
This change increase the dimension of spaces of stability parameters, which is similar to the situation in \cite{NY3}, and we have more complicated wall-and-chamber structure in $\R^{2}$ than \cite{O}, where we have only two kinds of generic stability conditions. 

Walls are defined by roots of the Dynkin diagram $A_{1}^{(1)}$, and we follow the arguments in \cite{NY3} for real roots and ones in \cite{O} for imaginary roots.
The difference of these two kinds of walls is mainly description of moduli spaces $M_{\mbi \alpha}^{p}$ of destabilizing objects appearing in iterated cohomology classes \eqref{itcoho} used for wall-crossing formula in Theorem \ref{thm:wcm}.
These moduli spaces $M_{\mbi \alpha}^{p}$ can be embedded into moduli spaces of ADHM data.
For imaginary roots, these embedding become identity, while for real roots these are non-identity and half dimensional.
But we can describe recursive procedure for both cases in a unified manner. 

 More technically, we need the obstruction theory for localization procedure over enhanced master spaces.
But we completely omit this argument, since it is similar as in \cite[\S 6]{O}, and our obstruction theories for various moduli spaces are naturally obtained from constructions.
Here we only mention that obstruction theory for $M_{\mbi \alpha}^{p}$ differs from one for the moduli space of ADHM data containing $M_{\mbi \alpha}^{p}$.
This causes unusual fundamental cycles of $M_{\mbi \alpha}^{p}$ for imaginary roots, while we get usual fundamental cycles for real roots, although $M_{\mbi \alpha}^{p}$ are smooth in both cases.
This phenomenon reflects the fact that enhanced master spaces are singular for imaginary roots, but smooth for real roots.
In this paper, we do not need detailed description of obstruction theory of $M_{\mbi \alpha}^{p}$ for a real root $\mbi \alpha$, since we only use vanishing of wall-crossing terms in this case.
On the other hand, when $\mbi \alpha$ is an imaginary root, they are computed in Proposition \ref{hilb} by using combinatorial descriptions in Appendix \ref{sec:comb} and reduced to the similar computations \cite[Proposition 8.1]{O}.

%%%%%%%%%%%%%%%%%%%%%%%%%%%%%%%%%%%%%%%%%%%%%%%%%%%%%%%%%%%%%%%%%%%%%%%%

%\subsection{Notation}
%\label{subsec:notation}

%%%%%%%%%%%%%%%%%%%%%%%%%%%%%%%%%%%%%%%%%%%%%%%%%%%%%%%%%%%%%%%%%%%%%%%%

\subsection{Moduli of ADHM data}
\label{subsec:modu0}
Let $Q = \C^{2}$ be the affine plane, and consider a finite sub-group $H=\lbrace \pm \id_{Q} \rbrace$ of $\SL(Q)$ with the action on $\PP^{2}=\PP(\C \oplus Q)$ induced by the natural $\SL(Q)$-action.
In \S \ref{subsec:modu1}, we consider {\it framed sheaves} on two resolutions $X_{0}$ and $X_{1}$ of the orbit space $\PP^{2} / H$, and describe them in terms of {\it ADHM data} as follows. 

The character group of $H$ is regarded as $\Z_{2}=\Z/2\Z$, where $0$ is regarded as the trivial character.
Let $W=W_{0} \oplus W_{1}$ and $V=V_{0} \oplus V_{1}$ be $\Z_{2}$-graded vector spaces.
These are regarded as weights decompositions of $H$-representations.
In this sense, we have $Q=Q_{0} \oplus Q_{1}$ with $Q_{0}=0$.

We put
\begin{eqnarray*}
\mb M(W, V)
&=&
\Hom_{\Z_{2}}(Q^{\vee} \otimes V, V) \oplus \Hom_{\Z_{2}}(W,V)\oplus \Hom_{\Z_{2}}(\det Q^{\vee} \otimes V,W),\\
\mb L(V)
&=&
\Hom_{\Z_{2}}(\det Q^{\vee} \otimes V, V).
\end{eqnarray*}
Here subscript $\Z_{2}$ means that these homomorphisms must respect $\Z_{2}$-gradings.
We also write $\mb M=\mb M(\mbi w, \mbi v)=\mb M(W,V)$ and $\mb L=\mb L(\mbi v)=\mb L(V)$ for $\mbi w =  (w_{0}, w_{1}), \mbi v =  (v_{0}, v_{1}) \in \Z^{2}$, where $w_{i} = \dim W_{i}, v_{i}=\dim V_{i}$ for $i=0,1$.
These vectors corresponds to Chern classes of framed sheaves later (cf. Theorem \ref{ale}).
In particular, $W$ corresponds to framings over the infinity line $\ell_{\infty} = [ \PP(Q) / H]$ via the natural identifications of coherent sheaves on $\ell_{\infty}$ with $H$-equivariant sheaves on $\PP^{1}=\PP(Q)$. 

In Definition \ref{defn:adhm}, we consider a map  
\begin{eqnarray*}
\mu \colon \mb M \to \mb L, \mc A= (B, z,w) \mapsto [B \wedge B]+zw,
\end{eqnarray*}
where $B \in \Hom_{\Z_{2}}(Q^{\vee} \otimes V, V)$, $z \in \Hom_{\Z_{2}}(W,V)$, and $w \in \Hom_{\Z_{2}}(\det Q^{\vee} \otimes V,W)$, and elements $(B, z,w)$ in $\mu^{-1}(\mbi 0)$ are called ADHM data on $(W, V)$.
For elements in $\mb M(W, V)$, we introduce stability conditions parametrized by $\mbi \zeta \in \R^{2}$ in Definition \ref{ADHM-stability}.
We put 
\begin{eqnarray*}
M^{\mbi \zeta}(\mbi w, \mbi v)
&=&
\lbrace \mc A= (B, z,w) \mid \text{ $\mbi \zeta$-stable ADHM data on }(W, V)\rbrace / G,
\end{eqnarray*}
where $G=\GL(V_{0}) \times \GL(V_{1}).$
We have natural $\GL(Q) \times \GL(W_{0}) \times \GL(W_{1})$-action on $M^{\mbi \zeta}(\mbi w, \mbi v)$.

%%%%%%%%%%%%%%%%%%%%%%%%%%%%%%%%%%%%%%%%%%%%%%%%%%%%%%%%%%%%%%%%%%%%%%%%

\subsection{Wall-and-chamber structure}
\label{subsec:wall}

As in \cite[1(ii)]{N3}, we introduce $\mbi \alpha_{m} =  (|m|, |m+1|)$ for $m \in \Z$, and $\mbi \delta = (1,1)$.
We call $\pm \mbi \alpha_{m}$ for $m \in \Z$ {\it real roots}, and $p \mbi \delta$ for $p \in \Z \setminus \lbrace \mbi 0 \rbrace$ {\it imaginary roots}.
We put 
$$
R_{+} = \lbrace \mbi \alpha_{m} \mid m \in \Z \rbrace \cup \lbrace p \mbi \delta \mid p >0 \rbrace.
$$
For $\mbi \zeta = (\zeta_{0}, \zeta_{1})$ and $\mbi \alpha=(\alpha_{0}, \alpha_{1})$, we put $(\mbi \zeta, \mbi \alpha)=\zeta_{0} \alpha_{0} + \zeta_{1} \alpha_{1}$, and 
$\mbi \alpha^{\perp} = \lbrace \mbi \zeta \in \R^{2} \mid (\mbi \zeta, \mbi \alpha) = 0 \rbrace.$

For a fixed dimension vector $\mbi v =(v_{0}, v_{1})$, we put $R_{+}(\mbi v) = \lbrace \mbi \alpha = (\alpha_{0}, \alpha_{1}) \in R_{+} \mid \alpha_{0} \le v_{0}, \alpha_{1} \le v_{1} \rbrace$, and call a connected component of 
$$
\R^{2} \setminus \bigcup_{\mbi \alpha \in R_{+}(\mbi v)} \mbi \alpha^{\perp}
$$
a {\it chamber}.

On these chambers, stability and semi-stability coincides, and all $\mbi \zeta$-stability conditions are equivalent when $\mbi \zeta$ lies in one fixed chamber $\mc C$ by \cite[2.8]{N1}.
Hence we can also write $M^{\mc C}(\mbi w, \mbi v) = M^{\mbi \zeta}(\mbi w, \mbi v)$ for $\mbi \zeta \in \mc C$.
On $M^{\mc C}(\mbi w, \mbi v)$, we have {\it tautological bundles} $\mc V_{i}=[\mu^{-1}(0)^{\mbi \zeta} \times V_{i} / G]$ and $\mc W_{i} =[\mu^{-1}(0)^{\mbi \zeta} \times W_{i} / G] \cong W_{i} \otimes \mo_{M^{\mc C}(\mbi w, \mbi v)}$ for $i=0,1$, and  write by $B \colon Q^{\vee} \otimes \mc V \to \mc V$, $z \colon \mc V \to \mc W$ and $w \colon (\det Q ) \otimes \mc W \to \mc V$ corresponding to components written by the same letter in ADHM data, where $\mc V=\mc V_{0} \oplus \mc V_{1}$ and $\mc W= \mc W_{0} \oplus \mc W_{1}$.
These homomorphisms are called {\it tautological homomorphisms}, and they are $\Z_{2}$-graded and $\GL(Q) \times \GL(W_{0}) \times \GL(W_{1})$-equivariant.

For suitable choice $\mbi \zeta^{0}$ and $\mbi \zeta^{1} \in \mb R^{2}$ as in Figure \ref{zeta+} in \S \ref{subsec:adhm3}, we see in Theorem \ref{ale} that these moduli spaces are isomorphic to moduli of framed sheaves on $X_{0}$ and $X_{1}$ respectively.
Our goal, Theorem \ref{main}, is to compare integrations over $M^{\mbi \zeta^{0}}(\mbi w, \mbi v)$ and $M^{\mbi \zeta^{1}}(\mbi w, \mbi v)$.
These parameters $\mbi \zeta^{0}$ and $\mbi \zeta^{1}$ are separated by finitely many walls defined by real roots.
Furthermore, we analyze wall-crossing across $\mbi \delta^{\perp}$ lying between $\mbi \zeta^{1}$ and $-\mbi \zeta^{1}[1]$ ( cf. \S \ref{subsec:symm} and Figure \ref{zeta+} ), and compare $-\mbi \zeta^{1}[1]$ and $\mbi \zeta^{1}[1]$ which leads to change $- \mbi \e$ of variables from $\mbi \e$. 
Then we also have finitely many walls defined by real roots between $\mbi \zeta^{1}[1]$ and $\mbi \zeta^{0}$.
These process give two kinds of functional equations in Theorem \ref{main}.

To analyze one fixed wall-crossing, that is, two chambers $\mc C, \mc C'$ adjacent to a common boundary ray $\mk D$, we divide $\mbi \alpha^{\perp}$ into two rays
\begin{eqnarray*}
\mathfrak{D}_{\mbi \alpha}=
\lbrace \mbi \zeta= (\zeta_{0}, \zeta_{1}) \in \mb R^{2} \mid (\mbi \zeta, \mbi \alpha) = 0, \zeta_{0} \le  \zeta_{1} \rbrace  
\end{eqnarray*}
and $- \mathfrak{D}_{\mbi \alpha} = \lbrace - \mbi \zeta \mid \mbi \zeta \in \mk D_{\mbi \alpha} \rbrace$ such that $\mbi \alpha^{\perp} = \mathfrak{D}_{\mbi \alpha} \cup ( - \mathfrak{D}_{\mbi \alpha})$.
We call subsets $\pm \mk D_{\mbi \alpha}$ for $\mbi \alpha \in R_{+}(\mbi v)$ a {\it wall}.

%%%%%%%%%%%%%%%%%%%%%%%%%%%%%%%%%%%%%%%%%%%%%%%%%%%%%%%%%%%%%%%%%%%%%%%%

\subsection{Enhanced master space}
\label{subsec:enha0}
In \S \ref{sec:Moch}, we introduce enhanced master spaces as follows.
We consider ADHM data in $\mb M(\mbi w, \mbi v)$ and fix $i_{0} \in \Z_{2}$.
We also fix chambers $\mc C, \mc C'$ and a wall $\mk D$ as above, and assume that $\mk D \subset \mbi \alpha^{\perp}$ for $\mbi \alpha \in R_{+}(\mbi v)$, that is, $\mk D= \mk D_{\mbi \alpha}$ or $-\mk D_{\mbi \alpha}$, and $(\mbi \zeta, \mbi \alpha) < 0$ for $\mbi \zeta \in \mc C$.

Roughly, the {\it enhanced master space} $\mc M$ is a moduli space parametrizing $G$-orbits of objects
\begin{eqnarray}
\label{enhanced}
( \mc A, F^{\bullet}, [x_{-}, x_{+}])
\end{eqnarray}
satisfying certain stability condition, where $\mc A=(B, z, w)$ are ADHM data in $\mb M(\mbi w, \mbi v)$, $F^{\bullet}$ are full flag of $V_{i_{0}}$ for fixed $i_{0} \in \Z_{2}$, and $[x_{-}, x_{+}]$ is the homogeneous coordinate of $\PP^{1}$.
For $g=(g_{0}, g_{1}) \in G$, we define $g[x_{-}, x_{+}]= [\det g_{0}^{k_{0}} \det g_{1}^{k_{1}} x_{-}, x_{+}]$ for $(k_{0}, k_{1}) \in \Z^{2}$ defined later in \eqref{kzeta}, and consider natural $G$-actions on the other components.

We consider an algebraic torus $\C^{\ast}_{\hbar}$, and $\C^{\ast}_{\hbar}$-action on $\mc M$ induced by a map $[x_{-}, x_{+}] \mapsto [ e^{\hbar} x_{-}, x_{+}]$ for $e^{\hbar} \in \C^{\ast}_{\hbar}$.
In \S \ref{sec:Moch}, we study $\C^{\ast}_{\hbar}$-fixed point set $\mc M^{\C^{\ast}_{\hbar}}$. 
For each $\ell \in [v_{i_{0}}] = \lbrace 1, \ldots, v_{i_{0}} \rbrace$, if we choose stability conditions suitably, then in \S \ref{subsec:dire}, it will be shown that we get a decomposition 
$$
\mc M^{\C^{\ast}_{\hbar}}= \mc M_{+} \sqcup \mc M_{-} \sqcup \bigsqcup_{\mk J \in \mc D^{\ell}(v_{i_{0}}, \alpha_{i_{0}})} \mc M_{\mk J}.
$$
Here $\mc M_{\pm}$ is an obvious component defined by a zero locus of $x_{\mp}$.
We have an isomorphism $\mc M_{+} \cong \wt{M}^{\mc C, \ell}(\mbi w, \mbi v)$ which parametrizes ADHM data with full flags satisfying $(\mc C, \ell)$-stability, and $\mc M_{-} \cong \wt M^{\mc C, 0}(\mbi w, \mbi v)$.

The other components of $\mc M^{\C^{\ast}_{\hbar}}$ parametrize objects \eqref{enhanced} whose components $(\mc A, F^{\bullet})$ have non-trivial stabilizer groups and satisfying $x_{\pm} \neq 0$.
Consequently, we have decompositions $V=V_{\flat} \oplus V_{\sharp}$ of $\Z_{2}$-graded vector spaces $V_{\flat} = V_{\flat 0} \oplus V_{\flat 1}$ and $V_{\sharp}=V_{\sharp 0} \oplus V_{\sharp 1}$ with $(\dim V_{\sharp 0}, \dim V_{\sharp 1} )= p \mbi \alpha$.
Furthermore \eqref{enhanced} is isomorphic to $(\mc A_{\flat} \oplus \mc A_{\sharp}, F_{\flat}^{\bullet} \oplus F_{\sharp}^{\bullet}, [1, \rho_{\flat} \rho_{\sharp}])$, where 
$$
\mc A_{\flat}=(B_{\flat}, z, w) \in \mb M(\mbi w,  \mbi v - p \mbi \alpha), 
\mc A_{\sharp}=(B_{\sharp}, 0, 0) \in \mb M(\mbi 0, p \mbi \alpha)
$$ 
are ADHM data, and $F_{\flat}^{\bullet}, F_{\sharp}^{\bullet}$ are full flags of $V_{\flat i_{0}}, V_{\sharp i_{0}}$, and $\rho_{\flat}, \rho_{\sharp} \in \C^{\ast}$.
Here we allow the indices of $F_{\flat}^{\bullet}$ and $F_{\sharp}^{\bullet}$ repetitions so that $F_{\flat}^{\bullet} \oplus F_{\sharp}^{\bullet}$ can be regarded as full flags of $V_{i_{0}} = V_{\flat i_{0}} \oplus V_{\sharp i_{0}}$.
Hence we get a disjoint union $I_{\flat} \sqcup I_{\sharp}$ of $[v_{i_{0}}]=\lbrace 1, \ldots, v_{i_{0}}\rbrace$ such that 
$$
I_{\flat} = \lbrace i \in [v_{i_{0}}] \mid F_{\flat}^{i_{0}} / F_{\flat}^{i-1} \neq 0 \rbrace, I_{\sharp} = \lbrace i \in [v_{i_{0}}] \mid F_{\sharp}^{i_{0}} / F_{\sharp}^{i-1} \neq 0 \rbrace.
$$
The index set $\mc D^{\ell}(v_{i_{0}}, \alpha_{i_{0}})$ consists of such {\it decomposition data} $\mk J =(I_{\flat}, I_{\sharp})$ with $I_{\flat} \sqcup I_{\sharp} = [v_{i_{0}}]$, $I_{\sharp} \in \alpha_{i_{0}} \Z_{>0}$, and $\min (I_{\sharp}) \le \ell$.

For each decomposition data $\mk J=(I_{\flat}, I_{\sharp}) \in \mc D^{\ell}(v_{i_{0}}, \alpha_{i_{0}})$, a component $\mc M_{\mk J}$ is \'etale locally written as  
$$
\wt M^{\mc C, \min(I_{\sharp})-1}(\mbi w, \mbi v - p \mbi \alpha) \times \wt M_{\mbi \alpha}^{p},
$$ 
where $p=|I_{\sharp}| / \alpha_{i_{0}}$, and $\wt M_{\mbi \alpha}^{p}$ parametrizes $(\mc A_{\sharp}, F_{\sharp}^{\bullet}, \rho_{\sharp})$, where ADHM data $\mc A_{\sharp}=(B_{\sharp}, 0, 0) \in \mb M(\mbi 0, p \mbi \alpha)$ with full flags $F_{\sharp}^{\bullet}$ of $V_{\sharp i_{0}}$ satisfy $(\mk D, +)$-stability conditions ( cf. Definition \ref{+}).

For such pairs $(\mc A_{\sharp}, F_{\sharp}^{\bullet})$, we take generators $z_{\sharp}$ in $F_{\sharp}^{1}$, and regard them as homomorphisms $z \colon W_{\sharp} \to V_{\sharp}$, where $W_{\sharp} = W_{\sharp 0} \oplus W_{\sharp 1}$ is a one-dimensional $\Z_{2}$-graded vector space with $W_{\sharp i_{0}} = \C$ for fixed $i_{0}$.
Then we get ADHM data $\mc A_{\sharp}^{+} = (B_{\sharp}, z_{\sharp}, 0) \in \mb M(\mbi w_{\sharp}, p \mbi \alpha)$, where $\mbi w_{\sharp} = (w_{\sharp 0}, w_{\sharp 1} )$ and $w_{\sharp 0} = \dim W_{\sharp 0}, w_{\sharp 1} = \dim W_{\sharp 1}$.
We see in Lemma \ref{d+} that $(\mk D, +)$-stability for $(\mc A_{\sharp}, F_{\sharp}^{\bullet})$ is equivalent to $\mbi \zeta$-stability for $\mc A_{\sharp}^{+} = (B_{\sharp}, z_{\sharp}, 0)$, where $\mbi \zeta \in \mc C$.
This motivates us to consider the zero locus $M_{\mbi \alpha}^{p}$ of the tautological homomorphism $w$ on $M^{\mc C}(\mbi w_{\sharp}, p \mbi \alpha)$.
Then we see that $\wt M_{\mbi \alpha}^{p}$ is a full flag bundle of \'etale covers $\hat M_{\mbi \alpha}^{p}$ of $M_{\mbi \alpha}^{p}$ which is obtained by forgetting $\rho_{\sharp}$.

Hence integrations over $\mc M_{\mk J}$ are essentially reduced to integrations over $\wt{M}^{\mc C, \min(I_{\sharp})-1}(\mbi w, \mbi v - p \mbi \alpha) \times M_{\mbi \alpha}^{p}$.
We note that when $\mbi \alpha = \mbi \delta$, we have $M_{\mbi \alpha}^{p} = M^{\mc C}(\mbi w_{\sharp}, p \mbi \alpha)$ since the tautological homomorphism $w$ vanishes by Lemma \ref{d+vszeta} (3).

%%%%%%%%%%%%%%%%%%%%%%%%%%%%%%%%%%%%%%%%%%%%%%%%%%%%%%%%%%%%%%%%%%%%%%%%

\subsection{Recursive procedure}
\label{subsec:recu}
 
By the localization theorem and the above decomposition of $\mc M^{\C^{\ast}_{\hbar}}$, we can calculate difference between integrations over $\mc M_{\pm}$ in terms of $\mc M_{\mk J}$ for $\mk J \in \mc D^{\ell}(v_{i_{0}}, \alpha_{i_{0}})$.

When $\ell = 0$, or $\ell=v_{i_{0}}$, the component $\mc M_{+} = \wt M^{\mc C, \ell}(\mbi w, \mbi v)$ is equal to the full flag bundle of the tautological bundle $\mc V_{i_{0}}$ over $M^{\mc C}(\mbi w, \mbi v)$, or $M^{\mc C'}(\mbi w, \mbi v)$, and integrations over $\mc M_{+}$ are reduced to ones over $M^{\mc C}(\mbi w, \mbi v)$, or $M^{\mc C'}(\mbi w, \mbi v)$ respectively.

Hence if we start from the case where $\ell=v_{i_{0}}$, we get difference between integrations over $M^{\mc C}(\mbi w, \mbi v)$ and $M^{\mc C'}(\mbi w, \mbi v)$ in terms of integrations over $\wt M^{\mc C, \ell'}(\mbi w, \mbi v') \times M_{\mbi \alpha}^{p}$ for $\ell' < v_{i_{0}}$ and $\mbi v' = \mbi v - p \mbi \alpha$ with $p =1, \ldots, \lfloor v_{i_{0}} / \alpha_{i_{0}} \rfloor$.
These can be viewed as integrations over $\wt M^{\mc C, \ell'}(\mbi w, \mbi v')$ whose integrands are defined by integrations along projections $\wt M^{\mc C, \ell'}(\mbi w, \mbi v') \times M_{\mbi \alpha}^{p} \to \wt M^{\mc C, \ell'}(\mbi w, \mbi v')$.
These integrands are also viewed as special cases of iterated cohomology classes \eqref{itcoho} in \S \ref{subsec:iter1}.

Then we continue putting $\mbi v = \mbi v' $ and $\ell = \ell'$, and get integrations over $\wt M^{\mc C, \ell'}(\mbi w, \mbi v')$ for $\ell' < \ell$.
Hence this recursion procedure terminates finally, and get Theorem \ref{thm:wcm}.
Applying this formula repeatedly to walls separating $\mbi \zeta^{0}$ and $\mbi \zeta^{1}$ as explained in \S \ref{subsec:wall}, we get formula in Theorem \ref{main} conjectured by \cite{IMO}.

%%%%%%%%%%%%%%%%%%%%%%%%%%%%%%%%%%%%%%%%%%%%%%%%%%%%%%%%%%%%%%%%%%%%%%%%

%%%%%%%%%%%%%%%%%%%%%%%%%%%%%%%%%%%%%%%%%%%%%%%%%%%%%%%%%%%%%%%%%%%%%%%%%%%
%%%%%%%%%%%%%%%%%%%%%%%%%%%%%%%%%%%%%%%%%%%%%%%%%%%%%%%%%%%%%%%%%%%%%%%%%%%
%%%%%%%%%%%%%%%%%%%%%%%%%%%%%%%%%%%%%%%%%%%%%%%%%%%%%%%%%%%%%%%%%%%%%%%%%%%

\section{Ito-Maruyoshi-Okuda conjecture}
We explain formula conjectured by \cite{IMO}.

\subsection{Two resolutions}
\label{sub:two}
Let $Q$ be an affine plane $\C^{2}$ and $H$ a finite sub-group of $\SL(Q)$, and consider the left $H$-action on $Q$ induced by the natural $\SL(Q)$-action.
We have a unique minimal resolution 
$$
g \colon X_{1}^{\circ} = H\Hilb (Q)\to Q/ H = \Spec \C[x_1, x_2]^{H}
$$ 
of the quotient singularity, and the exceptional curve $C \subset X_{1}^{\circ}$ with the dual graph equal to the Dynkin diagram corresponding to $ADE$ classification of $H$ (cf. \cite[Chapter 4]{N2}, or \cite[Chapter 12]{K}).

On the other hand, we consider the natural $H$-action on $\PP^2=\PP(\C \oplus Q)$, 
%and the homogeneous coordinate $[x_{0}, x_{1}, x_{2}] \in \PP^2$, 
and the quotient stack $X_{0}=[\PP^{2}/H]$, and the coarse moduli map $f \colon X_{0} \to \PP^{2}/H$.  
The only loci with non-trivial stabilizer groups are $O = [\lbrace x_{1}=x_{2}=0 \rbrace/H]$ and $\ell_{\infty}= [\lbrace x_{0}=0 \rbrace / H]$ in $X_{0}$.
Hence we have an isomorphism 
\begin{eqnarray}
\label{isom1}
X_{0}\setminus (O \sqcup \ell_{\infty}) \cong X_{1}^{\circ} \setminus C.
\end{eqnarray}
We patch $X_{1}^{\circ}$ and $X_{0} \setminus O$ to get a compactification $X_{1}$ of $X_{1}^{\circ}$ via this isomorphism:
$$
X_{1}= X_{1}^{\circ} \sqcup \ell_{\infty} = X_{1}^{\circ} \cup ( X_{0}\setminus O)
$$ 

Framed sheaves on $X_{\kappa}$ for $\kappa=0,1$ are pairs $(E, \Phi)$ of sheaves $E$ on $X_{\kappa}$ and isomorphisms $\Phi \colon E|_{\ell_{\infty}} \cong \mo_{\PP^{1}} \otimes {\rho}$ for some $H$-representation $\rho$, where we identify coherent sheavs on $\ell_{\infty}$ as $H$-equivariant coherent sheaves on $\PP^{1}=\PP(Q)$.
Moduli of framed sheaves are examples of quiver varieties constructed from quivers whose underlying graphs are Dynkin diagrams corresponding to $ADE$ classification of $H$ as shown in \cite{N3}.
We define Nekrasov partition functions by integrations on these moduli spaces.

In this paper, we consider the case where $H=\left \lbrace \pm \id_{Q} \right \rbrace, $
and prove relations between Nekrasov partition functions defined from $X_{0}$ and $X_{1}$, suggested by \cite{IMO}.
In this case, the corresponding graph is of $A_{1}^{(1)}$-type.
We write by $F$ the divisor class defined by $\lbrace x_{i}=0 \rbrace$ for $i=1,2$ on $X_{0}$.
By the same symbol $F$, we also write the divisor class on $X_{1}$ defined via the isomorphism \eqref{isom1}.
Similarly we write by $\mc R_{0}$ and $\mc R_{1}$ line bundles $\mo_{X_{\kappa}}$ and $\mo_{X_{\kappa}}(F - \ell_{\infty})$ on $X_{\kappa}$ for $\kappa=0,1$, and put $\mc R= \mc R_{0} \oplus \mc R_{1}$.
When $H=\lbrace \pm \id_{Q} \rbrace$ as above, our compactification $X_{1}$ coincides with the one used in the previous research \cite{BPSS}. 

\begin{figure}[h]
\includegraphics[scale=0.4]{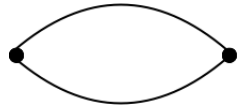}
\caption{graph of type $A_{1}^{(1)}$}
\end{figure}

%%%%%%%%%%%%%%%%%%%%%%%%%%%%%%%%%%%%%%%%%%%%%%%%%%%

\subsection{Derived equivalence}
\label{subsec:deri}
We consider the universal subscheme $U$ over $X_{1}^{\circ}=H\Hilb(Q)$.
\begin{eqnarray}
\label{diagram}
\xymatrix{U \ar[d] \ar[r] & Q\\
X_{1}^{\circ} &
}
\end{eqnarray}
Each fibre of $U \to X_{1}^{\circ}$ is a $H$-invariant subscheme of $Q$, and we have a natural $H$-action on $U$.
We put $W^{\circ} = [U/H]$. 
Then the above $\Gamma$-equivariant diagram \eqref{diagram} induces 
\begin{eqnarray*}
\xymatrix{W^{\circ} \ar[d] \ar[r] & X_{0}^{\circ}=[Q/H]\\
X_{1}^{\circ} &
}
\end{eqnarray*}
where $X_{0}^{\circ}=[Q/H]$ is an open subset $X_{0} \setminus \ell_{\infty}$ of $X_{0}$.
These morphisms are isomorphisms outside exceptional sets.
As in the previous section, we patch together to get
$$
W = W^{\circ} \cup (X_{0} \setminus O) = W^{\circ} \sqcup \ell_{\infty}
$$
and morphisms $p \colon W \to X_{0}, q \colon W \to X_{1}$.
We write by $D(X_{0}), D(X_{1})$ bounded derived categories of coherent sheaves on $X_{0}, X_{1}$.

\begin{prop}
\label{eq}
A functor $F=q_{\ast} p^{\ast} \colon D(X_{0}) \to D(X_{1})$ gives an equivalence of categories.
\end{prop}
\proof
First we show faithfulness of $F$.
By the similar argument to \cite[Example 2.2]{B}, we have a spanning class 
$$
\Omega=\lbrace \mo_{Z} \otimes \rho \mid Z \colon H\text{-orbit in } \PP^{2}, \rho \in \text{Irrep} (H) \rbrace
$$
of $D(X_{0})$, and it is enough to check on $\Omega$.
Here $\text{Irrep} (H)$ is the set of irreducible $H$-representations.
If supports of $L, L' \in \Omega'$ are contained in $X_{0} \setminus \ell_{\infty}$, then the derived McKay correspondence by \cite{KV} between $X_{0}^{\circ}$ and $X_{1}^{\circ}$ implies that $\Hom(L,L') \cong \Hom(F(L), F(L'))$.
Otherwise, we also have faithfulness by an isomorphism $X_{0} \setminus O \cong W \setminus p^{-1}(O) \cong X_{1} \setminus C$.

By \cite{Ni}, Serre functors of $D(X_{0}), D(X_{1})$ are given by canonical bundles $\omega_{X_{0}}, \omega_{X_{1}}$ of $X_{0}, X_{1}$.
Since we have an isomorphism
\begin{eqnarray*}
\label{canonical}
p^{\ast} \omega_{X_{0}} \cong q^{\ast} \omega_{X_{1}},
\end{eqnarray*}
Serre functors commute with $F$.
Hence, $F$ is an equivalence by \cite[Theorem 2.3]{BKR}. 
\endproof

We also write by $F \colon K(X_{0}) \to K(X_{1})$ the induced isomorphism of $K$-groups $K(X_{\kappa})$.

%%%%%%%%%%%%%%%%%%%%%%%%%%%%%%%%%%%%%%%%%%%%%%%%%%

\subsection{Chow ring of inertia stacks}
\label{subsec:chow}
We consider inertia stacks $IX_{\kappa} \to X_{\kappa}$, which parametrize pairs $(x, \sigma)$ of objects $x$ of $X_{\kappa}$ and automorphisms $\sigma$ of $x$.
We identify $X_{\kappa}$ as components of $IX_{\kappa}$ consisting of objects of $X_{\kappa}$ together with identities.
We have other components $\ell_{\infty}^{1}$ of $IX_{0}$ and $IX_{1}$, and $O^{1}$ of $IX_{0}$, consisting of objects of $\ell_{\infty}$ and $O$ together with non-trivial automorphisms. 
These are isomorphic to $\ell_{\infty}$ and $O$ as stacks.
We have $IX_{0}=X_{0} \sqcup \ell_{\infty}^{1} \sqcup O^{1}$ and $IX_{1}= X_{1} \sqcup \ell_{\infty}^{1}$.

By \cite[Theorem 2.2]{I}, Chow rings of $IX_{\kappa}$ are described as
\begin{eqnarray*}
A(IX_{0})
&=&
A(X_{0}) \oplus A(\ell_{\infty}^{1}) \oplus A(O^{1})\\
&=&\left( \C [X_{0}] \oplus \C [\ell_{\infty}] \oplus \C [O] \right) \oplus \left( \C [\ell_{\infty}^{1}] \oplus \C [Q^{1}] \right) \oplus \C [O^{1}],\\
A(IX_{1})
&=& A(X_{1}) \oplus A(\ell_{\infty}^{1})\\
&=&\left( \C[X_{1}] \oplus \C[\ell_{\infty}] \oplus \C [F] \oplus \C [P] \right) \oplus \left( \C[\ell_{\infty}^{1}] \oplus \C [Q^{1}] \right),
\end{eqnarray*}
where $P=F \cap C$ in $X_{1}$, and $Q=F \cap \ell_{\infty}$ in $X_{\kappa}$ for $\kappa=0,1$.
We have $Q= \frac{1}{2} P$ in $A(X_{1})$, and $Q^{1}$ is a sub-stack of $\ell_{\infty}^{1}$ consisting of a point in $Q$ with the non-trivial stabilizer group.

For $\alpha \in K(X_{\kappa})$, we define $\wt{\ch}(\alpha) \in A(IX_{\kappa})$ as follows.
Any vector bundle $E$ on $\ell_{\infty}^{1}$ has an eigen decomposition $E_{0} \oplus E_{1}$ for the action of non-trivial automorphisms.
We define $\rho(E)=E_{0}-E_{1}$, which gives an operation on $K(\ell_{\infty}^{1})$.
Similarly we define an operator $\rho$ on $K(O^{1})$.
On $K(X_{\kappa})$, we define $\rho=\id_{K(X_{\kappa})}$. 
Then we have operators $\rho$ on $K(IX_{\kappa})$ for $\kappa=0,1$.
We define $\wt{\ch}(\alpha)$ by $\ch (\rho(\alpha |_{IX_{\kappa}}))$, where $\alpha |_{IX_{\kappa}}$ is the pull-back of $\alpha \in K(X_{\kappa})$ to $K(IX_{\kappa})$.

We write by $r$ and $\ba r$ the coefficients of $[X_{\kappa}]$ and $[\ell_{\infty}^{1}]$ in $\wt{\ch}(\alpha)$ respectively.
We put 
$$
r_{0}=r_{0}(\alpha)=(r + \ba r)/2, r_{1}=r_{1}(\alpha)=(r - \ba r)/2,
$$ 
and $\mbi r =\mbi r(\alpha)= (r_{0}, r_{1})$.
For $\kappa=1$, we write by $k(\alpha)$ and $n(\alpha)$ the coefficient of $[C]=2([\ell_{\infty}] - [F])$ and $-[P]$ in $\wt{\ch}(\alpha)$ respectively. 
For $\kappa=0$, we put $k(\alpha)=k(F(\alpha))$ and $n(\alpha)=n(F(\alpha))$, where $F \colon K(X_{0}) \to K(X_{1})$ is an induced isomorphism from the derived equivalence $F \colon D(X_{0}) \to D(X_{1})$ in the previous subsection. 

%%%%%%%%%%%%%%%%%%%%%%%%%%%%%%%%%%%%%%%%%%%%%%%%%%%%%%%%%%%%%%%%%%%%%%%%%%%

\subsection{Moduli of framed sheaves}
\label{subsec:modu1}
A framed sheaf on $X_{\kappa}$ is a pair $(E, \Phi)$ of a torsion free sheaf $E$ on $X_{\kappa}$, and an isomorphism $\Phi \colon E|_{\ell_{\infty}} \cong \mc O_{X_{\kappa}} \otimes \rho$ called framing, where $\rho$ is a representation of $H=\lbrace \pm \id_{Q} \rbrace$. 
Such a representation $\rho$ is given by a $\Z_{2}$-graded vector space $W= W_{0} \oplus W_{1}$, where $W_{0}$ is trivial, and $W_{1}$ is a sum of non-trivial representations.
If we have $\mbi r(E)=(r_{0}, r_{1})$, then we have $\dim W_{0}=r_{0}, \dim W_{1}=r_{1}$ by definition.

For $c \in A(IX_{\kappa})$, we write by $M_{X_{\kappa}}(c)$ the moduli space of framed sheaves $(E, \Phi)$ on $X_{\kappa}$ with $\wt{\ch} (E)=c$ in $A(IX_{\kappa})$ for $\kappa =0,1$.
This moduli space is constructed, and shown to be smooth and have a universal sheaf $\E$ on $X_{\kappa} \times M_{X_{\kappa}}(c)$ in \cite{BPSS} at least for $\kappa=1$.
For $\kappa_0$, see \cite[Remark 2.2 ]{N2}.
We also construct moduli in terms of ADHM data in \S \ref{subsec:adhm2} and \S Appendix \ref{sec:const} for both $\kappa=0$ and $1$.
%
%\begin{prop}
%Moduli spaces $M_{X_{0}}(c)$ and $M_{X_{1}}(c)$ are smooth quasi-projective varieties.
%\end{prop}
Then smoothness of moduli spaces $M_{X_{\kappa}}(c)$ follows from the following lemma (cf. \cite[4.2]{BPSS}).
\begin{lem}
\label{smooth}
For $(E,\Phi)\in M_{X_{\kappa}}(c)$, we have
$$
\Hom(E,E(-\ell_{\infty}))=\Ext^2(E,E(-\ell_{\infty}))=0.
$$
\end{lem}
\proof
It is similarly proven as in \cite[Proposition 2.1]{NY1}.
\endproof

We introduce {\it tautological bundles} $\mc V_{0}= \mb R^{1} p_{\ast} \mc E (- \ell_{\infty}), \mc V_{1}=\mb R^{1} p_{\ast} \mc E (-F)$, where $p \colon X_{\kappa} \times M_{X_{\kappa}}(c) \to M_{X_{\kappa}}(c)$ is the projection.

%%%%%%%%%%%%%%%%%%%%%%%%%%%%%%%%%%%%%%%%%%%%%%%%%%

\subsection{Torus action on $M_{X_{\kappa}}(c)$}
\label{subsec:toru}
We put $\tilde{T}=T^2 \times T^r \times T^{2r}$, where $T=\C^{\ast}$ is the algebraic torus, and define $\tilde{T}$-action on moduli spaces as follows.
To do that, we consider $X_{0}$ and $X_{1}$ as quotient stacks
\begin{eqnarray}
\label{cpt1}
X_{0}&=&[\PP^{2} / \lbrace \pm 1\rbrace]\\
\label{cpt2}
X_{1}&=&[\lbrace (y, x_0, x_1, x_2) \in \C^{4} \mid (y, x_{0}) \neq 0, (x_1, x_2) \neq 0 \rbrace / \left( \C^{\ast}_{s} \times \C^{\ast}_{t} \right)],
\end{eqnarray}
where the $\C^{\ast}_{s} \times \C^{\ast}_{t}$-action is defined by 
$$
(s,t)(y, x_{0}, x_1, x_2) = \left(\frac{s^2}{t^2} y, s x_0, t x_1, t x_2 \right).
$$
Then we have $T^2$-action on $X_{0}, X_{1}$ defined by
\begin{eqnarray*}
F_{\mbi t} &\colon& X_{0} \to X_{0}, [x_0, x_1, x_2] \mapsto [x_0, t_1 x_1, t_2 x_2]\\
F_{\mbi t} &\colon& X_{1} \to X_{1}, [y, x_0, x_1, x_2] \mapsto [y, x_0, t_1 x_1, t_2 x_2]
\end{eqnarray*}
for $\mbi t=(t_1, t_2) \in T^2$.

We define $\tilde{T}$-action on moduli $M_{X_{\kappa}}(c)$ of framed torsion free sheaves $(E, \Phi)$ by
$$
(E, \Phi) \mapsto \left((F_{\mbi t}^{-1})^\ast E, \Phi'\right)
$$
for $(\mbi t, e^{\mbi a}, e^{\mbi m}) \in \tilde{T}$.
Here $\Phi' \colon (F_{\mbi t}^{-1})^\ast E |_{\ell_{\infty}} \to \mo_{\PP^1} \otimes \rho$ is the composition of the pull-back 
$$
(F_{\mbi t}^{-1})^\ast \Phi \colon (F_{\mbi t}^{-1})^\ast E \to (F_{\mbi t}^{-1})^\ast \left(\mo_{\PP^1}\otimes \rho \right) \cong \mo_{\PP^{1}} \otimes \rho,
$$ 
%$(F_{t})^\ast \Phi$, the natural isomorphism $(F_{t}^{-1})^\ast \left(\mo_{\ell_{\infty}} \otimes \rho \right) \cong \mo_{\PP^1} \otimes \rho $, 
and the diagonal action of $e^{\mbi a}=(e^{a_1},\ldots, e^{a_r})\in T^{r}$
$$
\id_{\mo_{\PP^1}} \otimes \text{diag}(e^{a_1},\ldots, e^{a_r}) \colon 
\mo_{\PP^1} \otimes \rho \to \mo_{\PP^1} \otimes \rho.
$$
Finally $e^{\mbi m} =(e^{m_{1}}, \ldots, e^{m_{2r}}) \in T^{2r}$ trivially acts on moduli spaces, but in the next subsection, we consider fiber-wise action of $T^{2r}$ on vector bundles on moduli spaces.

%%%%%%%%%%%%%%%%%%%%%%%%%%%%%%%%%%%%%%%%%%%%%%%%%%%%%%%%%%%%%%%%%%%%%%%%%%%
\subsection{Partition functions}
\label{subsec:part}
The $\tilde{T}$-equivariant Chow ring $A^{\ast}_{\tilde{T}}(M_{X_{\kappa}}(\alpha))$ is a module over the $\tilde{T}$-equivariant Chow ring $A_{\tilde{T}}^{\ast}(pt)$ of a point, which is isomorphic to $\mb Z[\mbi \e, \mbi a, \mbi m]$, where $\mbi \e=(\e_{1},\e_{2}), \mbi a=(a_{1}, \ldots, a_{r})$, and $\mbi m=(m_{1}, \ldots, m_{2r})$ correspond to the first Chern classes of characters of $\tilde T$ with eigen-values $\mbi t \in T^{2}, e^{\mbi a} \in T^{r}$, and $e^{\mbi m} \in T^{2r}$ respectively.
We write by $\mc S$ the quotient field $\mb Q(\mbi \e, \mbi a, \mbi m)$ of $\Z [ \mbi \e, \mbi a, \mbi m ]$.

We consider a $\tilde{T}$-equivariant vector bundle
$$
\mc F_{r}(\mc V_{0})=\left( \mc V_{0} \otimes \frac{e^{m_{1}}}{\sqrt{t_{1}t_{2}}} \right) \oplus \cdots \oplus \left( \mc V_{0} \otimes \frac{e^{m_{2r}}}{\sqrt{t_{1}t_{2}}} \right)
$$
on $M_{X_{\kappa}}(\alpha)$, and the $\tilde{T}$-equivariant Euler class $e(\mc F_{r}(\mc V_{0}))$, where $e^{\mbi m}=(e^{m_{1}}, \ldots, e^{m_{2r}})$ is an element in the last component $T^{2r}$ of $\tilde{T}$.
Here we consider a homomorphism $\tilde{T}' =\tilde{T} \to \tilde{T}$ defined by
$$
(t_{1}', t_{2}', e^{\mbi a'}, e^{\mbi m'}) \mapsto (t_1, t_{2}, e^{\mbi a}, e^{\mbi m})=((t_{1}')^2, (t_{2}')^2, e^{\mbi a'}, e^{\mbi m'}), 
$$
and use identification $A^{\ast}_{\tilde{T}'}(\text{pt}) \otimes \mc S \cong \mc S$ via $t_{1}'=\sqrt{t_{1}}, t_{2}'=\sqrt{t_{2}}$.

For fixed $\mbi r \in \Z^{2}_{\ge 0}\setminus \lbrace(0,0)\rbrace$ and $k \in \frac{1}{2} \Z$, we define partition functions for $\kappa=0,1$ by
\begin{eqnarray*}
Z_{X_{\kappa}}^{k}(\mbi \e, \mbi a, \mbi m, q)&=&\sum_{\substack{\alpha \in K(X_{\kappa})\\ \mbi r(\alpha)=\mbi r, k(\alpha)=k}} q^{n(\alpha)} \int_{M_{X_\kappa}( \wt{\ch} (\alpha))} e(\mc F_{r}(\mc V_{0})) \in \mc S[[q^{\frac{1}{4}}]].
\end{eqnarray*}
Precise definitions of integrations are explained in \S \ref{subsec:inte}.

The purpose of this paper is to prove the following statement conjectured by \cite{IMO}.
\begin{thm}
\label{main}
We have 
$$
Z_{X_{1}}^{k}(-\mbi \e, \mbi{a},  \mbi m, q) = 
\begin{cases}
(1 - (-1)^{r} q)^{u_{r}} Z_{X_{0}}^{k}(\mbi \e, \mbi a,  \mbi m, q) & \text{ for }k \ge 0,\\
Z_{X_{0}}^{k}(-\mbi \e, \mbi{a}, \mbi{m}, q) & \text{ for }k \le 0,
\end{cases}
$$
where 
$$
u_{r}=\frac{(\e_{1}+\e_{2})(2 \sum_{\alpha=1}^{r} a_{\alpha} + \sum_{f=1}^{2r} m_{f})}{2 \e_{1} \e_{2}}.
$$
\end{thm}

%%%%%%%%%%%%%%%%%%%%%%%%%%%%%%%%%%%%%%%%%%%%%%%%%%%%%%%%%%%%%%%%%%%%%%%%%%%
%%%%%%%%%%%%%%%%%%%%%%%%%%%%%%%%%%%%%%%%%%%%%%%%%%%%%%%%%%%%%%%%%%%%%%%%%%%
%%%%%%%%%%%%%%%%%%%%%%%%%%%%%%%%%%%%%%%%%%%%%%%%%%%%%%%%%%%%%%%%%%%%%%%%%%%
\section{ADHM description}
\label{sec:ADHM}
We introduce ADHM description of framed moduli spaces in the previous section.

%%%%%%%%%%%%%%%%%%%%%%%%%%%%%%%%%%%%%%%%%%%%%%%%%%%%%%%%%%%%%%%%%%%%%%%%%%%
\subsection{ADHM data}
\label{subsec:adhm1}
Let $W=W_{0} \oplus W_{1}, V=V_{0}\oplus V_{1}$ and $Q= Q_{0} \oplus Q_{1}$ be $\Z_{2}$-graded vector spaces with $Q_{0}=0$ and $Q_{1}=\C^2$.
We introduce ADHM data on $(W, V)$.

\begin{defn}
\label{defn:adhm}
ADHM data on $(W, V)$ are collections $(B,z,w)$ of $\Z_{2}$-graded linear maps such that $B \in \Hom_{\Z_{2}} (Q^{\vee} \otimes V, V), z \in \Hom_{\Z_{2}} (W, V)$ and $w \in \Hom_{\Z_{2}} (\det Q ^{\vee} \otimes V, W)$, satisfying 
$$
[B \wedge B]  + zw = 0 \in \Hom_{\Z_{2}} (\det Q ^{\vee} \otimes V, W), 
$$
where $[B \wedge B] $ is the restriction of $B \circ ( \id_{Q^{\vee}} \otimes B ) \colon Q^{\vee} \otimes Q^{\vee} \otimes V \to V$ to the subspace $\det Q ^{\vee} \otimes V$.
%$$
%\xymatrix{
%V \otimes \det Q ^{\vee} \ar[d]^{j} \ar[r]^{\wedge B} & V \otimes Q^{\vee} \ar[d]^{B} \\
%W \ar[r]^{i} & V.
%}
%$$
\end{defn}
If we write $Q^{\vee} = \C e_{1} \oplus \C e_{2}$ and $B( e_{1} \otimes v_{1} + e_{2} \otimes v_{2} )= B_{1} (v_{1}) + B_{2}(v_{2})$, then the equation $[B \wedge B ]+zw=0 $ is equivalent to $[B_1, B_2] + zw =0$. 
We take $\mbi \zeta=(\zeta_{0}, \zeta_{1}) \in \R^{2}$, and put $\mbi \zeta(V)= \zeta_{0} \dim V_{0} + \zeta_{1} \dim V_{1}$ for a $\Z_{2}$-graded vector space $V$.

\begin{defn}
\label{ADHM-stability}
We assume $W \neq 0$.
Elements $(B, z, w)$ in $\mb M(W,V)$ are said to be $\mbi \zeta$-semistable if the conditions (i) (ii) hold for any $\Z_{2}$-graded subspace $S$ of $V$ with $B(Q^{\vee} \otimes S) \subset S$.
\begin{enumerate}
\item If $S \subset \ker w$, then we have $\mbi \zeta(S) \le 0$. 
\item If $\im z \subset S$, then we have $\mbi \zeta(V/S) \ge 0$. 
\end{enumerate}
They are said to be $\mbi \zeta$-stable when the strict inequality always holds for $S \neq 0$ in (i), and $S \neq V$ in (ii).
\end{defn}

We put 
$$
\mbi w = ( \dim W_{0}, \dim W_{1} ), \mbi v = ( \dim V_{0}, \dim V_{1} ),
$$
and construct moduli $M^{\mbi \zeta}(\mbi w, \mbi v)$ of $\mbi \zeta$-semistable ADHM data on $(W, V)$ as follows.
We put 
\begin{eqnarray*}
\mb M &=& \mb M(W, V)=\Hom_{\Z_{2}}(Q^{\vee} \otimes V, V) \times \Hom_{\Z_{2}}(W, V) \times \Hom_{\Z_{2}}(\det Q ^{\vee} \otimes V , W),\\ 
\mb L &=& \mb L(W,V)=\Hom_{\Z_{2}}(\det Q ^{\vee} \otimes V, V),
\end{eqnarray*}
and define a map $\mu \colon \mb M \to \mb L$ by 
\begin{eqnarray}
\label{mu}
\mu(B, z, w) = [B \wedge B] + zw.
\end{eqnarray}
We take the $\mbi \zeta$-semistable locus $\mu^{-1}(0)^{\mbi \zeta}$ and define $M^{\mbi \zeta}(\mbi w, \mbi v)=[\mu^{-1}(0)^{\mbi \zeta} / G ]$, where $G=\GL(V_{0}) \times \GL(V_{1})$.

We have a natural $\GL(Q) \times \GL(W_{0}) \times \GL(W_{1})$-action on $M^{\mbi \zeta}(\mbi w, \mbi v)$. 
Hence via the diagonal embedding $T^{2} \times T^{r}$ into $\GL(Q) \times \GL(W_{0}) \times \GL(W_{1})$ and the projection $\tilde{T} = T^{2} \times T^{r} \times T^{2r}\to T^{2} \times T^{r}$, we get a $\tilde{T}$-action on $M^{\mbi \zeta}(\mbi w, \mbi v)$, where $r=\dim W_{0} +\dim W_{1}$.
Concretely $(\mbi t, e^{\mbi a}) \in \tilde{T}$ acts by $(t_{1} B_{1}, t_{2} B_{2}, z e^{-\mbi a}, e^{\mbi a} w t_{1} t_{2})$.

%%%%%%%%%%%%%%%%%%%%%%%%%%%%%%%%%%%%%%%%%%%%%%%%%%%%%%%%%%%%%%%%%%%%%%%%%%%
\subsection{Wall-and-chamber structure on $\zeta$-plane}
\label{subsec:cham}
As in \S \ref{subsec:wall}, we consider positive roots 
$$
R_{+}= \lbrace \mbi \alpha_{m}=(m, m+1), \mbi \alpha_{-m-1}=(m+1, m) \mid m \in \Z_{\ge 0}  \rbrace \cup \lbrace p \mbi \delta = (p,p)\mid p \in \Z_{>0} \rbrace,
$$
where $\pm \mbi \alpha_{m}$ for $m \in \Z$ are called {\it real roots}, and $p \mbi \delta$ for $p \in \Z \setminus \lbrace 0 \rbrace$ are called {\it imaginary roots}.
These roots $\mbi \alpha \in R^{+}$ define the hyperplane $\mbi \alpha^{\perp}= \lbrace \mbi \zeta \in \mb R^{2} \mid (\mbi \zeta, \mbi \alpha)=0 \rbrace$, where $( \mbi \zeta, \mbi \alpha)=\zeta_{0} \alpha_{0} + \zeta_{1} \alpha_{1}$ for $\mbi \alpha=(\alpha_{0}, \alpha_{1})$ and $\zeta= (\zeta_{0}, \zeta_{1})$.
For later purpose, we divide $\mbi \alpha^{\perp}$ into
\begin{eqnarray*}
\mathfrak{D}_{\mbi \alpha}=
\lbrace \mbi \zeta= (\zeta_{0}, \zeta_{1}) \in \mb R^{2} \mid (\mbi \zeta, \mbi \alpha) = 0, \zeta_{0} \le  \zeta_{1} \rbrace  
\end{eqnarray*}
and $- \mathfrak{D}_{\mbi \alpha}$ such that $\mbi \alpha^{\perp} = \mathfrak{D}_{\mbi \alpha} \cup ( - \mathfrak{D}_{\mbi \alpha})$.

Fix a dimension vector $\mbi v =(v_{0}, v_{1})$, and put $R_{+}(\mbi v) = \lbrace \mbi \alpha \in R_{+} \mid  \alpha_{0} \le v_{0}, \alpha_{1} \le v_{1} \rbrace$.
Subsets $\pm \mk D_{\mbi \alpha}$ for $\mbi \alpha \in R_{+}(\mbi v)$ are called {\it walls}, and {\it chambers} are connected components of 
$$
\R^{2} \setminus \bigcup_{\mbi \alpha \in R_{+}(\mbi v)} \mbi \alpha^{\perp}. 
$$

On these chambers, stability and semi-stability for ADHM data coincide, and all $\mbi \zeta$-stability conditions for ADHM data are equivalent when $\mbi \zeta$ lies in one fixed chamber $\mc C$ by \cite[2.8]{N1}.
Hence we can also write $M^{\mc C}(\mbi w, \mbi v) = M^{\mbi \zeta}(\mbi w, \mbi v)$ for $\mbi \zeta \in \mc C$.
On $M^{\mc C}(\mbi w, \mbi v)$, we have {\it tautological bundles} $\mc V_{i}=[\mu^{-1}(0)^{\mbi \zeta} \times V_{i} / G]$ for $i=0,1$, and {\it tautological homomorphisms} $B \colon Q^{\vee} \otimes \mc V \to \mc V$, where $\mc V=\mc V_{0} \oplus \mc V_{1}$.

For $m \in \Z$, we put
\begin{eqnarray*}
\mc C_{m}&=&\lbrace \zeta=(\zeta_{0}, \zeta_{1}) \in \mb R^{2} \mid m\zeta_{0} + (m+1) \zeta_{1} <0, (m-1) \zeta_{0} + m \zeta_{1}> 0 \rbrace.\\
\end{eqnarray*}

In the following, we take $\mbi \zeta^{0} \in \mc C_{0}$ and $\mbi \zeta^{1} \in \mc C_{m}$ with enough large $m \gg 0$.
By the above chamber structure,  we have chambers $\mc C_{\infty}$ containing $\bigcup_{m >\min \lbrace v_{0}, v_{1} -1 \rbrace} \mc C_{m}$, and $\mc C_{-\infty}$ containing $\bigcup_{m > \min \lbrace v_{0}-1, v_{1} \rbrace} \mc C_{-m}$.
In this notation, we take $\mbi \zeta^{1} \in \mc C_{\infty}$.
%We put $M^{\pm}(\mbi w, \mbi v) = M^{\mc C_{\pm}}(\mbi w, \mbi v)$. 

%%%%%%%%%%%%%%%%%%%%%%%%%%%%%%%%%%%%%%%%%%%%%%%%%%%%%%%%%%%%%%%%%%%%%%%%%%%
\subsection{ADHM description of framed moduli}
\label{subsec:adhm2}
We recall ADHM description from \cite[Chapter 2]{N2} and \cite[Theorem 2.2]{N3}.
\begin{thm}
\label{ale}
We have the following $\tilde{T}$-equivariant isomorphisms.
\begin{enumerate}
\item We have 
$M^{\mbi \zeta^{0}}(\mbi w, \mbi v) \cong M_{X_{0}}(c)$ for $\mbi \zeta^{0} \in \mc C_{0}$, where we put 
\begin{eqnarray}
\label{chern1}
c = \wt \ch \left(w_{0} \mc R_{0} + w_{1} \mc R_{1} - v_{0} \mo_{P} - v_{1} \mo_{P} \otimes (-1) \right) \in A(IX_{0}).
\end{eqnarray}
\item We have  
$M^{\mbi \zeta^{1}}(\mbi w, \mbi v) \cong M_{X_{1}}(c)$ for $\mbi \zeta^{1} \in \mc C_{m}$ with enough large $m \gg 0$, where we put 
\begin{eqnarray}
\label{chern2}
c=(w_{0} + w_{1})[X_{1}]+\left(-2v_{0} + 2v_{1} - w_{1}\right) \frac{[C]}{2} - \left( v_{0} + \frac{w_{1}}{4}\right) [P]  + (w_{0} - w_{1}) [\ell_{\infty}^{1}] \in A(IX_{1}).
\end{eqnarray}
We note that $[C] = 2 [\ell_{\infty}] - 2 [F]$.
\end{enumerate}
Furthermore via these isomorphisms (i) and (ii), tautological bundles $\mc V_{0}, \mc V_{1}$ on both sides coinside as $\tilde{T}$-equivariant vector bundles.
\end{thm}

\proof
We recall these isomorphisms as follows.

We put $\mc R_{0}=\mo_{X_{\kappa}}, \mc R_{1}=\mo_{X_{\kappa}}( F - \ell_{\infty})$ and $\mc R= \mc R_{0} \oplus \mc R_{1}$ as in the introduction.
For ADHM data $(B,z,w)$, we consider the following complex
\begin{equation}
\label{coh}
0\to 
\begin{matrix}
  \mc R_{0}(-\ell_{\infty})\otimes V_0
\\
  \oplus  
\\
  \mc R_{1}(-\ell_{\infty})\otimes V_1
\end{matrix}
\stackrel{\sigma}{\to} 
\begin{matrix}
  (\mc R_{1} \otimes V_{0} \oplus \mc R_{0} \otimes V_{1})^{\oplus 2}
\\
  \oplus  
\\
  \mc R_{0} \otimes W_0 \oplus \mc R_{1} \otimes W_1
\end{matrix}
\stackrel{\tau}{\to}
\begin{matrix}
  \mc R_{0}(\ell_{\infty}) \otimes V_0
\\
  \oplus  
\\
  \mc R_{1}(\ell_{\infty}) \otimes V_1 
\end{matrix}
\to 0
\end{equation} 
with
$$
\sigma=
\begin{bmatrix}
x_{0} B_{1} - x_{1}I_{y}\\
x_{0} B_{2} - x_{2}I_{y}\\
x_{0} w
\end{bmatrix},
\tau=
\begin{bmatrix}
x_{0} B_{2}-x_{2}I^{y}&x_{0} B_{1} - x_{1}I^{y}&x_{0} z
\end{bmatrix},
$$
where $I^{y}=
\begin{bmatrix}
y^{\kappa}&0\\
0&1
\end{bmatrix},
I_{y}=
\begin{bmatrix}
1&0\\
0&y^{\kappa}
\end{bmatrix}
$
for $\kappa=0,1$.
We take its cohomology $E=\ker \tau/\im \sigma$.
By restricting to $\ell_{\infty}$, we get a framing $\Phi$.

In the following, we will show that this map $(B,z,w) \mapsto (E, \Phi)$ gives the desired isomorphism.
For (ii), this follows from \cite[Theorem 2.2]{N3} and Appendix \ref{sec:const}. 

For (i), this follows from \cite[Chapter 2]{N2} as follows.
For a stability parameter $\mbi \zeta^{0}=(\zeta_{0}^{0}, \zeta_{1}^{0}) \in \mc C_{0}$ such that $\zeta_{0}^{0}, \zeta_{1}^{0}<0$, the $\mbi \zeta^{0}$-stability condition in Definition \ref{ADHM-stability} is equivalent to the condition that for any graded subspace $S=S_0\oplus S_1 \subset V$ such that $B$-invariant and $\im z\subset S_{k}$ for $k=0,1$, we have $S=V$.
This is equivalent to the condition that for any subspace $S' \subset V$ (without grading) such that $B$-invariant and $\im z\subset S'$, we have $S'=V$, since we can get a graded subspace $S=S' \cap V_{0} \oplus S' \cap V_{1}$ containing  homogeneous subspace $\im z$.
By \cite[Lemma 2.6]{N3}, this is equivalent to the condition that $\sigma$ is injective except finitely many points and $\tau$ is surjective for any point of $X_{0}$.
This implies that the middle cohomology is a torsion free sheaf on $X_{0}$, and we get the desired isomorphism.

The last assertion follows if we compute $\mc V_{0}= \mb R^{1} p_{\ast} \mc E (- \ell_{\infty}), \mc V_{1}=\mb R^{1} p_{\ast} \mc E (-F)$ on $M_{X_{\kappa}}(c)$ using \eqref{coh}.
\endproof

%%%%%%%%%%%%%%%%%%%%%%%%%%%%%%%%%%%%%%%%%%%%%%%%%%%%%%%%

\subsection{Symmetry}
\label{subsec:symm}
We consider isomorphisms among moduli of semistable ADHM data with various parameters $\mbi \zeta, \mbi w$ and $\mbi v$.
We consider $\Z_{2}$-graded vector spaces $V[1]$ and $W[1]$ and put $\mbi \zeta[1]=(\zeta_{1}, \zeta_{0}),$ 
$$
\mbi v[1]=
(\dim V_{1},
\dim V_{0}
),
\mbi w[1]=
(
\dim W_{1},
\dim W_{0}
).$$
Then any $\mbi \zeta$-semistable ADHM datum on $(W,V)$ is naturally identified with $\mbi \zeta[1]$-semistable ADHM datum on $(W[1], V[1])$, hence we have an isomorphism 
\begin{equation}
\label{shift}
[1]\colon M^{\mbi \zeta}(\mbi w, \mbi v)\cong M^{\mbi \zeta[1]}(\mbi w[1], \mbi v[1]).
\end{equation}

On the other hand, if we take dual vector spaces $W^{\vee}, V^{\vee}$, then $(B_{2}^{\vee}, B_{1}^{\vee}, w^{\vee}, z^{\vee})$ is a $(-\mbi \zeta)$-semistable ADHM datum on $(W^{\vee}, V^{\vee})$ for any $\mbi \zeta$-smistable ADHM datum $(B,z,w)$ on $(W, V)$.
This gives an isomorphism 
\begin{equation}
\label{dual}
M^{\mbi \zeta}(\mbi w, \mbi v) \cong M^{-\mbi \zeta}(\mbi w, \mbi v), 
\end{equation}
which is $\tilde{T}$-equivariant via $\tilde{T} \to  \tilde{T}, ((t_1,t_{2}), e^{\mbi a}, e^{\mbi m}) \mapsto ((t_2,t_1), e^{-a}, e^{\mbi m})$.
%We can also consider $(B_{2}^{\vee}, B_{1}^{\vee}, w^{\vee}, z^{\vee})$ as an elements of $\mb M( W^{\vee}, V^{\vee} \otimes \det Q )$.

Using this we have the following as a corollary of Theorem \ref{ale}.
\begin{cor}
We have an isomorphism
$$
M^{\mbi \zeta}(\mbi w, \mbi v) \cong M_{X_{1}}(c_{-})
$$
for $\mbi \zeta \in \mc C_{-m}$ with enough large $m \gg 0$, where we put 
\begin{eqnarray*}
c_{-}=(w_{0} + w_{1})[X_{1}] + \left(2v_{0} - 2v_{1} + w_{1}\right) \frac{[C]}{2} - \left( v_{0} + \frac{w_{1}}{4}\right)[P]  + (w_{0} - w_{1}) [\ell_{\infty}^{1}] \in A(IX_{1}).
\end{eqnarray*}
\end{cor}
\proof
From the assumption, we see that $\mbi \zeta$ lies in the chamber $\mc C_{-\infty}$ in the notation in \S \ref{subsec:cham}.
Hence we can take $\mbi \zeta^{1} = \mbi \zeta [1] \in \mc C_{\infty}$ so that we have an isomorphism $\pi \colon M^{\mbi \zeta^{1}}(\mbi w[1], \mbi v[1]) \cong M_{X_{1}}(c[1])$ in Theorem \ref{ale}, where 
$$
c[1]=c_{-} \cdot \wt \ch (\mo_{X_{1}}(F -\ell_{\infty})).
$$
%\begin{eqnarray*}
%c[1]=(w_{1} + w_{0})[X_{1}] + \left(2v_{0} - 2v_{1} -w_{0} \right) \frac{[C]}{2} - \left( v_{1} + \frac{w_{0}}{4}\right)[P]  + ( w_{1} - w_{0} ) [\ell_{\infty}^{1}] \in A(IX_{1}).
%\end{eqnarray*}

%We put $\mc R_{0}^{-}=\mo_{X_{1}}, \mc R_{1}^{-}=\mo_{X_{1}}(\ell_{\infty} - F)$ and $\mc R^{-}= \mc R_{0}^{-} \oplus \mc R_{1}^{-}$.
Replacing $\mc R_{i}=\mo_{X_{1}} (i(F - \ell_{\infty}))$ with $\mc R_{i}^{-} = \mo_{X_{1}}( i (\ell_{\infty} - F))$ for $i=0,1$ in the complex \eqref{coh} in the proof of Theorem \ref{ale}, we get a $\tilde{T}$-equivariant isomorphism $\pi_{-} \colon M^{-}(\mbi w, \mbi v) \cong M_{X_{1}}(c_{-})$.
This follows from the following commutative diagram:
$$
\xymatrix{
M^{\mbi \zeta^{1}}(\mbi w[1], \mbi v[1]) \ar[r]^-{[1]} \ar[d]^{\pi} &M^{\mbi \zeta}(\mbi w, \mbi v) \ar[d]^{\pi_{-}}\\
M_{X_{1}}(c[1]) \ar[r]& M_{X_{1}}(c_{-}) 
}
$$
%$\vec{r}[1]=(r_{1}, r_{0}), k'= -k-r$ and $n'=n+\frac{r}{4}+\frac{k}{2}$, 
Here the top horizontal arrow is the isomorphism $[1]$ in \eqref{shift}, and the bottom one is the isomorphism induced by tensoring $\mo_{X_{1}}(\ell_{\infty} - F)$.
\endproof

%%%%%%%%%%%%%%%%%%%%%%%%%%%%%%%%%%%%%%%%%%%%%%%%%%%%%%%%%%%%%%%%%%%%%%%%%%%

\subsection{Integrations}
\label{subsec:inte}
We take $\mbi \zeta = \mbi 0 \in \R^{2}$, and put $M_{0}(\mbi w, \mbi v) = M^{\mbi 0}(\mbi w, \mbi v)$.
Then we have a proper map $\Pi \colon M^{\mbi \zeta}(\mbi w, \mbi v) \to M_{0}(\mbi w, \mbi v)$ for any $\mbi \zeta \in \R^{2}$.

\begin{prop}
The $\tilde{T}$-fixed points set $M_{0}(\mbi w, \mbi v)^{\tilde{T}}$ consists of one point.
\end{prop}
\proof 
If we take a representative $\mc A=(B, z, w) \in \mu^{-1}(0)$ of a point $p \in M_{0}(\mbi w, \mbi v)^{\tilde{T}}$.
Then for any $f \in \Gamma(\mu^{-1}(0), \mo_{\mu^{-1}(0)})^{G}$, we have $f( \mc A ) = f ( t \mc A )$ for any $t \in T^{2}$.
Since $\lim_{t\to 0} t \mc A=(0,z,0)$, we have $f(\mc A)=f (0, z, 0)$. 
Furthemore the closure of $G$-orbit of $(0,z,0)$ contains $(0,0,0)$.
Hence we see that $(0,0,0) \in \mu^{-1}(0)$ represents the same point $p$.
\begin{NB}
It is similarly proven as in \cite[Proposition 2.9 (3)]{NY1}.
\end{NB}
\endproof

We write the inclusion by $\iota \colon M_{0}(\mbi w, \mbi v)^{\tilde{T}} \to M_{0}(\mbi w, \mbi v)$.
For $\mbi \zeta$ in some chambers and $\psi \in A_{\tilde{T}}^{\ast}(M^{\mbi \zeta}(\mbi w, \mbi v))$, we define the integration over $M^{\mbi \zeta}(\mbi w, \mbi v)$ by 
$$
\int_{M^{\mbi \zeta}(\mbi w, \mbi v)} \psi = (\iota_{\ast})^{-1} \circ \Pi_{\ast} ( \psi \cap [M^{\mbi \zeta}(\mbi w, \mbi v)]) \in \mc S=\mb Q(\mbi \e, \mbi a, \mbi m), 
$$
where $[M^{\mbi \zeta}(\mbi w, \mbi v)]$ is the fundamental cycle.
Here $\mc S$ is a fractional field of $A_{\tilde{T}}^{\ast}(pt)=\mb Z[\mbi \e, \mbi a, \mbi m]$ as in the introduction.

Integrations over $M_{X_{\kappa}}(c)$ in \S \ref{subsec:part} are defined by these integrations via isomorphisms in Theorem \ref{ale}.

%%%%%%%%%%%%%%%%%%%%%%%%%%%%%%%%%%%%%%%%%%%%%%%%
%%%%%%%%%%%%%%%%%%%%%%%%%%%%%%%%%%%%%%%%%%%%%%%%

\subsection{ADHM descriptions of Partition functions}
\label{subsec:adhm3}
For $\mbi \zeta \in \R^{2}$ on a certain chamber, we consider a $\tilde{T}$-equivariant bundle
$$
\mc F_{r}(\mc V_{0})=\left( \mc V_{0} \otimes \frac{e^{m_{1}}}{\sqrt{t_{1}t_{2}}} \right) \oplus \cdots \oplus \left( \mc V_{0} \otimes \frac{e^{m_{2r}}}{\sqrt{t_{1}t_{2}}} \right)
$$
 on $M^{\mbi \zeta}(\mbi w, \mbi v)$. 
For fixed $k \in \frac{1}{2} \Z$ and $\mbi w \in \Z^{2}_{\ge 0}$, we consider $\mbi v=(v_0, v_0 + \frac{w_1}{2} + k)$, and take sums over $v_{0} \in \Z_{\ge 0}$.
Then by Theorem \ref{ale}, we have
\begin{eqnarray}
\label{ADHMpartition}
Z_{X_{0}}^{k}(\mbi \e, \mbi a, \mbi m, q)=\sum q^{n} \int_{M^{\mbi \zeta^{0}}(\mbi w, \mbi v)} e(\mc F_{r}(\mc V_{0})),
Z_{X_{1}}^{k}(\mbi \e, \mbi a, \mbi m, q)=\sum q^{n} \int_{M^{\mbi \zeta^{1}}(\mbi w, \mbi v)} e(\mc F_{r}(\mc V_{0})),
\end{eqnarray}
where $n = v_{0} + \frac{w_{1}}{4}$, and $\mbi \zeta^{0}, \mbi \zeta^{1}$ as in Theorem \ref{ale}.

\begin{figure}[h]
\includegraphics[scale=0.4]{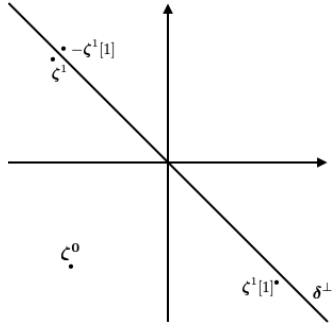}
\caption{stability parameter $\mbi \zeta^{0}$ and $\mbi \zeta^{1}$}
\label{zeta+}
\end{figure}

Our strategy to prove Theorem \ref{main} is the following.
For $k \le 0$, we show that wall-crossing across $\mbi \alpha_{m}^{\perp}$ for $m \ge 0$ lying between $\mbi \zeta^{0}$ and $\mbi \zeta^{1}$ does not change partition functions in the similar way to \cite{NY3}.

For $k \ge 0$, we analyze wall-crossing across $\mbi \delta^{\perp}$ lying between $\mbi \zeta^{1}$ and $-\mbi \zeta^{1}[1]$ in the similar way to \cite{O}.
Then we use an isomorphism $M^{-\mbi \zeta^{1}[1]}(\mbi w, \mbi v) \cong M^{\mbi \zeta^{1}[1]}(\mbi w, \mbi v)$ via the homomorphism $\tilde{T} \to \tilde{T}, (t_{1},t_{2}, e^{\mbi a}, e^{\mbi m}) \mapsto (t_{2}, t_{1}, e^{-\mbi a}, e^{\mbi m})$ by \eqref{dual}.
Furthermore, we can show that wall-crossing across $\mbi \alpha_{m}^{\perp}$ for $m<0$ lying between $\mbi \zeta^{1}[1]$ and $\mbi \zeta^{0}$ does not change partition functions, since this process is equivalent to the above wall-crossing between $\mbi \zeta^{0}$ and $\mbi \zeta^{1}$ for $k \le 0$ via the isomorphism $[1]$ in \S \ref{subsec:symm}.
This completes our proof of Theorem \ref{main}.

\section{Mochizuki method}
\label{sec:Moch}
We apply Mochizuki method \cite{Mo} to ADHM description in the previous section following \cite{NY3}.
It is the similar arguments as in \cite{O}.
Hence we often omit the detailed description.

%%%%%%%%%%%%%%%%%%%%%%%%%%%%%%%%%%%%%%%%%%%%%%%%%%%%%%%%%%%%%%%%%%%%%%%%

\subsection{Quiver description}
\label{subsec:quiv}

For later purpose, we modify the definition of ADHM data following \cite{C}.
As in the previous section, we consider ADHM data on $(W,V)$ for $\Z_{2}$-graded vector spaces $W= W_{0} \oplus W_{1}, V= V_{0} \oplus V_{1}$.

We introduce a quiver with relations $\Gamma = (\vec{\Gamma}, I)$ as follows, where $\vec{\Gamma}$ is a quiver, and $I$ is an ideal of the path algebra of $\vec{\Gamma}$.
The set $\vec{\Gamma}_{0}$ of vertex consists of $0, 1$ and $\infty$.
The set $\vec{\Gamma}_{1}$ of arrows consists of $\alpha \colon 0 \to 1, \beta \colon 1 \to 0, \gamma_{1}, \ldots, \gamma_{r_{0}} \colon \infty \to 0, \delta_{1}, \ldots, \delta_{r_{1}} \colon \infty \to 1$, and their converse $\alpha^{\ast} \colon 1 \to 0, \beta^{\ast} \colon 0 \to 1, \gamma_{1}^{\ast}, \ldots, \gamma_{r_{0}}^{\ast} \colon 0 \to \infty, \delta_{1}^{\ast}, \ldots, \delta_{r_{1}}^{\ast} \colon 1 \to \infty$. 
The ideal $I$ is generated by $\beta \beta^{\ast} - \alpha^{\ast} \alpha + \sum_{i=1}^{r_{0}} \gamma_{i} \gamma_{i}^{\ast}$ and $\alpha \alpha^{\ast} - \beta^{\ast} \beta + \sum_{i=1}^{r_{1}} \delta_{i} \delta_{i}^{\ast}$.

A $\Gamma$-representation $\mc A$ consists of finite dimensional vector spaces $\mc A_{0}, \mc A_{1}$ and $\mc A_{\infty}$ corresponding to each vertices in $\vec{\Gamma}_{0}$, and linear maps $\mc A_{a} \colon \mc A_{s(a)} \to \mc A_{t(a)}$ for each arrow $a \in \vec{\Gamma}_{1}$, where $s(a)$ and $t(a)$ are source and target of an arrow $a$.

We identify ADHM data and $\Gamma$-representations as follows.
We put 
\begin{equation}
\label{Gamma}
\mc A_{0}=V_{0}, \mc A_{1} = V_{1},  \mc A_{\infty}=
\begin{cases}
\C \text{ if }W \neq 0, \\
0 \text{ if }W=0.
\end{cases}
\end{equation}
We consider the graded vector space $V[1]$ with $V[1]_{0}=V_{1}$ and $V[1]_{1} = V_{0}$, and take a basis $\mbi e_{1} , \ldots, \mbi e_{r_{0}}$ and $\mbi e_{r_{0}+1}, \ldots, \mbi e_{r_{0} + r_{1}}$ of $W_{0}$ and $W_{1}$, and their dual basis $\mbi e_{1}^{\ast} , \ldots, \mbi e_{r_{0}}^{\ast}$ and $\mbi e_{r_{0}+1}^{\ast}, \ldots, \mbi e_{r_{0} + r_{1}}^{\ast}$.
Then from $\Gamma$-representations $\mc A$, we can assign ADHM data $(B_{1}, B_{2}, z, w)$ by
$$
B_{1} =
\begin{bmatrix}
0 & \mc A_{\beta}\\
\mc A_{\alpha} & 0
\end{bmatrix},
B_{2} =
\begin{bmatrix}
0 & \mc A_{\alpha^{\ast}}\\
\mc A_{\beta^{\ast}} & 0
\end{bmatrix}
\in \Hom_{\Z_{2}}(V, V[1]),
$$
$z= \sum_{i=1}^{r_{0}} \mc A_{\gamma_{i}} \mbi e_{i}^{\ast} + \sum_{j=1}^{r_{1}} \mc A_{\delta_{j}} \mbi e_{r_{0} + j}^{\ast} \in \Hom_{\Z_{2}}(W, V)$, and $w= \sum_{i=1}^{r_{0}} \mbi e_{i} \mc A_{\gamma_{i}^{\ast}} + \sum_{j=1}^{r_{1}} \mbi e_{r_{0}+j} \mc A_{\delta_{j}^{\ast}} \in \Hom_{\Z_{2}}( V, W )$.
Conversely, we can assign $\Gamma$-representation from ADHM data on $(W, V)$ by the above equations.

For $\mbi \zeta =(\zeta_{0}, \zeta_{1}) \in \mb R^{2}$, we define stability of $\Gamma$-representations as follows.
For any sub-representation $S=S_{0} \oplus S_{1} \oplus S_{\infty}$ of $\Gamma$-representations $\mc A = V_{0} \oplus V_{1} \oplus \C$, we put 
$$
\mbi \zeta(S)=\zeta_{0} \dim S_{0} + \zeta_{1} \dim S_{1} - (\zeta_{0} v_{0} + \zeta_{1} v_{1} ) \dim S_{\infty}, 
$$
where $v_{0}=\dim V_{0}, v_{1} = \dim V_{1}$.

\begin{defn}
A $\Gamma$-representation $\mc A$ is said to be $\mbi \zeta$-semistable if for any sub-representation $S$ of $\mc A$, we have $\mbi \zeta(S) \le 0$.
In addition, if we have $\mbi \zeta(S) < 0$ for any $S \neq \mbi 0, \mc A$, we say $\mc A$ is $\mbi \zeta$-stable.
\end{defn}

This coincides with Definition \ref{ADHM-stability} of stability for ADHM data via the above identification between ADHM data and $\Gamma$-representations.
We remark that $\mbi \zeta(S)$ is not equal to $(\mbi \zeta, (\dim S_{0}, \dim S_{1})) = \zeta_{0} \dim S_{0} + \zeta_{1} \dim S_{1}$ when $\dim S_{\infty} \neq 0$.
So we often consider the quotient space $V/ S$ in such a case.

%%%%%%%%%%%%%%%%%%%%%%%%%%%%%%%%%%%%%%%%%%%%%%%%%%%%%%%%%%%%%%%%%%%%%%%%

\subsection{ADHM data with full flags}
\label{subsec:adhm4}

In the rest of this section, we fix one of the walls $\mk D=\pm \mathfrak{D}_{\mbi \alpha}$ for $\mbi \alpha =(\alpha_{0}, \alpha_{1}) \in R^{+}$ defined in \S \ref{subsec:cham}.
We choose $i_{0} \in \Z_{2}$ such that $\alpha_{i_{0}} \neq 0$. 
When $\mbi \alpha \neq (1,0), (0,1)$, we can choose both $i_{0} = 0$ and $1$.
However, for simplicity, when $\mbi \alpha$ is an imaginary root $p \mbi \delta=(p,p)$ for $p \in \Z_{>0}$, we always choose $i_{0} =1$, and set $p=1$, that is, $\mbi \alpha = \mbi \delta$.

%In the following, we also consider pairs $(\mc A, F^{\bullet})$, where $\mc A=(B,z,w)$ does not necessarily satisfy the relation $[B \wedge B] +zw=0$, and $F^{\bullet}$ may have repetitions satisfying $\dim F^{i_{0}}/F^{i-1} =0$ or $1$ and $F^{i_{0}}=V$ for $i \gg 0$. 

We have two chambers adjacent to the wall $\mk D$, and write by $\mc C$ the one whose element $\mbi \zeta$ satisfy $(\mbi \zeta, \mbi \alpha)<0$, and by $\mc C'$ the other one.
We take $\mbi \zeta^{\mk D}$ on the wall $\mk D$.
For $\Z_{2}$-graded vector spaces $W=W_{0} \oplus W_{1}, V=V_{0} \oplus V_{1}$, we consider pairs $(\mc A, F^{\bullet})$ of ADHM data $\mc A=(B,z,w)$ on $(W, V)$ and full flags $F^{\bullet}$ of $V_{i_{0}}$ for fixed $i_{0}$ as above.
We consider ADHM data $\mc A$ as $\Gamma$-representations as in \S \ref{subsec:quiv}.

\begin{defn}
\label{ellstab}
For $\ell \ge 0$, a pair $(\mc A, F^{\bullet})$ is said to be $(\mc C, \ell)$-stable if $\mc A$ is $\mbi \zeta^{\mk D}$-semistable and any sub-representation $S = S_{0} \oplus S_{1} \oplus S_{\infty}$ of $\mc A$ with $\mbi \zeta^{\mk D}(S)=0$ satisfies the following two conditions:
\begin{enumerate}
\item[\textup{(1)}] If $S_{\infty}=0$ and $S \neq 0$, we have $S_{i_{0}} \cap F^{\ell} =0$.
\item[\textup{(2)}] If $S_{\infty}=\C$ and $S \neq \mc A$, we have $F^{\ell} \not\subset S_{i_{0}}$. 
\end{enumerate}
\end{defn}
We write by $\wt{M}^{\mc C,\ell}( \mbi w, \mbi v)$ moduli of $(\mc C, \ell)$-stable ADHM data on $(W, V)$ with full flags of $V_{i_{0}}$, which will be constructed in the next subsection.
%For $\mc C=\mc C_{m}$ and $\mc C_{+}$ in \S \ref{subsec:cham}, we put $\wt{M}^{m,\ell}( \mbi w, \mbi v)=\wt{M}^{\mc C_{m}, \ell}( \mbi w, \mbi v)$ and $\wt{M}^{+,\ell}( \mbi w, \mbi v)=\wt{M}^{\mc C_{+}, \ell}( \mbi w, \mbi v)$.
We remark that $\mbi \zeta^{\mk D}(S)=0$ implies that $(\dim S_{0}, \dim S_{1})$ is proportional to the root $\mbi \alpha$ defining the wall $\mk D$, and when $\ell=0$ (resp. $\ell=v_{i_{0}}$), an object $(\mc A, F^{\bullet})$ is $(\mc C, \ell)$-stable if and only if $\mc A=(B,z,w)$ is $\mc C$-stable (resp. $\mc C'$-stable).
Hence we see that $\wt{M}^{\mc C, 0}(\mbi w, \mbi v)$ and $\wt{M}^{\mc C, v_{i_{0}}}(\mbi w, \mbi v)$ are full flag bundles of tautological bundles $\mc V_{i_{0}}$ on $M^{\mc C}(\mbi w, \mbi v)$ and $M^{\mc C'}(\mbi w, \mbi v)$ respectively.  

We also interpret $(\mc C, \ell)$-stability in terms of $\Gamma$-representations as follows.
We take $\mbi \eta = (\eta_{1}, \ldots, \eta_{v_{i_{0}}}) \in (\mb Q_{>0})^{v_{i_{0}}}$, and for any sub-representation $S$ of $\mc A$, we put
$$
\mu_{\mbi \zeta, \mbi \eta}(S)= \frac{\mbi \zeta(S) + \sum_{j=1}^{v_{i_{0}}} \eta_{j} \dim (S_{i_{0}} \cap F^{j})}{\rk S},
$$
where $\rk S = \dim S_{0} + \dim S_{1} + \dim S_{\infty}$.
We say that $(\mc A, F^{\bullet})$ is $(\mbi \zeta, \mbi \eta)$-semistable if for any non-zero proper sub-representation $S$, we have
\begin{eqnarray}
\label{slope}
\mu_{\mbi \zeta, \mbi \eta}(S) \le \mu_{\mbi \zeta, \mbi \eta}(\mc A)= \frac{\sum_{j=1}^{v_{i_{0}}} j \eta_{j}}{\rk \mc A}.
\end{eqnarray}
If inequality is always strict, we say that $(\mc A, F^{\bullet})$ is $(\mbi \zeta, \mbi \eta)$-stable.

We consider the following condition
\begin{eqnarray}
\label{conda}
\sum_{j=1}^{v_{i_{0}}} j \eta_{j} < \underset{( \mbi \zeta^{\mk D}, \mbi s) \neq 0}{\min} \frac{ |( \mbi \zeta, \mbi s) | }{\rk \mc A},\\
\label{condb}
\rk \mc A \sum_{j=\ell +1}^{v_{i_{0}}} j \eta_{j} 
<
 \min \left( \sum_{j=1}^{\ell} j \eta_{j} - \frac{\rk \mc A}{\rk \mbi \alpha} (\mbi \zeta, \mbi \alpha), \right.
\left. -\sum_{j=1}^{\ell} j \eta_{j} + \frac{\rk \mc A}{\rk \mbi \alpha} (\mbi \zeta, \mbi \alpha)+\eta_{\ell}\right),\\
\label{condc}
\eta_{k} 
>
\rk \mc A \sum_{j=k+1}^{v_{i_{0}}} j \eta_{j} \text{ for } k=1, \ldots, v_{i_{0}},\\
\label{condd}
\sum_{j=1}^{v_{i_{0}}} k_{j} \eta_{j} \neq 0 \text{ for any  } (k_{1}, \ldots, k_{v_{i_{0}}}) \in \Z^{v_{i_{0}}} \setminus \lbrace 0\rbrace \text{ with } |k_{j}| \le n^{2},
\end{eqnarray}
where in \eqref{conda} minimum is taken over the set of all $\mbi s = ( s_{0}, s_{1} ) \in \Z^{2}$ with $0 \le s_{0} \le v_{0}, 0 \le s_{1} \le v_{1}$, and $(\mbi \zeta^{\mk D}, \mbi s ) \neq 0$.
We call \eqref{condd} $2$-stability condition following \cite{Mo}.

By the similar arguments as in \cite[\S 4.1]{O}, we have the following
\begin{prop}
[\protect{
\cite[Proposition 4.2.4]{Mo}, 
\cite[Lemma 5.6]{NY3}}]
\label{ell}
We take $\mbi \zeta \in \mc C'$, and assume that $(\mbi \zeta, \mbi \eta)$ satisfies \eqref{conda} and \eqref{condb}. 
Then for $(\mc A,F^{\bullet})$, the $(\mc C, \ell)$-stability is equivalent to the $(\mbi \zeta, \mbi \eta)$-stability.
Furthermore the $(\mbi \zeta, \mbi \eta)$-semistability automatically implies that the $(\mbi \zeta,\mbi \eta)$-stability.
\end{prop}
\proof
By \eqref{conda}, for any sub-representation $S$ of $\mc A$ with $\mbi \zeta^{\mk D} (S) \neq 0$, we have $\mu_{\mbi \zeta, \mbi \eta}(S) \neq \sum_{j=1}^{v_{i_{0}}} j \eta_{j} / \rk \mc A$, and $\mbi \zeta^{\mk D}(S) < 0$ if and only if $\mu_{\mbi \zeta, \mbi \eta}(S) < \sum_{j=1}^{v_{i_{0}}} j \eta_{j} / \rk \mc A$.
By \eqref{condb}, condition (1) and (2) in Definition \ref{ellstab} holds if and only if \eqref{slope} holds for any sub-representation $S$ of $\mc A$ with $\mbi \zeta^{\mk D}(S)=0$. 
Furthermore in both cases where $\mbi \zeta^{\mk D}(S) \neq 0$ and $\mbi \zeta^{\mk D}(S)=0$, we see that inequality \eqref{slope} is strict.
Hence the assertions hold.
\endproof

We can choose $\mbi \zeta, \mbi \eta$ satisfying \eqref{conda}, \eqref{condb}, and \eqref{condd} as follows.
First, we determine neighborhood of $(\mbi \zeta^{\mk D}, (0, \ldots, 0))$ in $\mb R^{2} \times \mb Q^{v_{i_{0}}}$ in which any $(\mbi \zeta, \mbi \eta)$ satisfy \eqref{conda} and \eqref{condd}. 
We choose $\mbi \zeta$ near $\mbi \zeta^{\mk D}$ and $\eta_{1}, \ldots, \eta_{\ell} > 0$ small enough such that $(\mbi \zeta, (\eta_{1}, \ldots, \eta_{\ell}, 0, \ldots, 0))$ satisfies \eqref{condc} and
$$
\sum_{j=1}^{\ell} j \eta_{j} - \eta_{\ell} < \frac{\rk \mc A }{\rk \mbi \alpha} \mbi \zeta(\mbi \alpha) < \sum_{j=1}^{\ell} j \eta_{j}.
$$
Finally we take $\eta_{\ell+1}, \ldots, \eta_{v_{i_{0}}} > 0$ satisfying \eqref{condb} and \eqref{condc} such that $(\mbi \zeta, \eta_{1}, \ldots, \eta_{v_{i_{0}}})$ belongs to the above neighborhood.
We use \eqref{condc} and \eqref{condd} in \S \ref{subsec:dire}.

%%%%%%%%%%%%%%%%%%%%%%%%%%%%%%%%%%%%%%%%%%%%%%%%%%%%%%%%%%%%%%%%%%%%%%%%

\subsection{Enhanced master space}
\label{subsec:enha}

We introduce an enhanced master space using ADHM description.
For $\Z_{2}$-graded vector spaces $W=W_{0} \oplus W_{1}, V=V_{0} \oplus V_{1}$, we consider pairs $(\mc A, F^{\bullet})$ of ADHM data $\mc A=(B,z,w)$ on $(W, V)$ and full flags $F^{\bullet}$ of $V_{i_{0}}$.

Let $[v_{i_{0}}]$ denote the set $\lbrace 1, \ldots, v_{i_{0}}\rbrace$ of integers, and $Fl=Fl(V_{i_{0}}, [v_{i_{0}}])$ denote the full flag variety of $V_{i_{0}}$, where $v_{i_{0}} = \dim V_{i_{0}}$.
We consider natural projections $\rho_{j} \colon Fl \to G_j=Gr(V_{i_{0}},j)$ to Grassmanian manifolds $G_j$ of $j$-dimensional subspace of $V_{i_{0}}$ and pull-backs $\rho_{j}^{\ast}\mo_{G_j}(1)$ of polarizations $\mo_{G_{j}}(1)$ by Plucker embeddings.

In the following, we fix $\ell \in [v_{i_{0}}]$, and choose $\mbi \zeta^{-} \in \mc C, \mbi \zeta \in \mc C'$, and $\mbi \eta \in (\Q_{>0})^{v_{i_{0}}}$ as follows.  
$|\mbi \zeta|, |\mbi \eta|$ are enough smaller than $|\mbi \zeta^{-}|$ so that any $(\mc A, F^{\bullet})$ is $(\mbi \zeta^{-}, \mbi \eta)$-stable if and only if $\mc A$ is $\mbi \zeta^{-}$-stable, and $(\mbi \zeta, \mbi \eta)$ satisfy the conditions \eqref{conda}, \eqref{condb}, and \eqref{condd}.
We take a positive integer $k$ enough divisible such that $k \mbi \zeta, k \mbi \zeta^{-}$ and $k\mbi \eta$ are all integer valued, and consider ample $G$-linearizations 
\begin{eqnarray*}
L_{+}
&=&
\left(\mo_{\mb M} \otimes (\det V)^{\otimes k \mbi \zeta} \right)\boxtimes \bigotimes_{j=1}^{n}\rho_{j}^{\ast}\mo_{G_{j}}(k\eta_{j}),\\ 
L_{-}
&=&
\left(\mo_{\mb M} \otimes (\det V)^{\otimes k \mbi \zeta^{-}} \right)\boxtimes \bigotimes_{j=1}^{n}\rho_{j}^{\ast}\mo_{G_{j}}(k\eta_{j})
\end{eqnarray*}
on $\wt{\mb M}=\wt{\mb M}(W,V)=\mb M(W, V) \times Fl$.
We consider the composition $\tilde{\mu} \colon \wt{\mb M} \to \mb L$ of the projection $\wt{\mb M} \to \mb M$ and $\mu \colon \mb M \to \mb L$ in \eqref{mu}, and semistable loci $\tilde{\mu}^{-1}(0)^{+}$ and $\tilde{\mu}^{-1}(0)^{-}$ with respect $L_{+}$ and $L_{-}$ respectively. 

We put $\wh{\mb M}=\wh{\mb M}(W,V)= \mb P(L_{-}\oplus L_{+})$ and consider a composition $\hat{\mu}\colon \wh{\mb M} \to \mb L$ of the projection $\wh{\mb M} \to \mb M$ and $\mu \colon \mb M \to \mb L$.
Then we have a natural $G=\GL(V_{0}) \times \GL(V_{1})$-action on $\wh{\mb M}$ compatible with $\hat{\mu}$.
We also write by $\mo(1)$ the restriction of the tautological bundle $\mo(1)$ to $\hat{\mu}^{-1}(0)$, which defines semistable locus $\hat{\mu}^{-1}(0)^{ss}$.  
We define an {\it enhanced master space} by $\mc M=[\hat{\mu}^{-1}(0)^{ss}/G]$.
The projection $\hat{\mu}^{-1} (0) \to \mu^{-1}(0)$ induces a proper morphism $\Pi \colon \mc M \to M_{0}(\mbi w, \mbi v)$.

We have a $\C^{\ast}_{\hbar}$-action on $\mc M$ defined by 
\begin{eqnarray}
\label{act0}
\left(\mc A, F^{\bullet}, [x_{-},x_{+}] \right) \mapsto \left(\mc A, F^{\bullet}, [e^{\hbar} x_{-},x_{+}] \right),
\end{eqnarray}
where $[x_{-}, x_{+}]$ is the homogeneous coordinate of $\PP(L_{-} \oplus L_{+})$.
We also note that we have 
\begin{eqnarray}
\label{kzeta}
(k_{0}, k_{1} )= k( \mbi \zeta^{-} - \mbi \zeta)
\end{eqnarray} 
in \S \ref{subsec:enha0}.
Finally, we get $\C^{\ast}_{\hbar}$-fixed points sets $\mc M^{\C^{\ast}_{\hbar}}$.

%%%%%%%%%%%%%%%%%%%%%%%%%%%%%%%%%%%%%%%%%%%%%%%%%%%%%%%%%%%%%%%%%%%%%%%%%%%

\subsection{Direct sum decompositions of fixed points sets}
\label{subsec:dire}
In this section, we follow the similar argument as in \cite[\S 4.3]{O}.
 
For ADHM data with full flags $(\mc A, F^{\bullet})$, if we have a direct sum decomposition $(\mc A, F^{\bullet})=(\mc A_{\flat}, F^{\bullet}_{\flat}) \oplus (\mc A_{\sharp}, F^{\bullet}_{\sharp})$, then we put $I_{\alpha}=\lbrace j \in [v_{i_{0}}] \mid F^{j}_{\alpha}/ F^{j-1}_{\alpha} \neq 0 \rbrace$ for $\alpha=\flat, \sharp$ so that $[v_{i_{0}}]=I_{\flat} \sqcup I_{\sharp}$.
Here we allow for flags $F_{\flat}^{\bullet}, F_{\sharp}^{\bullet}$ to have repetitions, and assume that $(\mc A_{\sharp})_{\infty}=0$.
The data $\mk J=(I_{\flat}, I_{\sharp})$ are called the {\it decomposition type}.

By $2$-stability condition \eqref{condd}, we see that $x \in \hat{\mu}^{-1}(0)^{ss} \setminus \left( \mb P(L_{-}) \sqcup \mb P(L_{+}) \right)$ over $(\mc A,F^{\bullet}) \in \wt{\mb M}(r,n)$ represents a $\C^{\ast}_{\hbar}$-fixed point in $\mc M$ if and only if we have a decomposition $(\mc A, F^{\bullet})=(\mc A_{\flat}, F^{\bullet}_{\flat}) \oplus (\mc A_{\sharp}, F^{\bullet}_{\sharp})$ safisfying the following conditions (cf \cite[Lemma 5.16]{NY3}).
The decomposition type $\mk J=(I_{\flat}, I_{\sharp})$ satisfies $\min(I_{\sharp}) \le \ell$, and there exists a $\mbi \zeta'$ on the segment connecting $\mbi \zeta^{-}$ and $\mbi \zeta$ such that 
\begin{eqnarray}
\label{slope2}
\mu_{\mbi \zeta', \mbi \eta}(\mc A_{\flat}) = \mu_{\mbi \zeta', \mbi \eta}(\mc A_{\sharp}) =  \mu_{\mbi \zeta', \mbi \eta}(\mc A), 
\end{eqnarray}
and both $(\mc A_{\flat}, F_{\flat}^{\bullet})$ and $(\mc A_{\sharp}, F_{\sharp}^{\bullet})$ are $(\mbi \zeta', \mbi \eta)$-stable.

Since $\mbi \eta$ is smaller enough than $|\mbi \zeta^{-}|, |\mbi \zeta'|$, there exists such a $\mbi \zeta'$ such that the last equation holds if and only if $\mbi \zeta^{\mk D}(\mc A_{\sharp})=0$.
Hence the dimension vector of $\mc A_{\sharp}$ is multiple of $\mbi \alpha$.
For a fixed wall $\mathfrak{D}=\pm \mathfrak{D}_{\mbi \alpha}$ for $\mbi \alpha=(\alpha_{0}, \alpha_{1}) \in R_{+}$ and $\ell \in [v_{i_{0}}]$, we introduce the set 
\begin{eqnarray}
\label{decdata}
\mc D^{\ell}(v_{i_{0}}, \alpha_{i_{0}}) = \left\lbrace \mk J =(I_{\flat}, I_{\sharp}) \mid [v_{i_{0}} ]= I_{\flat} \sqcup I_{\sharp}, |I_{\sharp}| \in \alpha_{i_{0}} \Z_{>0}, \min (I_{\sharp}) \le \ell \right\rbrace.
\end{eqnarray}
For $\mk J =(I_{\flat}, I_{\sharp}) \in \mc D^{\ell}(v_{i_{0}}, \alpha_{i_{0}})$, we put $p = |I_{\sharp}| / \alpha_{i_{0}}$.

\begin{defn}
\label{+}
$(\mc A_{\sharp}, F_{\sharp}^{\bullet})$ is said to be $(\mathfrak{D}, +)$-stable if $(\mc A_{\sharp})_{\infty}=0$, $\mbi \zeta^{\mk D}(\mc A_{\sharp})=0$, $\mc A_{\sharp}$ is $\mbi \zeta^{\mk D}$-semistable, and for any proper sub-representation $S =S_{0} \oplus S_{1}$ of $\mc A_{\sharp}$ with $\mbi \zeta^{\mk D}(S)=0$, we have $F_{\sharp}^{1} \cap S_{i_{0}} = 0$.
\end{defn}

Suppose that we are given a pair $(\mc A_{\flat}, F^{\bullet}_{\flat}), (\mc A_{\sharp}, F^{\bullet}_{\sharp})$ with the decomposition type $\mk J=(I_{\flat}, I_{\sharp}) \in \mc D^{\ell}(v_{i_{0}}, \alpha_{i_{0}})$ satisfying $(\mc A_{\flat})_{\infty}=\C, (\mc A_{\sharp})_{\infty}=0$, and $\mbi \zeta^{\mk D}(\mc A_{\sharp})=0$. 
We take $\mbi \zeta'$ on the segment connecting $\mbi \zeta^{-}$ and $\mbi \zeta$ satifsying \eqref{slope2}.
By \eqref{condc}, we have the following.
\begin{lem}
[\protect{\cite[Proposition 4.4.4]{Mo},\cite[Lemma 5.26]{NY3}}]
\label{lem4}
We have the following.
\begin{enumerate}
\item[\textup{(1)}] $(\mc A_{\flat}, F_{\flat}^{\bullet})$ is $(\zeta', \mbi \eta)$-stable if and only if it is $(\mc C, \min(I_{\sharp})-1)$-stable.

\item[\textup{(2)}] $(\mc A_{\sharp}, F_{\sharp}^{\bullet})$ is $(\zeta', \mbi \eta)$-stable if and only if it is $(\mk D, +)$-stable.
\end{enumerate}
\end{lem}
\proof
It is similarly proven as in \cite[Lemma 5.26]{NY3}.
\endproof

%%%%%%%%%%%%%%%%%%%%%%%%%%%%%%%%%%%%%%%%%%%%%%%%%%%%%%%%%%%%%%%%%%%%%%%%%%%

\subsection{$(\mk D, +)$-stability and $\mbi \zeta$-stability}
\label{subsec:plusstab}
We put $\mbi w_{\sharp}=(w_{\sharp 0}, w_{\sharp 1})$, where 
\begin{eqnarray}
\label{sharp}
w_{\sharp i_{0}}=1, w_{\sharp i_{0}+1}=0.
\end{eqnarray}
We compare $(\mk D, +)$-stability on $\mb M(\mbi 0, p \mbi \alpha) \times Fl(V_{i_{0}}, [v_{i_{0}}])$ and $\mbi \zeta$-stability on $\mb M(\mbi w_{\sharp}, p \mbi \alpha)$, where we recall that $\mbi \zeta$ is in a chamber adjacent to the wall $\mk D \subset \mbi \alpha^{\perp}$, and satisfies $(\mbi \zeta, \mbi \alpha) < 0$.

\begin{lem}
\label{d+vszeta}
For ADHM data $\mc A_{\sharp} = (B, z, w) \in \mb M( \mbi w_{\sharp}, p \mbi \alpha)$, the following hold. \\
(1) $\mc A_{\sharp} = (B, z, w)$ are $\mbi \zeta$-stable if and only if they are $\mbi \zeta^{\mk D}$-semistable and there exists no proper sub-graded vector space $S$ of $V$ such that $B(Q^{\vee} \otimes S) \subset S$, $\im z \subset S$ and $\mbi \zeta^{\mk D}(S)=0$. \\
(2) When $\mbi \alpha = \mbi \alpha_{m}$ for $m \in \Z$, if $w=0$, then $(B, 0, 0) \in \mb M( 0, p \mbi \alpha)$ is $\mbi \zeta^{\mk D}$-semistable ADHM data.\\
(3) When $\mbi \alpha = \mbi \delta$, if $\mc A_{\sharp} = (B, z, w)$ are stable, then we have $w=0$.
\end{lem}
\proof
(1) This follows directly from Definition \ref{ADHM-stability} and choice of $\mbi \zeta$ as noted above. \\
(2) This follows from (1).\\
(3) For $\mbi a \in Q^{\vee}$, we consider a linear map $B_{\mbi a} \colon  V \to V$ defined by $B_{\mbi a}(\mbi v) = B(\mbi a \otimes \mbi v)$ for $\mbi v \in V$.
Then $\mbi \zeta^{\mk D}$-semistablity implies that $B_{\mbi a}|_{V_{0}} \colon V_{0} \to V_{1}$ is an isomorphism for some $\mbi a \in Q^{\vee}$.
In particular, we can take $(1,0)$, or $(0,1)$ as $\mbi a$.
Via this $B_{\mbi a}$ for $\mbi a= (1,0)$, or $(0,1)$, we identify $V_{0}$ and $V_{1}$ to get ADHM data on $(W_{i_{0}}, V_{i_{0}})$, and the condition in (1) implies that they are stable.
Then $w=0$ follows from \cite[Proposition 2.8]{N2}.
\endproof

For ADHM data $(\mc A_{\sharp}, F_{\sharp}^{\bullet})$ on $(0, V_{\sharp})$ with full flags of $V_{\sharp i_{0}}$, we take a generator of $F_{\sharp}^{1}$, that is, a non-zero element $f_{\sharp}^{1}$ in $F_{\sharp}^{1}$. 
We consider a new ADHM data $\mc A_{\sharp}^{+}=(B, z, w)$ on $(W_{\sharp}, V_{\sharp})$ as follows, where $W_{\sharp}=W_{\sharp 0} \oplus W_{\sharp 1}$ is a $\Z_{2}$-graded vector space such that 
$$
\mbi w_{\sharp}=
(
\dim W_{\sharp 0}, 
\dim W_{\sharp 1} 
).
$$
We define $B \in \Hom_{\Z_{2}}(Q ^{\vee} \otimes V_{\sharp}, V_{\sharp})$ by the same data as in $\mc A_{\sharp}$, and $z \in \Hom_{\Z_{2}}(W_{\sharp}, V_{\sharp})$ by $z(1)=f_{\sharp}^{1}$, where $1$ is a generator of $W_{\sharp i_{0}}$, and $w=0 \in \Hom_{\Z_{2}}(\det Q^{\vee} \otimes V_{\sharp}, W_{\sharp})$.

\begin{lem}
\label{d+}
ADHM data $(\mc A_{\sharp}, F_{\sharp}^{\bullet})$ are $(\mathfrak{D},+)$-stable if and only if ADHM data $\mc A_{\sharp}^{+}=(B,z,0)$ are $\mbi \zeta$-stable.
\end{lem}
\proof This follows from Definition \ref{+} and Lemma \ref{d+vszeta} (1).

From this lemma and Lemma \ref{d+vszeta} (2), $(\mathfrak{D},+)$-stable objects $(\mc A_{\sharp}, F_{\sharp}^{\bullet})$ are parametrized by the full flag bundle of a quotient of a tautological homomorphism $W_{\sharp i_{0}} \otimes \mo_{M^{\mbi \zeta}(\mbi w_{\sharp}, p \mbi \alpha)} \to \mc V_{\sharp i_{0}}$ over $M^{\mbi \zeta}(\mbi w_{\sharp}, p \mbi \alpha)$.

%%%%%%%%%%%%%%%%%%%%%%%%%%%%%%%%%%%%%%%%%%%%%%%%%%%%%%%%%%%%%%%%%%%%%%%%%%%

\subsection{Moduli stacks parametrizing destabilizing objectts}
\label{subsec:modu2}
For $\mk J = (I_{\flat}, I_{\sharp}) \in \mc D^{\ell}(v_{i_{0}}, \alpha_{i_{0}})$, we see from the previous subsection that elements $(\mc A, F^{\bullet} ) \in \mc M_{\mk J}$ are decomposed into direct sums of $(\mc C, \min(I_{\sharp})-1)$-stable objects $(\mc A_{\flat}, F^{\bullet}_{\flat})$, which are parametrized by $\wt M^{\mc C, \min(I_{\sharp}) -1}(\mbi w, \mbi v - p \mbi \alpha)$,  and $(\mk D, +)$-stable objects $(\mc A_{\sharp}, F^{\bullet}_{\sharp})$.
In this subsection, we consider moduli spaces parametrizing latter objects.

We fix a direct sum decomposition $V=V_{\flat} \oplus V_{\sharp}$ of $\Z/2\Z$-graded vector spaces such that $V_{\sharp}=V_{\sharp 0} \oplus V_{\sharp 1}$ with $V_{\sharp 0}=\C^{p \alpha_{0}}$ and $V_{\sharp 1} = \C^{p \alpha_{1}}$.
We consider a moduli stack 
$$
\wt M_{\mbi \alpha}^{p}= \left[ \left( \tilde{\mu}^{-1}(0)^{(\mathfrak{D}, +)} \times \C^{\ast}_{\rho_{\sharp}}\right) / \GL(V_{\sharp 0}) \times \GL(V_{\sharp 1} ) \right]
$$ 
parametrizing tuples $(X_{\sharp}, F^{\bullet}_{\sharp}, \rho_{\sharp})$ of $(\mathfrak{D}, +)$-stable pairs $(X_{\sharp}, F_{\sharp}^{\bullet})$ and orientations 
$$
\rho_{\sharp} \colon \det V_{\sharp 0}^{\otimes k(\zeta_{0} - \zeta_{0}^{-})} \otimes \det V_{\sharp 1}^{\otimes k(\zeta_{1} - \zeta_{1}^{-})} \cong \C.
$$ 
Here $\tilde{\mu} \colon \wt{\mb M}(0, V_{\sharp}) \to \mb L(V_{\sharp})$ is defined in \S \ref{subsec:enha}, $\tilde{\mu}^{-1}(0)^{(\mathfrak{D}, +)}$ is the $(\mathfrak{D}, +)$-stable locus of $\tilde{\mu}^{-1}(0)$, and $\GL(V_{\sharp 0}) \times \GL(V_{\sharp 1})$ acts naturally.
In the following, we put 
$$
D= ( k(\mbi \zeta - \mbi \zeta^{-}),  \mbi \alpha) =k \left( (\zeta_{0} - \zeta_{0}^{-}) \alpha_{0} + (\zeta_{1} - \zeta_{1}^{-} ) \alpha_{1}\right).
$$

We consider a vector bundle $\mc V_{\sharp i_{0}}^{\vee} \otimes \mc W_{\sharp i_{0}}$ and write by $w \in \Gamma( M^{\mc C}(\mbi w_{\sharp}, p \mbi \alpha), \mc V_{\sharp i_{0}}^{\vee} \otimes \mc W_{\sharp i_{0}})$ the induced section from $w \in \Hom_{\C} (V_{\sharp i_{0}}, W_{\sharp i_{0}})$.
We put
\begin{eqnarray}
\label{malpha}
M_{\mbi \alpha}^{p} =  w^{-1} (\mbi 0).
\end{eqnarray}
When $\mbi \alpha$ is an imaginary root, we have $M_{\mbi \alpha}^{p} = M^{\mc C}(\mbi w_{\sharp}, p \mbi \alpha)$ by Lemma \ref{d+vszeta} (3).
When $\mbi \alpha$ is a real root, we describe $M_{\mbi \alpha}^{p}$ in the next subsection. 

For $k=0,1$, we also write by $\mc V_{\sharp k}$ the restriction of tautological bundle $\mc V_{\sharp k}$ to $M_{\mbi \alpha}^{p}$.
By Lemma \ref{d+}, we have the following proposition.

\begin{prop}
%[\protect{\cite[Proposition 5.9]{NY3}}]
\label{destab}
$\wt M_{\mbi \alpha}^{p}$ is the full flag bundle of $\mc V_{\sharp i_{0}} /\mo_{\hat M_{\mbi \alpha}^{p}}$ over the quotient stack 
$$
\hat M_{\mbi \alpha}^{p} =
\left[ \left( \left( \det \mc V_{\sharp 0} \right) ^{\otimes k (\zeta_{0} - \zeta^{-}_{0})} \otimes \left( \det \mc V_{\sharp 1} \right) ^{\otimes k (\zeta_{1} - \zeta^{-}_{1})} \right) ^{\times} / \C^{\ast}_{u} \right]
$$ 
where $D=( k(\mbi \zeta - \mbi \zeta^{-}), \mbi \alpha)$, and $\C^{\ast}_{u}$ acts by fiber-wise multiplication of $u^{p D}$.
\end{prop}
\proof
It is proven similarly as in 
%\cite[Proposition 5.9]{NY3} for a real root $\mbi \alpha$, and 
\cite[Proposition 6.1]{O}.
%for $\mbi \alpha=\mbi \delta$.
\endproof

The homomorphism $\C^{\ast}_{u} \to \C^{\ast}_{s}$ given by $s=u^{p D}$ induces an \'{e}tale and finite morphism $\hat M_{\mbi \alpha}^{p} \to M_{\mbi \alpha}^{p}$ of degree $1/pD$.

%%%%%%%%%%%%%%%%%%%%%%%%%%%%%%%%%%%%%%%%%%%%%%%%%%%%%%%%%%%%%%%%%%%%%%%%%%%

\subsection{Destabilizing objects for real roots}
\label{subsec:destab}
We describe $M_{\mbi \alpha_{m}}^{p}$ for $m \in \Z$.
But we omit proof since we do not use this description later, and it is proven similarly to \cite[\S 5.4]{NY3}.  
By the symmetry \eqref{shift}, we can reduce to walls on the region defined by 
$$
\zeta_{0} < 0, \text{ or } \zeta_{1} > 0.
$$
So we may assume that $\mathfrak{D}=\mathfrak{D}_{\mbi \alpha_{m}}$.

Then $\mbi \zeta^{\mk D}$-stable ADHM data in $M(\mbi 0, \mbi \alpha_{m})$ is unique object, written by $P^{(m)} = P^{(m)}_{0} \oplus P^{(m)}_{1}$, up to isomorphisms.
Hence vector spaces $P^{(m)}_{0}$ and $P^{(m)}_{1}$ have $\tilde{T}$-module structures.
Every $\mbi \zeta^{\mk D}$-semistable ADHM data in $\mb M(\mbi 0, p \mbi \alpha_{m})$ is isomorphic to $P^{(m)} \otimes \C^{p}$.
An object $(P^{(m)} \otimes \C^{p}, F^{\bullet})$ is $(\mk D, +)$-stable if and only if there exists no proper subspace $S$ of $\C^{p}$ such that $F^{1} \subset P^{(m)}_{i_{0}} \otimes S$.
If we take a generator $z$ of $F^{1}$, this is equivalent to saying that $z \colon (P^{(m)}_{i_{0}})^{\vee} \to \C^{p}$ is surjective via $P^{(m)}_{i_{0}} \otimes \C^{p} \cong \Hom_{\C}((P^{(m)}_{i_{0}})^{\vee}, \C^{p})$.
Hence $M_{\mbi \alpha_{m}}^{p}$ is isomorphic to the Grassmannian $Gr((P^{(m)}_{i_{0}})^{ \vee}, p)$ of surjections $z \in \Hom_{\C}((P^{(m)}_{i_{0}})^{ \vee}, \C^{p}) = V_{\sharp i}$. 

If we write by $\mc Q$ the universal quotient of $(P^{(m)}_{i_{0}})^{ \vee} \otimes \mo_{Gr((P^{(m)}_{i_{0}})^{ \vee}, p)}$ on $Gr((P^{(m)}_{i_{0}})^{ \vee}, p)$, then we have $\mc V_{\sharp} \cong P^{(m)} \otimes \mc Q$ via the above isomorphism $M_{\mbi \alpha_{m}}^{p} \cong Gr((P^{(m)}_{i_{0}})^{ \vee}, p)$.

%%%%%%%%%%%%%%%%%%%%%%%%%%%%%%%%%%%%%%%%%%%%%%%%%%%%%%%%%%%%%%%%%%%%%%%%%%%

\subsection{Moduli stacks of fixed points sets}
\label{subsec:modu3}
By the observation in \S \ref{subsec:dire}, we have a decomposition, 
$$
\mc M^{\C^{\ast}_{\hbar}}= \mc M_{+} \sqcup \mc M_{-} \sqcup \bigsqcup_{\mk J \in \mc D^{\ell}(v_{i_{0}}, \alpha_{i_{0}})} \mc M_{\mk J},
$$
where $\mc M_{\pm}=\lbrace x_{\mp} = 0 \rbrace$, and $\mc M_{\mk J}$ is described as follows.
We fix a direct sum decomposition $V=V_{\flat} \oplus V_{\sharp}$ of a $\Z_{2}$-graded vector space corresponding to decomposition data $\mk J$, where $V_{\flat} = V_{\flat 0} \oplus V_{\flat 1}$ and $V_{\sharp}= V_{\sharp 0} \oplus V_{\sharp 1}$ are also $\Z_{2}$-graded.

We write by $\mc V^{\mk J}_{\flat}, \mc V^{\mk J}_{\sharp}$ $\Z_{2}$-graded tautological bundles on $\mc M_{\mk J}$
% corresponding to $V_{\flat}=\C^{n-p}, V_{\sharp}=\C^{p}$ via the isomorphism \eqref{mkj}, 
such that we have $\mc V|_{\mc M_{\mk J}} = \mc V^{\mk J}_{\flat} \oplus \mc V^{\mk J}_{\sharp}$ for the $\Z_{2}$-graded tautological bundle $\mc V$ on $\mc M$.
We also write by $\mc V_{\flat}, \mc V_{\sharp}$ $\Z_{2}$-graded tautological bundles on $M^{\mc C, \min(I_{\sharp})-1}(\mbi w, \mbi v - p \mbi \alpha), \wt M_{\mbi \alpha}^{p}$ corresponding to $V_{\flat}, V_{\sharp}$, and write by the same letters their pull-backs to $M^{\mc C, \min(I_{\sharp})-1}(\mbi w, \mbi v - p \mbi \alpha) \times \wt M_{\mbi \alpha}^{p}$ by projections.

\begin{thm}
[\protect{\cite[Theorem 5.18]{NY3}}]
\label{decomp}
For $\C^{\ast}_{\hbar}$-action on $\mc M$ defined by \eqref{act0}, we have a decomposition
$$
\mc M^{\C^{\ast}_{\hbar}}=\mc M_{+} \sqcup \mc M_{-} \sqcup \bigsqcup_{\mk J \in \mc D^{\ell}(v_{i_{0}}, \alpha_{i_{0}})}\mc M_{\mk J}
$$
such that the followings hold.\\
(i) We have $\mc M_{+} \cong \wt{M}^{\mc C, \ell}(\mbi w, \mbi v)$ and $\mc M_{-} \cong \wt{M}^{\mc C, 0}(\mbi w, \mbi v)$, that is, the full flag bundle $Fl(\mc V_{i_{0}}, [v_{i_{0}}])$ of the tautological bundle $\mc V_{i_{0}}$ over $M^{\mc C}(\mbi w, \mbi v)$.\\
(ii) For each $\mk J=(I_{\flat}, I_{\sharp}) \in \mc D^{\ell}(v_{i_{0}}, \alpha_{i_{0}})$, we have finite \'etale morphisms 
$$
F \colon \mc S_{\mk J} \to \mc M_{\mk J}, G \colon \mc S_{\mk J} \to \wt M^{\mc C, \min(I_{\sharp})-1}(\mbi w, \mbi v - p \mbi \alpha) \times \wt M_{\mbi \alpha}^{p}
$$ 
of degree $\frac{1}{ p D }$, where $p=|I_{\sharp}| / \alpha_{i_{0}}$. \\
(iii) There exists a line bundle $L_{\mc S_{\mk J}}$ on $\mc S_{\mk J}$ such that we have isomorphisms 
$$
L_{\mc S_{\mk J}}^{\otimes pD} \cong 
\bigotimes_{j =0, 1} G^{\ast}( \det \mc V_{\flat j}  \otimes  \det \mc V_{\sharp j})^{k (\zeta_{j} - \zeta_{j}^{-})},
$$ 
$F^{\ast} \mc V^{\mk J}_{\flat} \cong G^{\ast} \mc V_{\flat}$, and $F^{\ast} \mc V^{\mk J}_{\sharp} \cong G^{\ast} \mc V_{\sharp} \otimes e^{\frac{\hbar}{pD}} \otimes L_{\mc S_{\mk J}}^{\vee}$ as $\C^{\ast}_{\hbar}$-equivariant vector bundles on $\mc S_{\mk J}$.
\end{thm}

\proof
This is similarly proven as in \cite[Theorem 5.18]{NY3}.
\endproof

We also recall that an obstruction theory of enhanced master space $\mc M$ is given by $ob_{\mc M} = \id_{L_{\mc M}}$, where $L_{\mc M}$ is the cotangent complex of $\mc M$.
Then we have induced obstruction theories on $\mc M_{\pm}$ and $\mc M_{\mk J}$ as in \cite[\S 6]{O}.
Furthermore these obstruction theories gives an obstruction theory on $M_{\mbi \alpha}^{p}$, which is different from usual one for $\mbi \alpha=\mbi \delta$.
For this reason, we must consider {\it virtual} fundamental cycles for $M_{\mbi \delta}^{p}$, while we consider usual fundamental cycles for the other moduli spaces.

%%%%%%%%%%%%%%%%%%%%%%%%%%%%%%%%%%%%%%%%%%%%%%%%%%%%%%%%%%%%%%%%%%%%%%%%%%%
%%%%%%%%%%%%%%%%%%%%%%%%%%%%%%%%%%%%%%%%%%%%%%%%%%%%%%%%%%%%%%%%%%%%%%%%%%%
%%%%%%%%%%%%%%%%%%%%%%%%%%%%%%%%%%%%%%%%%%%%%%%%%%%%%%%%%%%%%%%%%%%%%%%%%%%

\section{Proof of Theorem \ref{main}}
\label{sec:Proo}
In this section, we deduce functional equations of Nekrasov partition functions as an application of the previous section.
These are the similar calculations to \cite[\S 6]{NY3}, hence we omit detail explanation.

In the following, we use wall-crossing formula deduced from analysis in the previous section.
For simplicity, we always assume that $\mbi \alpha$ is equal to $\mbi \alpha_{m}$ for $m \ge 0$, or $\mbi \delta$, and $i_{0}=1$, that is, we use full flags of $V_{1}$.
%%%%%%%%%%%%%%%%%%%%%%%%%%%%%%%%%%%%%%%%%%%%%%%%%%%%%%%%%%%%%%%%%%%%%%%%%%%

\subsection{Iterated cohomology classes}
\label{subsec:iter1}
We prepare notation for iteration of wall-crossing formula.
We consider a $\tilde{T}$-manifold $M$ and $\tilde{T}$-equivariant $\Z_{2}$-graded vector bundles $\mc V=\mc V_{0} \oplus \mc V_{1}$ and $\mc W= \mc W_{0} \oplus \mc W_{1}$ on $M$. 
We take $f(\mbi x, \mbi y , t) \in \mb Q (\mbi \e, \mbi a, \mbi m) [\mbi x, \mbi y][[t]]$, where $\mbi x=(x_{1}, \ldots, x_{v_{0}+v_{1}}), \mbi y = (y_{1}, \ldots, y_{r_{0}+r_{1}})$.
We put
$$
f(\mc V, \mc W, t) = f(c (\mc V_{0}), c(\mc V_{1}), c(\mc W_{0}), c (\mc W_{1}), t) \in A_{\tilde{T}}^{\ast} (M) [[t]],
$$
where $c(\mc V_{i}) =(c_{1} (\mc V_{i}), \ldots, c_{v_{i}} (\mc V_{i}))$ and $c(\mc W_{i}) = ( c_{1} (\mc W_{i}), \ldots, c_{r_{i}} (\mc W_{i}))$ are $\tilde{T}$-equivariant Chern classes for $i=0,1$.
For example, we put 
\begin{eqnarray}
\label{matter}
\mc F_{r}(\mc V_{0}) = \bigoplus_{f=1}^{2r} \mc V_{0} \otimes \frac{e^{m_{f}}}{\sqrt{t_{1}t_{2}}},
\end{eqnarray}
and consider the top Chern class $e(\mc F_{r}(\mc V_{0})) \in A^{\ast}_{\tilde{T}}(M)$.
Then we can take $f(x_{1}, \ldots, x_{v_{0}}) = f(\mbi x, \mbi y , t) \in \mb Q (\mbi \e, \mbi a, \mbi m) [\mbi x]$ such that $e(\mc F_{r}(\mc V_{0})) =  f(\mc V_{0} ) = f(\mc V, \mc W, t)$, where we regard $f(x_{1}, \ldots, x_{v_{0}})$ as a constant function with respect to $x_{v_{0}+1}, \ldots, x_{v_{0} + v_{1}}, y_{1}, \ldots, y_{r}$ and $t$.

We write by $\Theta^{rel}$ the pull backs to various moduli stacks of the relative tangent bundle of $[\tilde{\mu}^{-1}(0) / G] \to [\mu^{-1} (0) / G]$.
We put 
\begin{eqnarray*}
\wt f =\wt f (\mc V, \mc W, t)= \frac{ f (\mc V, \mc W, t) \cup e (\Theta^{rel})}{v_{1}!} \in A_{\C^{\ast}_{\hbar}}^{\ast}( \mc M \times_{\tilde{T}} E_{m}) \otimes \mb Q [[t]], 
\end{eqnarray*}
where $E_{m} \to E_{m}/\tilde{T}$ is any approximation of the universal bundle $E \tilde{T} \to B \tilde{T}$ over the classifying space.

For $j > 0 $ and $\vec{p} =(p_{1}, \ldots, p_{j})\in \Z_{>0}^{j}$, we consider a product $M_{\vec{p}} = \prod_{i=1}^{j} M_{\mbi \alpha}^{p_{i}}$, and $M \times M_{\vec{p}}$ for a $\tilde{T}$-equivariant manifold $M$.
We endow $\tilde{T}$-equivariant obstruction theories $ob_{M_{\mbi \alpha}^{p_{i}}}$ on $M_{\mbi \alpha}^{p_{i}}$, and the canonical $\tilde{T}$-equivariant obstruction theory $ob_{M}=\id_{M}$ on a $\tilde{T}$-equivariant smooth manifold $M$.
We pull-back these obstruction theories to $M \times M_{\vec{p}}$ by the projections and take direct sum.  
Then we have an obstruction theory $ob_{M \times M_{\vec{p}}}$ on $M \times M_{\vec{p}}$ and a virtual fundamental cycle $[M \times M_{\vec{p}}]^{vir} \in A^{\tilde{T}}_{\ast}(M \times M_{\vec{p}})$. 
For $\alpha \in A_{\tilde{T}}^{\ast}(M \times M_{\vec{p}})$, we write by $\int_{[M_{\vec{p}}]^{vir}} \alpha  \in A_{\tilde{T}}^{\ast}(M)$ the Poincare dual of the push-forward of $\alpha \cap [M \times M_{\vec{p}}]^{vir}$ by the projection $M \times M_{\vec{p}} \to M$. 
Here $M_{\mbi \alpha}^{p}$ is defined in \eqref{malpha} and appears in wall-crossing formula by Proposition \ref{destab}.
We consider the restriction $\mc V_{\sharp} = \mc V_{\sharp 0} \oplus \mc V_{\sharp 1}$ of $\Z_{2}$-graded tautological bundle  on $M^{\mc C}(\mbi w_{\sharp}, p \mbi \alpha)$ to the zero locus $M_{\mbi \alpha}^{p}$.
We write by $\mc V_{\sharp 0}^{(i)}, \mc V_{\sharp 1}^{(i)}$ the pull-backs to $M \times M_{\vec{p}}$ of $\mc V_{\sharp 0}, \mc V_{\sharp 1}$ on $i$-th component $M_{\mbi \alpha}^{p_{i}}$ of $M_{\vec{p}}$.

For $\Z_{2}$-graded $\tilde{T}$-equivariant vector bundles $\mc W, \mc V$ on $M$, we write by the same letters the pull-backs to the product $M \times M_{\vec{p}}$, and put
\begin{eqnarray}
\label{itcoho}
f_{\vec{p}} (\mc V, \mc W, t)
&=& 
\int_{[M_{\vec{p}}]^{vir}} \Res_{\hbar_{1} = \infty} \cdots \Res_{\hbar_{j}=\infty} \frac{ f \left( \mc V \oplus \bigoplus_{i=1}^{j} \mc V_{\sharp}^{(i)} \otimes e^{\hbar_{i}}, \mc W, t \right) e \left( \bigoplus_{i=1}^{j} (\mc V_{\sharp 1}^{(i)}/ \mo_{M \times M_{\vec{p}}} )^{\vee} \right) }{e\left( \bigoplus_{i=1}^{j} \mathfrak{N} (\mc W, \mc V \oplus \bigoplus_{k=i+1}^{j} \mc V_{\sharp}^{(k)} \otimes e^{\hbar_{k}}, \mc V_{\sharp}^{(i)} \otimes e^{\hbar_{i}})\right)}
\end{eqnarray}
in $A^{\ast}_{\tilde{T}}(M)$.
Here $e^{\hbar}$ is a trivial bundle with $e^{\hbar}$-weight, and $\mathfrak{N}(\mc W', \mc V', \mc V'')$ is defined by
\begin{eqnarray}
\label{normal}
\mathfrak{N}(\mc W', \mc V', \mc V'')
&=&
\mc H om_{\Z_{2}}( Q^{\vee} \otimes \mc V', \mc V'') + \mc H om_{\Z_{2}}( \mc W', \mc V'') \\  
&-&
\mc H om_{\Z_{2}}( \mc V', \mc V'') - \mc H om_{\Z_{2}}( \det Q^{\vee} \otimes \mc V', \mc V'') \notag\\
&+& 
\mc H om_{\Z_{2}}( Q^{\vee} \otimes \mc V'', \mc V' ) + \mc H om_{\Z_{2}}( \det Q^{\vee} \otimes \mc V'', \mc W' )  \notag \\
&-&
\mc H om_{\Z_{2}}( \mc V'', \mc V')  - \mc H om_{\Z_{2}}( \det Q^{\vee} \otimes \mc V'', \mc V') \notag
\end{eqnarray}
for $\Z_{2}$-graded vector bundles $\mc W', \mc V', \mc V''$.

To describe each component $\mc M_{\mk J}$ of $\mc M^{\C^{\ast}_{\hbar}}$, we consider a group action $\C^{\ast}_{\frac{\hbar}{pD}} \times \mc M \to \mc M$ defined by 
\begin{eqnarray}
\label{act}
\left( X,F^{\bullet}, [x_{-},x_{+}] \right) \mapsto \left( \id_{V_{\flat}} \oplus e^{\frac{h}{pD}} \id_{V_{\sharp}} \right) \left(X,F^{\bullet}, [e^{\hbar}x_{-},x_{+}] \right). 
\end{eqnarray}
This action is equal to the original $\C^{\ast}_{\hbar}$-action \eqref{act0}, since the difference is absorbed in $G$-action.
Then $(X_{\flat} \oplus X_{\sharp}, F^{\bullet}_{\flat} \oplus F^{\bullet}_{\sharp}, [x_{-}, x_{+}])$ is fixed by this $\C^{\ast}_{\hbar}$-action, and represents a $\C^{\ast}_{\hbar}$ -fixed point in $\mc M$.
Hence we need to multiply $\mc V_{\sharp}^{(i)}$ with $e^{\hbar_{i}}$ in \eqref{itcoho}.
On the other hand, $e \left( \bigoplus_{i=1}^{j} (\mc V_{\sharp 1}^{(i)}/ \mo_{M \times M_{\vec{p}}} )^{\vee} \right)$ in \eqref{itcoho} is obtained by integrations of $\Theta^{rel}$ over $M^{p_{i}}_{\mbi \alpha}$.
This does not include $e^{\hbar}$ since 
$$
\Theta^{rel} = \sum_{i > j}  \mc Hom\left( \mc F^{j} / \mc F^{j-1}, \mc F^{i} / \mc F^{i-1} \right).
$$

By the projection formula, we have $f_{\vec{p}}(\mbi x, \mbi y, \mbi t) \in \mb Q(\mbi \e, \mbi a, \mbi m)[\mbi x, \mbi y][[t]]$ independent of $M$ since we have finite $\mb T$-fixed points sets of $M_{\vec{p}}$.

%%%%%%%%%%%%%%%%%%%%%%%%%%%%%%%%%%%%%%%%%%%%%%%%%%%%%%%%%%%%%%%%%%%%%%%%%%%

\subsection{Localizations}
\label{subsec:loca}
By the main result \cite[(1)]{GP} and Theorem \ref{decomp}, we have the following diagram
$$
\xymatrix{
\varprojlim_{m} A^{\ast}_{\C^{\ast}_{\hbar}} (\mc M \times_{\tilde{T}} E_{m}) \otimes_{\C[\hbar]} \C[\hbar^{\pm 1}] \ar[d]_{\Pi_{\ast} (\cdot) \cap {[\mc M]^{vir}}} \ar[r]^{\cong}& \varprojlim_{m} A^{\ast} (\mc M^{\C^{\ast}_{\hbar}} \times_{\tilde{T}} E_{m}) \otimes_{\C} \C[\hbar^{\pm 1}] \ar[d]^{\Pi_{\ast} (\cdot) \cap \left( [\mc M_{+}]^{vir} +  [\mc M_{-}]^{vir} + \sum_{\mk J} [\mc M_{\mk J}]^{vir} \right)} \\
\varprojlim_{m} A_{\ast} (M_{0}(\mbi w, \mbi v) \times_{\tilde{T}} E_{m}) \otimes_{\C} \C[\hbar, \hbar^{-1}] \ar@{=}[r] & \varprojlim_{m} A_{\ast} (M_{0}(\mbi w, \mbi v) \times_{\tilde{T}} E_{m}) \otimes_{\C} \C[\hbar,\hbar^{-1}]\\
%\mc S[\hbar, \hbar^{-1}]] \ar@{=}[r]& \mc S[\hbar, \hbar^{-1}]]
}
$$
where the upper horizontal arrow is given by 
$$
\frac{\iota_{+}^{\ast}}{e(\mk N(\mc M_{+}))} + \frac{\iota_{-}^{\ast}}{e(\mk N(\mc M_{-}))} + \sum_{\mk J \in \mc D^{\ell}(v_{1}, \alpha_{1})} \frac{\iota_{\mk J}^{\ast}}{e(\mk N(\mc M_{\mk J}))}.
$$
Here $\hbar$ corresponds to the first Chern class in $A^{\C^{\ast}_{\hbar}}(\text{pt})$ of the weight $e^{\hbar} \in \C^{\ast}_{\hbar}$, and $\iota_{\pm}$ and $\iota_{\mk J}$ are embeddings of $\mc M_{\pm}$ and $\mc M_{\mk J}$ into $\mc M$.

For $j_{0} > 0 $ and $\vec{p}_{0} = (p_{01}, \ldots, p_{0j_{0}}) \in \Z_{>0}^{j_{0}}$, we take equivariant classes $\varphi = \wt f_{\vec{p}_{0}}(\mc V, \mc W, t)$ on $\mc M$.
For the convenience, we also put $f^{()}=f$ for $j_{0}=0$. 
By the above diagram, we have
\begin{eqnarray*}
\Pi_{\ast} \left( [\mc M]^{vir} \cap \varphi \right)
& = &
\Pi_{\ast} \left( \frac{ [\mc M_{+}]^{vir} \cap \iota_{+}^{\ast} \varphi }{e(\mk N(\mc M_{+}))}  + \frac{ [\mc M_{-}]^{vir} \cap \iota_{-}^{\ast} \varphi }{ e(\mk N(\mc M_{-})) } + 
\sum_{\mk J \in \mc D^{\ell}(v_{1}, \alpha_{1})} \frac{ [\mc M_{\mk J}]^{vir} \cap \iota_{\mk J}^{\ast} \varphi }{e(\mk N(\mc M_{\mk J}))} \right).
\end{eqnarray*}
The left hand side is a limit of polynomials in $\hbar$, hence taking coefficients of $\hbar^{-1}$ we have 
\begin{eqnarray}
\label{wc}
\int_{\wt M^{\mc C, \ell}(\mbi w, \mbi v)} \wt f_{\vec{p}_{0}}(\mc V, \mc W, t) - \int_{M^{\mc C}(\mbi w, \mbi v)} f_{\vec{p}_{0}}(\mc V, \mc W, t) = 
\Res_{\hbar =\infty} \sum_{\mk J \in \mc D^{\ell}(v_{1}, \alpha_{1})} \int_{[\mc M_{\mk J}]^{vir}}\frac{ \iota_{\mk J}^{\ast} \varphi }{e(\mk N(\mc M_{\mk J}))}
\end{eqnarray}
by Theorem \ref{decomp} (i) and $e(\mk N(\mc M_{\pm}))=\pm (\hbar - c_{1}(L_{+}^{\vee} \otimes L_{-})))$.
Here $\Res_{\hbar =\infty}$ denotes the operator taking the minus of coefficients of $\hbar^{-1}$.

Furthermore by Theorem \ref{decomp} (ii), Proposition \ref{destab}, the right hand side is equal to
\begin{eqnarray}
\label{exc}
\Res_{\hbar =\infty} \sum_{\mk J \in \mc D^{\ell}(v_{1}, \alpha_{1})} \frac{(v_{1} - p \alpha_{1})! (p \alpha_{1} - 1)! }{v_{1} !} \int_{\wt M^{\mc C, \min(I_{\sharp})-1}(\mbi w, \mbi v - p \mbi \alpha) } \wt f^{ ( \vec{p}_{0}, p )}(\mc V, \mc W, t),
\end{eqnarray}
where $p = |I_{\sharp}|/\alpha_{1}$ is determined from $\mk J=(I_{\flat}, I_{\sharp})$, and $(\vec{p_{0}}, p)=(p_{01}, \ldots, p_{0j_{0}}, p) \in \Z_{>0}^{j_{0} + 1}$.

In this expression, we deleted some line bundles and a parameter $pD$, since we have $\Res_{\hbar = \infty} f(\hbar) = pD \Res_{\hbar = \infty} f( pD \hbar+a)$ (cf. \cite[\S 8.2]{O}), and $\wt M_{\mbi \alpha}^{p}$ are full flag bundles of $1/pD$-degree \'etale covering of $M_{\mbi \alpha}^{p}$ by Proposition \ref{destab}.

%%%%%%%%%%%%%%%%%%%%%%%%%%%%%%%%%%%%%%%%%%%%%%%%%%%%%%%%%%%%%%%%%%%%%%%%%%%

\subsection{Wall-crossing formula}
\label{subsec:iter2}
For $j > 0 $, we put 
$$
S_{j}(v_{1}, \alpha_{1}) = \left \lbrace \vec{p}= (p_{1}, \ldots, p_{j}) \in \Z_{>0} ^{j} \ \Big| \ \sum_{k=1}^{j} p_{k} \alpha_{1}  \le v_{1} \right \rbrace.
$$
For $\vec{p}= (p_{1}, \ldots, p_{j}) \in \Z_{>0} ^{j}$, we put $|\vec{p}| = p_{1} + \cdots + p_{j}$.

\begin{thm}
\label{thm:wcm}
We have
\begin{eqnarray}
\notag
&&
\int_{M^{\mc C'}(\mbi w, \mbi v)}  f(\mc V, \mc W, t)  - \int_{M^{\mc C}(\mbi w, \mbi v)}  f (\mc V, \mc W, t) \\
\label{main>0}
&=&
\sum_{j=1}^{\lfloor \frac{v_{1}}{{\alpha_{1}}} \rfloor} 
\frac{1}{\alpha_{1}^{j}} \sum_{\vec{p} \in S_{j}(v_{1}, \alpha_{1})} \frac{1}{\prod_{i=1}^{j} \sum_{1 \le k \le i} p_{k}} \int_{M^{\mc C}(\mbi w, \mbi v - |\vec{p}| \mbi \alpha )} f_{\vec{p}} (\mc V, \mc W, t),
\end{eqnarray}
where integrands $f_{\vec{p}} (\mc V, \mc W, t)$ are defined in \eqref{itcoho}.
\end{thm}

\proof
Let $\Dec_{j}(v_{1}, \alpha_{1})$ be the set of collections $\mk J=(I_{\flat}, I_{\sharp 1}, \ldots, I_{\sharp j})$ such that
\begin{enumerate}
\item[$\bullet$] $[v_{1}]=I_{\flat} \sqcup I_{\sharp 1} \sqcup \cdots \sqcup I_{\sharp j}$,
\item[$\bullet$] $|I_{\sharp i}| = p_{i} \alpha_{1}$ for $p_{i} \in \Z_{>0}$ ($i=1, \ldots, j$), and
\item[$\bullet$] $\min(I_{\sharp 1}) > \cdots > \min(I_{\sharp j})$.
\end{enumerate}
We note that $\Dec_{1}(v_{1}, \alpha_{1}) = \mc D^{v_{1}}(v_{1}, \alpha_{1})$.
We consider maps $\sigma_{j} \colon \Dec_{j+1}(v_{1}) \to \Dec_{j}(v_{1}, \alpha_{1}), $
$$
\mk J= (I_{\flat}, I_{\sharp 1}, \ldots, I_{\sharp j+1}) \mapsto \sigma_{j}(\mk J) = (I_{\flat} \sqcup I_{\sharp j+1}, I_{\sharp 1}, \ldots, I_{\sharp j}),
$$
and $\rho_{j} \colon \Dec_{j}(v_{1}, \alpha_{1}) \to S_{j}(v_{1}, \alpha_{1}), \mk J \mapsto \vec{p}_{\mk J} = (\frac{|I_{\sharp i}|}{\alpha_{1}}, \ldots, \frac{|I_{\sharp j}|}{\alpha_{1}})$.

\begin{lem}
We have
\begin{eqnarray}
\notag
&&
\int_{M^{\mc C'}(\mbi w, \mbi v)} f(\mc V, \mc W, t)  - \int_{M^{\mc C}(\mbi w, \mbi v )} f(\mc V, \mc W, t) \\
\notag
&=&
\sum_{i=1}^{j-1} \sum_{\mk J \in \Dec_{i}(v_{1}, \alpha_{1})} \frac{ |I_{\flat}|! \prod_{k=1}^{i} (|I_{\sharp k}|-1)! }{ v_{1}! }
\int_{ M^{\mc C}(\mbi w, \mbi v - |\vec{p}_{\mk J}| \mbi \alpha)} f_{\vec{p}_{\mk J}} (\mc V, \mc W, t)\\
\label{formula2>0}
&+&
\sum_{\mk J \in \Dec_{j}(v_{1}, \alpha_{1})} \frac{|I_{\flat}|! \prod_{k=1}^{j} (|I_{\sharp k}|-1)! }{ v_{1}! }
\int_{\widetilde{M}^{\mc C, min(I_{\sharp j})-1}(\mbi w, \mbi v - |\vec{p}_{\mk J}| \mbi \alpha)} \wt f_{\vec{p}_{\mk J}} (\mc V, \mc W, t).
\end{eqnarray}
\end{lem}
\proof
We prove by induction on $j$.
In fact, for $j=1$, \eqref{formula2>0} is nothing but \eqref{wc} and \eqref{exc} for $j_{0}=0$ and $\ell = v_{1}$.
For $j \ge 1$, we assume the formulas \eqref{formula2>0}.
Then again by \eqref{wc} and \eqref{exc}, the last summand for each $\mk J \in \Dec_{j}(v_{1}, \alpha_{1})$ is equal to 
\begin{eqnarray*}
&&
 \frac{|I_{\flat}|! \prod_{k=1}^{j} (|I_{\sharp k}|-1)! }{ v_{1}! }
\left( \int_{M^{\mc C}(\mbi w, \mbi v - |\vec{p}_{\mk J}| \mbi \alpha)}  f_{\vec{p}_{\mk J}} (\mc V, \mc W, t)  \right. \\
&+&
\left. \sum_{\mk J \in \sigma_{j}^{-1} (\mk J)} 
\frac{|I_{\flat}|! (|I_{\sharp j+1}|-1)!}{|I_{\flat} \sqcup I_{\sharp j+1}|!} 
\int_{\widetilde{M}^{\mc C, min (I_{\sharp j+1} )-1 }(\mbi w, \mbi v - |\vec{p}_{\mk J}| \mbi \alpha)} \wt f_{\vec{p}_{\mk J}} (\mc V, \mc W, t)  \right ),
\end{eqnarray*}
where $p_{j+1} = |I_{\sharp j+1}| / \alpha_{1}$, and $(\vec{p}_{\mk J}, p_{j+1})=\vec{p}_{\mk J}$.
Hence we have \eqref{formula2>0} for general $j \ge 1$.
\endproof

For $j > \frac{v_{1}}{\alpha_{1}}$, the set $\Dec_{j}(v_{1}, \alpha_{1})$ is empty.
Thus we have
\begin{eqnarray}
\notag
&&
\int_{M^{\mc C'}(\mbi w, \mbi v)} f (\mc V, \mc W, t)  - \int_{M^{\mc C}(\mbi w, \mbi v )} f (\mc V, \mc W, t) \\
\label{formula3}
&=&
\sum_{j=1}^{\lfloor \frac{v_{1}}{\alpha_{1}} \rfloor} \sum_{\mk J \in \Dec_{j}(v_{1}, \alpha_{1})} \frac{ |I_{\flat}|! \prod_{k=1}^{j} (|I_{\sharp k}|-1)! }{ v_{1}! }
\int_{ M^{\mc C}(\mbi w, \mbi v - |\vec{p}_{\mk J}| \mbi \alpha)} f_{\vec{p}_{\mk J}} (\mc V, \mc W, t).
\end{eqnarray}
We note that each summand in the last sum depends only on $\vec{p}_{\mk J}$ at this stage.

\begin{lem}
For $\vec{p} \in S_{j}(v_{1}, \alpha_{1})$ and $\mk J=(I_{\flat}, I_{\sharp 1}, \ldots, I_{\sharp j}) \in \rho_{j}^{-1}(\vec{p})$, we have 
\begin{eqnarray*}
|\rho_{j}^{-1}(\vec{p})| 
&=&
\frac{1}{\alpha_{1}^{j} \prod_{i=1}^{j} \sum_{1 \le k \le i} p_{k}} 
\frac{ v_{1}! }{ |I_{\flat}|! \prod_{k=1}^{j} (|I_{\sharp k}|-1)! }.
\end{eqnarray*}
\end{lem}
\proof
This follows from \cite[Lemma 6.8]{NY3} since $|I_{\sharp k}| = p_{k} \alpha_{1}$.
\endproof

By this lemma and \eqref{formula3}, we get \eqref{main>0} and complete the proof of Theorem \ref{thm:wcm}.
\endproof

In the rest of the paper, we consider the case where 
$$
f(\mc V, \mc W, t) = f (\mc V_{0}) = e(\mc F_{r}(\mc V_{0})), 
$$
where $\F_{r}(\mc V_{0})$ is defined in \eqref{matter}.
%%%%%%%%%%%%%%%%%%%%%%%%%%%%%%%%%%%%%%%%%%%%%%%%%%%%%%%%%%%%%%%%%%%%%%%%%%%

\subsection{Wall-crossing across $\mathfrak{D}_{\mbi \alpha_{m}}$}
\label{subsec:wall1}
We take a real root $\mbi \alpha_{m}=(m, m+1) \in R_{+}$, and apply \eqref{main>0} to the case where $\mc C=\mc C_{m}, \mc C'=\mc C_{m+1}$, and $\mk D= \mk D_{\mbi \alpha_{m}}$.
\begin{thm}
\label{vanish2}
We put $k = - 2 v_{0} + 2 v_{1} - w_{1} $.
Then we have 
$$
\int_{M^{\mc C_{m+1}}(\mbi w, \mbi v)} e(\mc F_{r}(\mc V_{0})) = \int_{M^{\mc C_{m}}(\mbi w, \mbi v)} e(\mc F_{r}(\mc V_{0})) \text{ for } 
\begin{cases}
m \ge 0 \text{ if } k \le 0,\\
m \le 0 \text{ if }  k \ge 0.
\end{cases}
$$
\end{thm}
\proof
By Theorem \ref{thm:wcm}, it is enough to show that all summands in the right hand of \eqref{main>0} vanish.
For $k \le 0$, it follows from the similar arguments in \cite[Theorem 2.1]{NY3}.
In fact, we have cohomological degrees
\begin{eqnarray*}
\deg \int_{M^{\mc C_{m+1}}(\mbi w, \mbi v)} e(\mc F_{r}(\mc V_{0})) 
&=& \deg \int_{M^{\mc C_{m}}(\mbi w, \mbi v)} e(\mc F_{r}(\mc V_{0}))\\
&=& \dim M^{\mc C_{m}}(\mbi w, \mbi v) - 2r v_{0} \\
&=& 2(w_{0}v_{0}+w_{1}v_{1}) -2 (v_{0} -v_{1})^{2} - 2rv_{0}.
\end{eqnarray*}
On the other hand, we can write summands in the right hand of \eqref{main>0} as
$\int_{M^{\mc C_{m}}(\mbi w, \mbi v - j \mbi \alpha_{m})} e\left(\mc F_{r} (\mc V_{0} ) \right) \cup ?$
for some cohomology class $?$.
Therefore its degree is at most
\begin{eqnarray*}
&&
\dim M^{\mc C_{m}}(\mbi w, \mbi v - j \mbi \alpha_{m}) - 2r (v_{0} - jm)\\ 
&=& 2w_{0}(v_{0} - j m) + 2w_{1}(v_{1} -jm -j) -2 ( v_{0} - v_{1} +j)^{2} - 2r v_{0} + 2r j m\\
&=& 2(w_{0}v_{0}+w_{1}v_{1}) -2 (v_{0} -v_{1})^{2} - 2rv_{0} + 2j ( -2v_{0} + 2 v_{1} - w_{1} - j)\\
&<& \dim M^{\mc C_{m}}(\mbi w, \mbi v) - 2rv_{0}
\end{eqnarray*}
if $k=-2v_{0} + 2v_{1} - w_{1} \le 0$ and $j>0$. 
Hence it is zero.

For the case where $k \ge 0$, we can use isomorphisms \eqref{shift} to reduce the case where $k \le 0$.

\endproof
When $k \le 0$, repeating this theorem, we have $\int_{M^{\mc C_{\infty}}(\mbi w, \mbi v)} e(\mc F_{r}(\mc V_{0})) = \int_{M^{\mc C_{0}}(\mbi w, \mbi v)} e(\mc F_{r}(\mc V_{0})) $, and get a proof of Theorem \ref{main} for $k \le 0$.

%%%%%%%%%%%%%%%%%%%%%%%%%%%%%%%%%%%%%%%%%%%%%%%%%%%%%%%%%%%%%%%%%%%%%%%%%%%

\subsection{Wall-crossing across $\mk D_{\mbi \delta}$}
\label{subsec:wall2}
We assume $k \ge 0$, and apply \eqref{main>0} to the case where $\mk D=\mk D_{\mbi \delta}$, $\mc C = \mc C_{\infty}$ and $\mc C'= - \mc C_{-\infty}$ (cf. \S \ref{subsec:adhm3}).
Here $\mc C_{\infty}$ and $\mc C_{-\infty}$ are defined in \S \ref{subsec:cham}. 

We take $\mbi \zeta=\mbi \zeta^{1} \in \mc C=\mc C_{\infty}$ and $\mbi \zeta'=-\mbi \zeta^{1}[1] \in \mc C'= - \mc C_{-\infty}$ as in Figure \ref{zeta+}.
For $n = v_{0} + \frac{w_{1}}{4}$, we put 
$$
\alpha_{n}=\int_{M^{\mbi \zeta}(\mbi w, \mbi v)} e(\mc F_{r}(\mc V_{0})), \beta_{n}=\int_{M^{\mbi \zeta'}(\mbi w, \mbi v)} e(\mc F_{r}(\mc V_{0}))
$$ 
By \eqref{ADHMpartition}, we have $Z_{X_{1}}^{k}(\mbi \e, \mbi a, \mbi m, q)= \sum \alpha_{n} q^{n}$.

On the other hand, by Theorem \ref{vanish2}, the integrations over $M^{\mbi \zeta^{0}}(\mbi w, \mbi v)$ for $\mbi \zeta^{0} \in \mc C_{0}$ and $M^{\mbi \zeta^{1}[1]}(\mbi w, \mbi v)$ are same.
Furthermore, since we have a $\tilde{T}$-equivariant isomorphism $M^{\mbi \zeta^{1}[1]}(\mbi w, \mbi v) \cong M^{-\mbi \zeta^{1}[1]}(\mbi w, \mbi v)$ by \eqref{dual} modulo the automomorphism of $G \times \tilde{T}$, we can show that $Z_{X_{0}}^{k}(\mbi \e, -\mbi a, -\mbi m, q) = \sum \beta_{n} q^{n}$
%$$
%Z_{X_{0}}^{k}(\mbi \e, -\mbi a, -\mbi m, q)= \sum_{v_0 \in \Z} q^{v_0 + \frac{w_1}{4}} \int_{M^{ - \mc C_{-}}(\mbi w, \mbi v)} e(\mc F_{r}(\mc V_{0}))
%$$
similarly as in \cite[\S 3.2]{O}.

From the combinatorial description in \S \ref{sec:comb}, we can show $Z_{X_{\kappa}}^{k}( - \mbi \e, -\mbi a, -\mbi m, q) = Z_{X_{\kappa}}^{k}( \mbi \e, \mbi a, \mbi m, q)$ as in \cite[Lemma 6.3 (3)]{NY1}.
Hence to prove Theorem \ref{main} for $k \ge 0$, we must show
\begin{eqnarray}
\label{must}
\alpha_{n} = \sum_{k=0}^{v_{0}} \frac{(-1)^{kr} u_{r} (u_{r}+1) \cdots (u_{r} + k -1)}{k!} \beta_{n-k},
\end{eqnarray}
which will be proved in the rest of this section by the similar argument as in \cite{O}.

We compute $f_{\vec{p}}(\mc V, \mc W, t) = f_{\vec{p}}(\mc V_{0})$, and substitute it into \eqref{main>0}.
This is reduced to integrations over $M^{p}_{\mbi \delta} = M^{\mc C_{\infty}}(\mbi w_{\sharp}, p \mbi \delta)$ as in \S \ref{subsec:iter1}, where $\mbi w_{\sharp} = (0,1)$ as defined in \eqref{sharp}.
The $\tilde{T}$-fixed points set $M^{\mc C_{\infty}}(\mbi w_{\sharp}, p \mbi \delta)^{\tilde{T}}$ has a bijection to the set of pairs $\vec{Y} = (Y^{1}, Y^{2})$ of Young diagrams such that $| \vec{Y}| = |Y^{1}| + |Y^{2}| = p$ as in Proposition \ref{fixedpt1} in Appendix.
We write by $I_{\vec{Y}}( - \frac{1}{2} C)$ the corresponding element of $M^{\mc C_{\infty}}(\mbi w_{\sharp}, p \mbi \delta)^{\tilde{T}}$.
We also consider the embedding 
$$
\iota_{\vec{Y}} \colon \wt M^{\mc C_{\infty},\min(I_{\sharp}) -1}(\mbi w, \mbi v - p \mbi \delta) \times \left \lbrace I_{\vec{Y}}( - \frac{1}{2} C) \right \rbrace \to \wt M^{\mc C_{\infty},\min(I_{\sharp}) -1}(\mbi w, \mbi v - p \mbi \delta) \times M^{\mc C_{\infty}}(\mbi w_{\sharp}, p \mbi \delta).
$$

\begin{prop}
\label{hilb}
For $\vec{p}=(p_{1}, \ldots, p_{j}) \in \Z_{>0}^{j}$, we have
\begin{eqnarray*}
f_{\vec{p}}( \mc V_{0})&=&
e(\mc F_{r}(\mc V_{0})) \cup (-1)^{r|\vec{p}|+j}  u_{r}^{j}
\end{eqnarray*}
where $|\vec{p}|=p_{1} + \cdots + p_{j}$, and 
$$
u_{r}=
\frac{(\e_{1}+\e_{2})\left( 2 \sum_{\alpha=1} ^{r} a_{k} + \sum_{f=1} ^{2r} m_{f} \right)}{2\e_{1}\e_{2}}.
$$
\end{prop}
\proof
It is enough to show for $j=1$ and $\vec{p}=p \in \Z_{>0}$, since we can similarly prove for general $j>0$.
We have $f_{p}(\mc V_{0}) = f(\mc V_{0}) \cup \Psi_{p}=e(\mc F_{r}(\mc V_{0})) \cup \Psi_{p}$, where
$$
\Psi_{p}=\sum_{|\vec{Y}|=p} \Res_{\hbar = \infty } \frac{\iota_{\vec{Y}}^{\ast} e(\F_{r}(\mc V_{\sharp 0} \otimes e^{\hbar}))}{ \iota_{\vec{Y}}^{\ast} e( \mk N(\mc W, \mc V, \mc V_{\sharp} \otimes e^{\hbar})) }\cdot \frac{\iota_{\vec{Y}}^{\ast}e(\mc V_{\sharp 1} / \mo_{M^{\mc C_{\infty}}(\mbi w_{\sharp}, p \mbi \delta)})}{\iota_{\vec{Y}}^{\ast}e(\mk{N}^{\sharp}_{\vec{Y}})}.
$$
Here $\mk{N}^{\sharp}_{\vec{Y}} = T_{I_{\vec{Y}}}M^{\mc C_{\infty}}(\mbi w_{\sharp}, p \mbi \delta) / \Hom_{\Z_{2}}(\det Q ^{\vee} \otimes V_{\sharp}, W_{\sharp})$ is the virtual normal bundle induced from the unusual obstruction theories on $M_{\mbi \delta}^{p}$ defined similarly as in \cite[\S 7.2]{O}, and $T_{I_{\vec{Y}}}M^{\mc C_{\infty}}(\mbi w_{\sharp}, p \mbi \delta)$ is the tangent space at $I_{\vec{Y}}$ in $M^{\mc C_{\infty}}(\mbi w_{\sharp}, p \mbi \delta)$.
Since $\det Q  = t_{1}t_{2}$ as a $\tilde{T}$-module, we have
\begin{eqnarray*}
\frac{\iota_{\vec{Y}}^{\ast}e(\mc V_{\sharp 1} / \mo_{M^{\mc C_{\infty}}(\mbi w_{\sharp}, p \mbi \delta) } ) }{\iota_{\vec{Y}}^{\ast}e( \mk{N}^{\sharp}_{\vec{Y}} )}
&=&
\frac{\iota_{\vec{Y}}^{\ast}e(\mc V_{\sharp 1} / \mo_{M^{\mc C_{\infty}}(\mbi w_{\sharp}, p \mbi \delta) } \oplus \mc V_{\sharp 1}^{\vee} \otimes t_{1}t_{2} )}{\iota_{\vec{Y}}^{\ast}e(TM^{\mc C_{\infty}}(\mbi w_{\sharp}, p \mbi \delta)) }.
\end{eqnarray*}
Here we may only take pairs of Young diagrams $\vec{Y}=(Y^{1}, Y^{2})$ such that either one of $Y^{1}$ or $Y^{2}$ is the empty set.
For if neither one is not empty, we see  by Proposition \ref{tautological1} that there is a two dimensional trivial $\tilde{T}$-submodule in the fiber of $\mc V_{\sharp}$ over the fixed point corresponding to $\vec{Y}$, and the Euler class $e(\mc V_{\sharp} / \mo_{M^{\mc C_{\infty}}(\mbi w_{\sharp}, p \mbi \delta)})$ vanishes there. 

On the other hand, we have
\begin{multline*}
\Res_{\hbar = \infty } \frac{\iota_{\vec{Y}}^{\ast} e(\F_{r}(\mc V_{\sharp 0} \otimes e^{\hbar}))}{ \iota_{\vec{Y}}^{\ast} e( \mk N(\mc W, \mc V, \mc V_{\sharp} \otimes e^{\hbar})) } =
p\left( 2 \sum_{\alpha=1} ^{r} a_{k} + \sum_{f=1} ^{2r} m_{f} \right) + \\
4 ( \rank \mc V_{1} - \rank \mc V_{0} - w_{1}) \iota_{\vec{Y}}^{\ast} (c_{1} (\det \mc V_{\sharp 1}) - c_{1} (\det \mc V_{\sharp 0})) 
\end{multline*}
since $\rank \mc V_{\sharp 1} - \rank \mc V_{\sharp 0} = p-p=0$.
By Proposition \ref{tautological1}, we have 
$$
\iota_{\vec{Y}}^{\ast} (c_{1}(\det \mc V_{\sharp 1}) - c_{1}(\det \mc V_{\sharp 0})) =
\begin{cases} 
p \e_{1} & \text{ if } \vec{Y} = (Y, \emptyset) \\
p \e_{2} & \text{ if } \vec{Y} = (\emptyset, Y). 
\end{cases}
$$
Furthermore by \cite[Proposition 8.1]{O} and Proposition \ref{tautological1} and \ref{tangent1}, we have
\begin{eqnarray*}
\sum_{|Y^{1}|+|Y^{2}|=p, Y^{j}=\emptyset} \frac{\iota_{\vec{Y}}^{\ast}e(\mc V_{\sharp 1} / \mo_{M^{\mc C_{\infty}}(\mbi w_{\sharp}, p \mbi \delta) } \oplus \mc V_{\sharp 1}^{\vee} \otimes t_{1}t_{2} )}{\iota_{\vec{Y}}^{\ast}e(TM^{\mc C_{\infty}}(\mbi w_{\sharp}, p \mbi \delta)) } 
&=&
\begin{cases}
\frac{\e_{+}}{2p\e_{1}(\e_{2}-\e_{1})} & \text{ if } j=2\\
\frac{\e_{+}}{2p\e_{2}(\e_{1}-\e_{2})} & \text{ if } j=1.
\end{cases}
\end{eqnarray*}
Combining these together, we get the assertion.
\endproof

By \eqref{main>0} and Proposition \ref{hilb}, we get \eqref{must}, and complete the proof of Theorem \ref{main}.

%%%%%%%%%%%%%%%%%%%%%%%%%%%%%%%%%%%%%%%%%%%%%%%%%%%%%%%%%%%%%%%%%%%%%%%%%%%
%%%%%%%%%%%%%%%%%%%%%%%%%%%%%%%%%%%%%%%%%%%%%%%%%%%%%%%%%%%%%%%%%%%%%%%%%%%
%%%%%%%%%%%%%%%%%%%%%%%%%%%%%%%%%%%%%%%%%%%%%%%%%%%%%%%%%%%%%%%%%%%%%%%%%%%

\appendix

\section{Construction of framed moduli on $X_{\kappa}$}
\label{sec:const}
We show that moduli of ADHM data gives framed moduli on $X_{\kappa}$ for $\kappa=0,1$.
This is just a slight modification of the proof of \cite[Theorem 2.2]{N3} to the relative setting.

%%%%%%%%%%%%%%%%%%%%%%%%%%%%%%%%%%%%%%%%%%%%%%%%%%%%%%%%%%%%%%%%%%%%%%%%%%%
\subsection{Beilinson complex}
\label{subsec:beil1}
We consider tautological bundles $\mc R_{0}=\mo_{X_{\kappa}}$ and $\mc R_{1}=\mo_{X_{\kappa}}(F - \ell_{\infty})$, and put $\mc R= \mc R_{0} \oplus \mc R_{1}$ as in \S \ref{sub:two}.
For $i \in \Z/2 \Z$, we define tautological homomorphisms $\xi_{i} \colon \mc R_{i} \to \mc R_{i+1}( \ell_{\infty})$ and $\bar{\xi}_{i} \colon \mc R_{i+1} \to \mc R_{i}(\ell_{\infty})$ so that $\bar{\xi}_{i} \xi_{i} = \xi_{i-1} \bar{\xi}_{i-1}$ on both $X_{0}$ and $X_{1}$ as follows.
On $X_{0}$, we define $\xi_{i}$ and $\bar{\xi}_{i}$ by the multiplication of $x_{1}$ and $x_{2}$ respectively, where $[x_{0}, x_{1}, x_{2}]$ is the homogeneous coordinate of $\PP^{2}$. 
On $X_{1}$, we define $\xi_{0}, \xi_{1}, \bar{\xi}_{0}$ and $\bar{\xi}_{1}$ by the multiplication of $x_{1}, y x_{1}, yx_{2}$ and $ x_{2}$ respectively, where we use description of $X_{1}$ in \S \ref{subsec:toru}.
By the construction \eqref{cpt1} and \eqref{cpt2}, these define homomorphisms $\xi_{i}, \bar{\xi}_{i}$ on both $X_{0}$ and $X_{1}$, and satisfy $\bar{\xi}_{i} \xi_{i} = \xi_{i-1} \bar{\xi}_{i-1}$.
We can also check easily that these homomorphisms coincide on $X_{0} \setminus O \cong X_{1} \setminus C$. 
We put 
$$
\Xi =
e \otimes 
\begin{bmatrix}
0 & \xi_{1} \\
\xi_{0} & 0
\end{bmatrix}
+
\ba e \otimes
\begin{bmatrix}
0 & \bar{\xi}_{0} \\
\bar{\xi}_{1} & 0
\end{bmatrix}
\in \Hom_{\Z_{2}}(Q^{\vee} \otimes \mc R, \mc R(\ell_{\infty})),
$$
where $e, \ba e$ is a basis of $Q$.

We take $\mbi \zeta =(\zeta_{0}, \zeta_{1}) \in \mathbb{R}^{2}$ such that
\begin{eqnarray}
\label{zeta}
\zeta_{0} + \zeta_{1} < 0, \zeta_{0} < 0
\end{eqnarray}
and does not lie on any wall.
For a moduli $M^{\mbi \zeta}(\mbi w, \mbi v)$ of ADHM data on $\Z/2\Z$-graded vector spaces $W, V$, we also have tautological bundle $\mc W_{i}, \mc V_{i}$ for $i \in \Z/ 2\Z$ corresponding to $W_{i}, V_{i}$, and put $\mc V = \mc V_{0} \oplus \mc V_{1}, \mc W=\mc W_{0} \oplus \mc W_{1}$.
We write by $B \colon Q^{\vee} \otimes \mc V \to \mc V, z \colon \mc W \to \mc V$ and $w \colon \det Q \otimes \mc V \to \mc W$ the tautological homomorphism corresponding to $B, z$ and $w$ in Definition \ref{defn:adhm}.

We consider the following complex on $X_{\kappa} \times M^{\mbi \zeta}(\mbi w, \mbi v)$:
\begin{eqnarray}
\label{bei}
\mc Hom_{\Z_{2}}(\det Q \otimes p_{1}^{\ast} \mc R^{\vee} ( \ell_{\infty}), p_{2}^{\ast} \mc V ) 
%\bigoplus_{i \in \Z/ 2\Z} \mc R_{i}^{\vee} (- \ell_{\infty})\boxtimes \mc V_{i} 
\stackrel{\sigma}{\to} 
\begin{matrix}
\mc Hom_{\Z_{2}}(Q \otimes p_{1}^{\ast} \mc R^{\vee}, p_{2}^{\ast} \mc V)\\
\oplus \\
\mc Hom_{\Z_{2}}(p_{1}^{\ast} \mc R^{\vee}, p_{2}^{\ast} \mc W)
\end{matrix}
\stackrel{\tau}{\to}
\mc Hom_{\Z_{2}}( p_{1}^{\ast} \mc R^{\vee}(-\ell_{\infty}), p_{2}^{\ast} \mc V),
\end{eqnarray}
where $p_{1}$ and $p_{2}$ are projections to the first and second components respectively.
We call this {\it Beilinson complex}.
The differentials $\sigma$ and $\tau$ are defined by
$$
\sigma(\eta) =
\begin{bmatrix}
B \eta x_{0} - \eta \Xi^{\vee} \\
w \eta x_{0}
\end{bmatrix},
\tau(\eta', \gamma) =
\begin{bmatrix}
B \eta' x_{0} - \eta' \Xi^{\vee} \\
z \gamma x_{0} 
\end{bmatrix}.
$$
Here $\eta \in \Hom_{\Z_{2}}(\det Q \otimes p_{1}^{\ast} \mc R^{\vee}(\ell_{\infty}), p_{2}^{\ast} \mc V), \eta' \in \Hom_{\Z_{2}}(Q \otimes p_{1}^{\ast} \mc R^{\vee}, p_{2}^{\ast} \mc V)$, $\gamma \in \Hom_{\Z_{2}}(p_{1}^{\ast} \mc R^{\vee}, p_{2}^{\ast} \mc W)$, and $\Xi^{\vee}$ is the dual of $\Xi$, and we regard $x_{0}$ as a homomorphism among suitable line bundles.
We write by the same symbol $B$ the composition of the natural injection $\det Q \otimes \mc V \to Q \otimes Q \otimes \mc V$ and $\id_{Q} \otimes B \colon Q \otimes \mc V \to Q \otimes Q \otimes \mc V$.
We use same notation for $\Xi^{\vee}$.
Then this is just a relative version of \eqref{coh}, since $Q=Q_{0} \oplus Q_{1}$ is a $\Z/ 2 \Z$-graded vector space with $Q_{0}=0$ and $Q_{1}=\C^{2}$.

By the condition \eqref{zeta}, we can show that the restriction of $\sigma$ and $\tau$ to 
$X_{1} \times \lbrace m \rbrace$ for any closed point 
$m \in M^{\mbi \zeta}(\mbi w, \mbi v)$ is injective and 
surjective respectively as a sheaf homomorphism by the similar argument as in \cite{NY1}.
Hence, using \cite[Theorem 22.5]{Ma}, we see that $\E = \ker \tau / \im \sigma$ is flat over $M^{\mbi \zeta}(\mbi w, \mbi v)$.
This gives a family of framed sheaves $(\E,  \Phi)$ on $X_{\kappa} \times M^{\mbi \zeta^{\kappa}}(\mbi w, \mbi v)$, where $\Phi$ is naturally given by restricting the Beilinson complex \eqref{bei} to $\ell_{\infty} \times M^{\mbi \zeta}(\mbi w, \mbi v)$.
We call $(\E, \Phi)$ a universal framed sheaf, and will check that this is true in the rest of this section.

We take $\mbi \zeta^{0} \in \mc C_{0}$ and $\mbi \zeta^{1} \in \mc C_{+}$.
Then we have a morphism from $M^{\mbi \zeta^{\kappa}}(\mbi w, \mbi v)$ to framed moduli $M_{X_{\kappa}}(c)$ on $X_{\kappa}$ for $\kappa = 0, 1$, where $c \in A(IX_{\kappa})$ is defined by \eqref{chern1} and \eqref{chern2}.
That is, if we have a family of ADHM data on a scheme $S$, we have a morphism  $S \to M^{\mbi \zeta^{\kappa}}(\mbi w, \mbi v)$, and the pull-back $(\E_{S}, \Phi_{S})$ of universal framed sheaf $(\E, \Phi)$ to $X_{\kappa} \times S$ is a flat family over $S$ of framed sheaves on $X_{\kappa}$. 
This gives a morphism $M^{\mbi \zeta^{\kappa}}(\mbi w, \mbi v) \to M_{X_{\kappa}}(c)$, and by \cite{N3}, this is bijection at least set theoretically.

%%%%%%%%%%%%%%%%%%%%%%%%%%%%%%%%%%%%%%%%%%%%%%%%%%%%%%%%%%%%%%%%%%%%%%%%%%%

\subsection{Resolution of diagonal on $X_{\kappa} \times X_{\kappa}$}
\label{subsec:reso}
To construct the converse $M_{X_{\kappa}}(c) \to M^{\mbi \zeta^{\kappa}}(\mbi w, \mbi v)$, we consider resolutions of diagonals $\Delta_{\ast} \mo_{X_{\kappa}}$ on $X_{\kappa} \times X_{\kappa}$, where $\Delta \colon X_{\kappa}^{\circ} \to X_{\kappa}^{\circ} \times X_{\kappa}^{\circ}$ is the diagonal embedding.

First we construct resolutions on $X_{\kappa}^{\circ} \times X_{\kappa}^{\circ}$.
We regard $X_{0}^{\circ} \times X_{0}^{\circ}$ as $[Q \times Q / H_{1} \times H_{2}]$, where $H_{1}$ and $H_{2}$ are copies of $H=\lbrace \pm \id_{Q} \rbrace$ and acts on the first and second components respectively.
Then $\Delta_{\ast} \mo_{X_{0}}$ can be identified with $\mo_{Q} \oplus \mo_{-Q}$, where 
$$
-Q=\lbrace ((x_{1}, x_{2}), -(x_{1}, x_{2})) \in Q \times Q \mid (x_{1}, x_{2}) \in Q=\C^{2} \rbrace.
$$
We consider the following complex on $X_{\kappa}^{\circ} \times X_{\kappa}^{\circ}$ :
\begin{multline}
\mc Hom_{\Z_{2}}(p_{1}^{\ast} \mc R^{\vee}, p_{2}^{\ast} \mc R^{\vee} ) |_{X_{0}^{\circ} \times X_{0}^{\circ}} 
\stackrel{d_{-2}}{\to} 
\mc Hom_{\Z_{2}}(Q^{\vee} \otimes p_{1}^{\ast} \mc R^{\vee}, p_{2}^{\ast} \mc R^{\vee})  |_{X_{0}^{\circ} \times X_{0}^{\circ}} \\
\stackrel{d_{-1}}{\to}
\mc Hom_{\Z_{2}}(\det Q^{\vee} \otimes p_{1}^{\ast} \mc R^{\vee}, p_{2}^{\ast} \mc R^{\vee}) |_{X_{0}^{\circ} \times X_{0}^{\circ}} 
\stackrel{d^{0}}{\to} 
\Delta_{\ast} \mc O_{X_{\kappa}^{\circ}}.
\label{diag}
\end{multline}
Here $d^{-2}, d^{-1}$ is defined by replacing $B$ in first components of $\sigma, \tau$ in \eqref{bei} with $\Xi^{\vee}$.
The last differential $d_{0}$ is defined by the restriction to the diagonal and taking the contraction, where $p_{2}^{\ast} \mc R|_{-Q}$ on $-Q \cong Q$ is identified with $\mc R$ by multiplication of $-1$.
This gives a resolution of the diagonal $\Delta_{\ast}\mo_{X_{\kappa}}$ by \cite[Lemma 4.10]{N2} and \cite[3(iii)]{N3}, and they are identified on $X_{0} \setminus O \cong X_{1} \setminus C$.
This complex without the last component can be viewed as a special case of \eqref{bei} where $W=0$.

But when $W=0$, we need to change definition of stability.
We take $\mbi \zeta$ such that $(\mbi \zeta, \mbi v)=0$.
\begin{defn}
\label{ADHM-stability2}
We say that ADHM data $(B, 0, 0)$ on $W=0, V$ are $\mbi \zeta$-semistable if the following conditions hold:
For any $\Z_{2}$-graded subspace $S$ of $V$, if $B_{i}(S) \subset S$, then we have $\mbi \zeta(S) \le 0$. 
They are said to be $\mbi \zeta$-stable when the strict inequality always holds for non-trivial proper subspace $S$.
\end{defn}
For $\mbi \zeta = \zeta (1, -1)$ with $\zeta > 0$ and $\mbi \delta=(1,1)$, we have an isomorphism $M^{\mbi \zeta}(\mbi 0, \mbi \delta) \cong X_{1}^{\circ}=X_{1} \setminus \ell_{\infty}$ such that tautological bundles $\mc V_{0}, \mc V_{1}$ coincide with $\mc R_{0}, \mc R_{1}$, and $B = \xi$.
Hence $d_{-1}$ and $d_{-2}$ in \eqref{diag} are naturally obtained from \eqref{bei}.

To extend the complex \eqref{diag} on $X_{\kappa}^{\circ}$ to the compactification $X_{\kappa}$, we identify $X_{0} \times X_{0}$ with $[\PP^{2} \times \PP^{2} / H_{1} \times H_{2}]$ as above, where $H_{1}$ and $H_{2}$ are copies of $H$ and acts on the first and second components respectively.
We construct a resolution of the diagonal on $X_{0} \times X_{0}$ as complexes of $H_{1} \times H_{2}$-equivariant vector bundles on $\PP^{2} \times \PP^{2}$.
We construct vector bundles $\mc Q$ on $X_{0}$ and $X_{1}$ by
$$
\mc Q = \coker ( 
x_{0} \oplus \xi_{0} \oplus \bar{\xi}_{1} \colon \mo_{X_{\kappa}}(- \ell_{\infty}) \to \mo_{X_{\kappa}} \oplus \mc R_{1} \oplus \mc R_{1}).
$$

On $X_{0} \times X_{0}$, we consider a map $p_{2}^{\ast} \mo(-\ell_{\infty}) \to \C[H] \otimes p_{1}^{\ast}\mc Q$, defined by $x_{2} \mapsto \sum_{\gamma \in H} \gamma \otimes \varphi(\gamma x_{2})$, where $\varphi$ is defined by the compositions
$$
p_{2}^{\ast} \mo_{\PP^{2}}(-1)  \to p^{\ast}_{1} \mo_{\PP^{2}} \otimes (\C \oplus Q) \to p_{1}^{\ast} \mo_{\PP^{2}} \otimes (\C \oplus Q) / p_{2}^{\ast} \mo_{\PP^{2}}(-1)
$$
on $\PP^{2} \times \PP^{2}$.
Here $(\gamma_{1}, \gamma_{2} )\in H_{1} \times H_{2}$ acts on $\C[H]$ by $\gamma_{1} \gamma \gamma_{2}^{-1}$.
Then this map is $H_{1} \times H_{2}$-equivariant.

Hence we have a section $s \in H(X_{0} \times X_{0},  \C[H] \otimes \mo_{X_{0}}(\ell_{\infty}) \boxtimes \mc Q)$, whose zero locus is equal to $\Delta_{\ast} \mo_{X_{0}}=\mo_{\PP^{2}} \oplus \mo_{- \PP^{2}}$, and the Koszul resolution
\begin{eqnarray}
\label{koszul0}
0 \to C_{0}^{-2} \to C_{0}^{-1} \to C_{0}^{0} \to \Delta_{\ast} \mo_{X_{0}}=\mo_{\PP^{2}} \oplus \mo_{-\PP^{2}}  \to 0.
\end{eqnarray}
Here $- \PP^{2} =\lbrace ( [x_{0}, x_{1}, x_{2}], [x_{0}, -x_{1}, -x_{2}] ) \in \PP^{2} \times \PP^{2}\mid [x_{0}, x_{1}, x_{2}] \in \PP^{2} \rbrace$, and 
$$
C_{0}^{0} = \C[H] \otimes \mo_{X_{0} \times X_{0}}, C_{0}^{-1} = \C[H] \otimes \mc Q^{\vee} \boxtimes \mo_{X_{0}}(-\ell_{\infty}), C_{0}^{-2} = \C[H] \otimes \det \mc Q^{\vee} \boxtimes \mo_{X_{0}}(-2 \ell_{\infty}).
$$
We have $\mc Q|_{X_{0}^{\circ}} \cong \mo_{\C^{2}} \otimes Q$, and $\C[H] \cong \bigoplus_{i \in \Z/2\Z} \Hom(R_{i}^{\vee}, R_{i}^{\vee})$, where  $R_{0}$ and $R_{1}$ are trivial and non-trivial irreducible $H$-representations.
We can identify $\mc R_{0}$ and $\mc R_{1}$ with $H$-equivariant line bundles $\mo_{\PP^{2}} \otimes R_{0}$ and $\mo_{\PP^{2}} \otimes R_{1}$.
Then we have isomorphisms
$$
C_{0}^{-2} |_{X_{0}^{\circ} \times X_{0}^{\circ}} \cong \mc Hom_{\Z_{2}}(p_{1}^{\ast} \mc R^{\vee}, p_{2}^{\ast} \mc R^{\vee} )|_{X_{0}^{\circ} \times X_{0}^{\circ}} , C_{0}^{-1}|_{X_{0}^{\circ} \times X_{0}^{\circ}}  \cong \mc Hom_{\Z_{2}}(Q^{\vee} \otimes p_{1}^{\ast} \mc R^{\vee}, p_{2}^{\ast} \mc R^{\vee})  |_{X_{0}^{\circ} \times X_{0}^{\circ}},
$$
and $C_{0}^{0} |_{X_{0}^{\circ} \times X_{0}^{\circ}}  \cong \mc Hom_{\Z_{2}}(\det Q^{\vee} \otimes p_{1}^{\ast} \mc R^{\vee}, p_{2}^{\ast} \mc R^{\vee}) |_{X_{0}^{\circ} \times X_{0}^{\circ}} $.
Via these isomorphisms, we can check that the complex \eqref{koszul0} coincides with \eqref{diag}.

Hence we can patch the restriction of \eqref{koszul0} to $(X_{0} \setminus O) \times (X_{0} \setminus O)$ and \eqref{diag} on $X_{1}^{\circ} \times X_{1}^{\circ}$ to get the resolution of the diagonal on $X_{1} \times X_{1}$:
\begin{eqnarray}
\label{koszul1}
0 \to C_{1}^{-2} \to C_{1}^{-1} \to C_{1}^{0} \to \Delta_{\ast} \mo_{X_{1}} \to 0.
\end{eqnarray}
Here we can write $C_{\kappa}^{-2} = \bigoplus_{i \in \Z / 2 \Z} \mc R_{i}^{\vee} (-2 \ell_{\infty}) \boxtimes \mc R_{i} \otimes \det \mc Q^{\vee}, C_{\kappa}^{-1}= \bigoplus_{i \in \Z / 2 \Z} \mc R_{i}^{\vee} (-\ell_{\infty}) \boxtimes \mc R_{i} \otimes \mc Q^{\vee}$, and $C_{\kappa}^{0}= \bigoplus_{i \in \Z / 2 \Z}  \mc R_{i} \boxtimes \mc R_{i}$ for $\kappa = 0, 1$.

%%%%%%%%%%%%%%%%%%%%%%%%%%%%%%%%%%%%%%%%%%%%%%%%%%%%%%%%%%%%%%%%%%%%%%%%%%%

\subsection{Beilinson spectral sequence}
\label{subsec:beil2}
For a finitely generated $\C$-algebra $A$, we put $S= \Spec A$.
Here we consider families of framed sheaves $(\E_{S}, \Phi_{S})$ on $X_{\kappa} \times S$, where $\E_{S}$ is a torsion free sheaf on $X_{\kappa} \times S$ flat over $S$, and $\Phi_{S} \colon \E_{S}|_{\ell_{\infty} \times S} \cong W \otimes \mo_{\PP^{1} \times S}$ is an isomorphism.
Pulling back resolutions \eqref{koszul0} and \eqref{koszul1}, we get resolutions of diagonals on $X_{\kappa} \times X_{\kappa} \times S$.

We construct Beilinson complexes from framed sheaves $(\E_{S}, \Phi_{S})$.
We consider 
$$
\mb R p_{1\ast} ( p_{2}^{\ast} \E_{S} (- \ell_{\infty}) \otimes \Delta_{\ast} \mo_{X_{\kappa}} )= \mb R p_{1\ast} ( p_{2}^{\ast} \E_{S} (- \ell_{\infty}) \otimes C_{\kappa}^{\bullet})
$$ 
as a double complex, where $p_{1}, p_{2} \colon X_{\kappa} \times X_{\kappa} \times S \to X_{\kappa} \times S$ are the first and second projections, and $C_{\kappa}^{\bullet}$ is the complex in \eqref{koszul0} and \eqref{koszul1}.
We have a spectral sequence associated to this double complex, whose $E_{2}$-term is given by
$$
E_{2}^{pq}= \mb R^{q} p_{1\ast} (p_{2}^{\ast} \E_{S}(-\ell_{\infty}) \otimes C_{\kappa}^{p}).
$$
Explicitly, we have 
$$
E_{2}^{p q} = 
\begin{cases}
\bigoplus_{i \in \Z / 2\Z} \mc R_{i}^{\vee} (-2 \ell_{\infty}) \boxtimes \mb R^{q} p_{S \ast} ( \E_{S}( -2 \ell_{\infty}) \otimes \mc R_{i} ) & \text{ for } p=-2,\\
\bigoplus_{i \in \Z / 2\Z} \mc R_{i}^{\vee}(-\ell_{\infty}) \boxtimes \mb R^{q} p_{S \ast} ( \E_{S} (- \ell_{\infty}) \otimes \mc R_{i} \otimes \mc Q^{\vee} ) & \text{ for } p=-1,\\
\bigoplus_{i \in \Z / 2\Z} \mc R_{i}^{\vee} \boxtimes \mb R^{q} p_{S \ast} ( \E_{S} ( - \ell_{\infty}) \otimes \mc R_{i} ) & \text{ for }p=0,
\end{cases}
$$
where $p_{S} \colon X_{\kappa} \times S \to S$ is the projection.

We need the following vanishing lemma.
\begin{lem}
\label{vanish}
For $i=0,1$, we have
$$
\begin{cases}
\mb R^{q} p_{S \ast} ( \E_{S} ( - k \ell_{\infty}) \otimes \mc R_{i} ) =0 & \text{ for } k=1,2, q=0, 2\\
\mb R^{q} p_{S \ast} ( \E_{S} ( - \ell_{\infty}) \otimes  \mc R_{i} \otimes \mc Q^{\vee} ) =0 & \text{ for } q= 0, 2.\\
\end{cases}
$$
\end{lem}
\proof
This follows from \cite[Lemma 2.4]{N2}.
\endproof

From this lemma, $\mc V_{i}=\mb R^{1} p_{S \ast} ( \E_{S} (- 2 \ell_{\infty}) \otimes \mc R_{i} ), \mc V_{i}'=\mb R^{1} p_{S \ast} ( \E_{S} (-  \ell_{\infty}) \otimes \mc R_{i} )$, and $\wt{\mc W}=\mb R^{1} p_{S \ast} ( \E_{S} (-  \ell_{\infty}) \otimes \mc R_{i} \otimes \mc Q^{\vee} )$ are vector bundles on $S$.
Furthermore, the complex $E^{\bullet 1}_{2}( \ell_{\infty} )$ on $X_{\kappa} \times S$ 
\begin{eqnarray}
\label{bei2}
0 \to 
\bigoplus_{i \in \Z/ 2\Z} \mc R_{i}^{\vee} (- \ell_{\infty}) \boxtimes \mc V_{i}
\stackrel{a}{\to}
\bigoplus_{i \in \Z/ 2\Z} \mc R_{i}^{\vee} \boxtimes \wt{\mc W}
\stackrel{b}{\to}
\bigoplus_{i \in \Z/ 2\Z} \mc R_{i}^{\vee} (\ell_{\infty}) \boxtimes \mc V_{i}'
\to 0
\end{eqnarray}
satisfies $\ker a =0, \coker b =0$, and $\E_{S} \cong \ker b / \im a$.

We can write
$$
a=
\begin{bmatrix}
x_{0} \id_{\mc R_{0}} \boxtimes a_{0}^{0} & \xi_{1} \boxtimes a_{1}^{1} + \bar{\xi}_{1} \boxtimes a_{2}^{1} \\ 
\xi_{0} \boxtimes a_{1}^{0} + \bar{\xi}_{0} \boxtimes a_{2}^{0} & x_{0} \id_{\mc R_{1}} \boxtimes a_{0}^{1}
\end{bmatrix},
b=
\begin{bmatrix}
x_{0} \id_{\mc R_{0}} \boxtimes b_{0}^{0} & \xi_{1} \boxtimes b_{1}^{1} + \bar{\xi}_{1} \boxtimes b_{2}^{1} \\ 
\xi_{0} \boxtimes b_{1}^{0} + \bar{\xi}_{0} \boxtimes b_{2}^{0} & x_{0} \id_{\mc R_{1}} \boxtimes b_{0}^{1}
\end{bmatrix}.
$$
We put $\mc V= \bigoplus_{i=0}^{1} \mc V_{i},  \wt{\mc W}= \bigoplus_{i=0}^{1} \wt{\mc W}_{i}$ and $\mc V'= \bigoplus_{i=0}^{1} \mc V_{i}'$, and 
$$
a_{0} = 
\begin{bmatrix}
a_{0}^{0} & 0 \\
0 & a_{0}^{1}
\end{bmatrix},
a_{k} = 
\begin{bmatrix}
0& a_{k}^{1} \\
a_{k}^{0} & 0
\end{bmatrix},
b_{0} = 
\begin{bmatrix}
b_{0}^{0} & 0 \\
0 & b_{0}^{1}
\end{bmatrix},
b_{k} = 
\begin{bmatrix}
0& b_{k}^{1} \\
b_{k}^{0} & 0
\end{bmatrix}
$$
for $k =1,2$.
Then we have
$a_{0}, a_{1}, a_{2} \in \Hom_{S}(\mc V, \wt{\mc W})$ and $b_{0}, b_{1}, b_{2} \in \Hom_{S}(\wt{\mc W}, \mc V')$.

Over $X_{0} \setminus O \cong X_{1} \setminus C$, we can write
$$
a = a_{0} x_{0} + a_{1} x_{1} + a_{2} x_{2}, 
b = b_{0} x_{0} + b_{1} x_{1} + b_{2} x_{2}.
$$
Since $ba=0$, we have $b_{i} a_{i} = 0, b_{i}a_{i+1} + b_{i+1} a_{i} = 0$ for $i \in \Z/ 3 \Z$.
Restricting the complex \eqref{bei2} to $\ell_{\infty}$, we have
$$
0 \to \mo_{\PP^{1}} (-1) \boxtimes \mc V \stackrel{a|_{\ell_{\infty}}}{\to} \mo_{\PP^{1}}  \boxtimes \wt{\mc W} \stackrel{b|_{\ell_{\infty}}}{\to} \mo_{\PP^{1}} (1) \boxtimes \mc V' \to 0, 
$$ 
where $a|_{\ell_{\infty}} = x_{1} \boxtimes a_{1} + x_{2} \boxtimes a_{2}, b|_{\ell_{\infty}} = x_{1} \boxtimes b_{1} + x_{2} \boxtimes b_{2}$, and we regard $\mc V, \wt{\mc W}$ and $\mc V'$ as $H$-equivariant vector bundles on $S$ by $\Z/2\Z$-grading.

By the similar arguments as in \cite[\S 2.1]{N2} using framing $\Phi \colon E|_{\ell_{\infty} \times S} \cong W \otimes \mo_{\PP^{1} \times S}$, we see that $b_{1} a_{2} = - b_{2} a_{1}$ gives an isomorphism $\mc V \cong \mc V'$, and $\wt{W} = \im a_{1} \oplus \im a_{2} \oplus \mc W \cong \mc V \oplus \mc V \oplus \mc W$, where $\mc W= \ker b_{1} \cap \ker b_{2}$. 
Via these identifications, we have
$$
a_{1} = 
\begin{bmatrix}
- \id_{\mc V}\\
0\\
0
\end{bmatrix},
a_{2} = 
\begin{bmatrix}
0\\
- \id_{\mc V}\\
0
\end{bmatrix},
b_{1} =
\begin{bmatrix}
0 & - \id_{\mc V} & 0
\end{bmatrix},
b_{2} =
\begin{bmatrix}
\id_{\mc V} & 0 & 0
\end{bmatrix}.
$$
From the condition $ba=0$, we can write
$$
a_{0} = 
\begin{bmatrix}
B_{1}\\
B_{2}\\
w
\end{bmatrix},
b_{1} =
\begin{bmatrix}
-B_{2} & B_{1} & z
\end{bmatrix}
$$
with $[B_{1}, B_{2}] + zw =0$.
By \cite[Proposition 4.1]{N3} and \cite[Lemma 2.7]{N2}, these are family of $\mbi \zeta^{\kappa}$-stable ADHM data with $\mbi \zeta^{0} \in \mc C_{0}$ and $\mbi \zeta^{1} \in \mc C_{+}$. 
%
%\begin{lem}
%Suppose that a family $(B_{1} , B_{2} , z, w)$ of ADHM data on $\Z / 2 \Z$-graded vector bundles $\mc W, \mc V$ on $S=\Spec A$ are given, and define homomorphisms $a, b$ as above.
%Then $a$ is injective possible except finitely many points and $b$ is surjective at any point. 
%\end{lem}

%%%%%%%%%%%%%%%%%%%%%%%%%%%%%%%%%%%%%%%%%%%%%%%%%%%%%%%%%%%%%%%%%%%%%%%%%%%

\subsection{Isomorphisms of moduli spaces}
\label{subsec:isom}

We summarize results in the previous subsection.
If we have a family of framed sheaves $(\E_{S}, \Phi_{S})$ on $X_{\kappa} \times S$, then we have a family of $\mbi \zeta^{\kappa}$-stable ADHM data $(B_{1}, B_{2}, z, w)$ on $S$.
This defines a morphism to $S \to M^{\mbi \zeta^{\kappa}}(\mbi w, \mbi v)$ such that the pull back of the complex \eqref{bei} coincides with the complex \eqref{bei2}.
This means that the pull-back of the universal framed sheaf $(\E, \Phi)$ on $X_{\kappa} \times M^{\mbi \zeta^{\kappa}}(\mbi w, \mbi v)$ is isomorphic to $(\E_S, \Phi_S)$.
Furthermore such a morphism is unique, since isomorphisms of framed sheaves induces isomorphisms of ADHM data.
Hence we get a morphism $M_{X_{\kappa}}(c) \to M^{\mbi \zeta^{\kappa}}(\mbi w, \mbi v)$.

Together with \S A.1, we have two morphisms between $M^{\mbi \zeta^{\kappa}}(\mbi w, \mbi v)$ and $M_{X_{\kappa}}(c)$.
To check whether these are converse to each other, it is enough to see set theoretically that the compositions are identities, and this is proven in \cite[Theorem 2.2]{N3}.  
Hence $M^{\mbi \zeta^{\kappa}}(\mbi w, \mbi v)$ is a moduli of framed sheaves on $X_{\kappa}$, and we complete the proof of Theorem \ref{ale}. 

To extend Theorem \ref{ale}, we introduce $m$-stability for framed sheaves on $X_{1}$.
For the following definition, we put $\wt{C} = \ell_{\infty} - F=\frac{1}{2}C$.  

\begin{defn}
\label{m-stability}
For $m \in \Z_{\ge 0}$, a framed sheaf $(E,\Phi)$ on $X_{1}$ is said to be $m$-stable if $E(-m \widetilde{C} )$ is perverse coherent, i.e.,
\begin{enumerate}
\item $\Hom_{X_{1}}(E(-m \widetilde{C}), \mo_{C}(-1))=0$,
\item $\Hom_{X_{1}}(\mo_{C},E(-m \widetilde{C}))=0$,
\item $E(-m \widetilde{C})$ is torsion free outside $C$.
\end{enumerate} 
\end{defn}
For $c \in A(IX_{1})$, we write by $M^{m}_{X_{1}}(c)$ the moduli space of $m$-stable framed sheaves $(E, \Phi)$ on $X_{1}$ with $\wt{\ch} (E)=c$ in $A(IX_{1})$.
We can show the following theorem similarly as in \cite{NY3}, but we will give a proof elsewhere.
\begin{thm}
\label{perv}
We have an isomorphism between $M^{\mc C_{m}}(\mbi w, \mbi v)$ and $M^{m}_{X_{1}}(c)$,
where we put 
\begin{eqnarray*}
c=(w_{0} + w_{1})[X_{1}]+\left(-2v_{0} + 2v_{1} - w_{1}\right) \frac{C}{2} - \left( v_{0} + \frac{w_{1}}{4}\right)P  + (w_{0} - w_{1}) \ell_{\infty}^{1} \in A(IX_{1}).
\end{eqnarray*}
\end{thm}

By this theorem and Theorem \ref{ale}, we have $M^{m}_{X_{1}}(c) \cong M_{X_{1}}(c)$ for $m \gg 0$, and $M^{0}_{X_{1}}(c) \cong M_{X_{0}}(c)$.

%%%%%%%%%%%%%%%%%%%%%%%%%%%%%%%%%%%%%%%%%%%%%%%%%%%%%%%%%%%%%%%%%%%%%%%%%%%
%%%%%%%%%%%%%%%%%%%%%%%%%%%%%%%%%%%%%%%%%%%%%%%%%%%%%%%%%%%%%%%%%%%%%%%%%%%
%%%%%%%%%%%%%%%%%%%%%%%%%%%%%%%%%%%%%%%%%%%%%%%%%%%%%%%%%%%%%%%%%%%%%%%%%%%

\section{Combinatorial description of partition functions}
\label{sec:comb}
Following the same arguments as in \cite{NY1}, we give a combinatorial description of Nekrasov partition functions $Z_{X_{\kappa}}$ for $\kappa = 0,1$, and compare with the original Nekrasov partition function $Z$ defined from framed moduli $M(r,n)$ of torsion free sheaves on the plane $\PP^{2}$ with the rank $r$ and $c_{2}=n$.
%This computation is also written in \cite[4.6]{BPSS} for $\kappa=1$. 

In the following, we consider $M_{X_{\kappa}}(c)$ for $\kappa=0,1$, where
$$
c=r [X_{1}]+k [C] - n [P]  + \ba r [\ell_{\infty}^{1}] \in A(IX_{1}),
$$
and for $\kappa=0$, $c \in A(IX_{0})$ is defined via the equivalence $F \colon D(X_{0}) \cong D(X_{1})$ in Proposition \ref{eq}.
We put $r_{0}= (r + \ba r)/2$ and $r_{1} = (r-\ba r)/2$.

%%%%%%%%%%%%%%%%%%%%%%%%%%%%%%%%%%%%%%%%%%%%%%%%%%%%%%%%%%%%%%%%%%%%%%%%%%%%

\subsection{Fixed point sets of framed moduli}
\label{subsec:fixe}
For $\tilde{T}$-action on $M_{X_{\kappa}}(c)$ defined in \S \ref{subsec:toru}, we consider fixed point sets $M_{X_{\kappa}}(c)^{\tilde{T}}$ in $M_{X_{\kappa}}(c)$.

\begin{prop}
\label{fixedpt0}
For $c \in A(IX_{0})$ as in Theorem \ref{ale}, the set of fixed points of $M_{X_{0}}(c)$ consists of pairs
$$
\left(\bigoplus_{\alpha=1}^{r_{0}}I_{\alpha} \oplus \bigoplus_{\alpha=r_{0}+1}^{r}I_{\alpha}\left( F - \ell_{\infty} \right), \Phi \right),
$$
where $I_{\alpha}$ are ideal sheaves supported on $P$, and $\Phi$ is a direct sum of natural isomorphisms $I_{\alpha}|_{\ell_{\infty}}\cong \mo_{\ell_{\infty}}$ for $\alpha=1, \ldots, r_{0}$ and $I_{\alpha}( F - \ell_{\infty})|_{\ell_{\infty}} \cong \mo_{\ell_{\infty}}\otimes (-1)$ for $\alpha=r_{0}+1, \ldots, r_{0} + r_{1}$.
Furthermore each $I_{\alpha}$ corresponds to a Young diagram $Y_{\alpha}$ for $\alpha=1, \ldots, r$ such that 
$$
v_{s}= \sum_{\alpha = 1}^{r} \sharp \lbrace (i,j) \in Y_{\alpha} \mid l_{\alpha} + (i-1) +(j-1) \equiv s \text{ mod } 2 \rbrace,
$$
for $s=0,1$, where $l_1 = \cdots = l_{r_0} =0$ and $l_{r_{0}+1} = \cdots = l_{r_{0}+ r_{1}} =1$.
\end{prop}
\proof
For a $\tilde{T}$-fixed point $(E, \Phi) \in M_{X_{0}}(c)$, we consider $E \in \Coh X_{0}$ as a $\Z_{2}$-equivariant torsion free sheaf on $\PP^{2}$.
As in \cite[Proposition 2.9]{NY1}, we have an eigen-vector decomposition of $E$ for $\tilde{T}$-action.
Since $\tilde{T}$-action is compatible with $\Z_{2}$-action, this gives a decomposition of sheaves on $X_{0}$.
\endproof

\begin{prop}
\label{fixedpt1}
For $c \in A(IX_{1})$ as in Theorem \ref{ale}, the set of fixed point of $M_{X_{1}}(c)$ consists of pairs
$$
\left(\bigoplus_{\alpha=1}^{r}I_{\alpha}(k_{\alpha}C), \Phi \right),
$$
where $I_{\alpha}$ are ideal sheaves supported on $\lbrace P_1, P_2 \rbrace$, and $\Phi$ is a direct sum of natural isomorphisms $I_{\alpha}(k_{\alpha}C)|_{\ell_{\infty}}\cong \mo_{\PP^1} \otimes (-1)^{2k_{\alpha}}$.
Vectors $\vec{k}=(k_{1}, \ldots, k_{r} )\in \frac{1}{2}\Z^{r}$ satisfy $k_{1}, \ldots, k_{r_{0}} \in \Z$, $k_{r_0+1}, \ldots, k_{r} \in \frac{1}{2} + \Z$, and $\sum_{\alpha =1}^{r} k_{\alpha}=k$.
Furthermore each $I_{\alpha}$ corresponds to a pair of Young diagrams $(Y^1_{\alpha}, Y^{2}_{\alpha})$ for $\alpha=1, \ldots, r$ such that $\sum_{\alpha=1}^{r} \left( k^2_{\alpha}+|Y^1_{\alpha}|+|Y^2_{\alpha}| \right) = n$.
\end{prop}
\proof
$(E,\Phi)\in M_{X_{1}}(c)$ is fixed by $T^r$-action if and only if it has eigenvector decomposition $E=I_1\oplus \cdots \oplus I_{r}$ and $\Phi$ is direct sum of isomorphisms $I_i \cong \mo_{\ell_{\infty}}$ of the $i$-th factor for each $i=1, \dots, r$. 
\endproof

In the following, we compute $\tilde{T}$-modulue structures of fibers of $\tilde{T}$-equivariant vector bundles on framed moduli.
These are considered as elements of the representation ring $R(\tilde{T})$ of the torus $\tilde{T}$.
We identify it with the Laurent polynomial ring $\Z[t_{1}^{\pm}, t_{2}^{\pm}, e_{1}^{\pm}, \ldots, e_{r}^{\pm}, m_{1}^{\pm}, \ldots, m_{2r}^{\pm} ]$.
We also consider $\Z_{2}$-grading of $R(\tilde{T})$ as $\mc S$ in \S \ref{subsec:adhm2}, that is, defined by
$\deg t_{1} = \deg t_{2} = \deg e_{r_{0} + 1} = \cdots = \deg e_{r} = 1$, and $\deg e_{1} = \deg e_{r_{0}} = \deg \mu_{1} = \cdots = \deg \mu_{2r} = 0$.
For an element $F \in \mc R(\tilde{T})$, the degree $s$ part is denoted by $[F]_{s}$ for $s \in \Z_{2}$.

%For $X_{1}$, we have
%\begin{eqnarray}
%\label{comb1}
%Z_{X_{1}}^{k} (\mbi \e, \mbi a, \mbi m, q) = \sum_{(\vec{k}, \vec{Y}^{1}, \vec{Y}^{2})} \frac{\prod_{f=1}^{2r} \iota_{(\vec{k}, \vec{Y}^{1}, \vec{Y}^{2})}^{\ast} e\left( \mc V_{0} \otimes \frac{e^{m_{f}}}{\sqrt{t_{1}t_{2}}} \right)}{\iota_{(\vec{k}, \vec{Y}^{1}, \vec{Y}^{2})}^{\ast} e \left( TM_{X_{1}}(c) \right)} q^{\sum_{\alpha=1}^{r} k_{\alpha}^{2} + |Y^{1}_{\alpha}| + |Y^{2}_{\alpha}|},
%\end{eqnarray}
%where $(\vec{k}, \vec{Y}^{1}, \vec{Y}^{2})$ are data $\vec{k}=(k_{\alpha})_{\alpha = 1}^{r}, \vec{Y}^{i}=(Y^{i}_{\alpha})_{\alpha = 1}^{r}$ for $i=1,2$ as in Proposition \ref{fixedpt1}, and $\iota_{(\vec{k}, \vec{Y}^{1}, \vec{Y}^{2})} \colon \lbrace (E, \Phi)\rbrace \to M_{X_{1}}(c)$ is the embedding of the corresponding fixed point.

%%%%%%%%%%%%%%%%%%%%%%%%%%%%%%%%%%%%%%%%%%%%%%%%%%%%%%%%%%%%%%%%%%%%%%%%%%%

\subsection{$\tilde{T}$-module structures of tautological bundles on framed moduli}
\label{subsec:tild1}
We compute $\tilde{T}$-module structures of tautological bundles $\mc V_{0}= \mb R^{1} p_{\ast} \mc E (- \ell_{\infty}), \mc V_{1}=\mb R^{1} p_{\ast} \mc E (-F)$ over $\tilde{T}$-fixed points of framed moduli,
where $\mc E$ are universal sheaves on $X_{\kappa} \times M_{X_{\kappa}}(\alpha)$, and $p \colon X_{\kappa} \times M_{X_{\kappa}}(\alpha) \to M_{X_{\kappa}}(\alpha)$ is the projection.

\begin{prop}
\label{tautological0}
For  a fixed point $(E, \Phi) \in M_{X_{0}}(c)$ corresponding to a datum $\vec{Y} = ( Y_{1}, \ldots, Y_{r} )$, we have isomorphisms $\mc V_{s}|_{(E, \Phi)} \cong \bigoplus_{\alpha = 1}^{r} \bigoplus_{(i,j) \in Y_{\alpha}} [e_{\alpha} t_{1}^{-i+1} t_{2}^{-j+1}]_{s}$ of $\tilde{T}$-modules, where $ [e_{\alpha}t_{1}^{-i+1}t_{2}^{-j+1}]_{s}$ is the degree $s$ parts of $e_{\alpha}t_{1}^{-i+1}t_{2}^{-j+1}$ in $R(\tilde{T})$ for $s=0,1$.
\end{prop}
\proof
We compute $\mc V_{s}|_{(E, \Phi)} = H^{1}(X_{0}, E( s(\ell_{\infty} - F) - \ell_{\infty}))$ as follows.
We have a decomposition $E= \bigoplus_{\alpha = 1}^{r} e_{\alpha} I_{Y_{\alpha}}$.
For each $\alpha$, we consider an exact sequence
$$
0 \to \mo_{Z_{\alpha}} (  s(\ell_{\infty} - F)) \to I_{Y_{\alpha}} (  s(\ell_{\infty} - F) - \ell_{\infty}) \to \mo_{X_{0}}(  s(\ell_{\infty} - F) - \ell_{\infty} ) \to 0,
$$
where $Z_{\alpha}$ is a $0$-dimensional sub-scheme of $\PP^{2}$ defined by $I_{Y_{\alpha}}$.
Hence we have $H^{1} (X_{0}, I_{Y_{\alpha}}) \cong H^{0}(X_{0}, \mo_{Z_{\alpha}}( s(\ell_{\infty} - F) - \ell_{\infty}))$.
This is the space of $\Z_{2}$-invariant sections of $\mo_{Z_{\alpha}} \otimes (-1)^{s}$, where $(-1)^{s}$ denote representations $H=\lbrace \pm 1 \rbrace \to \C^{\ast}, m \mapsto m^{s}$ for $s=0,1$.
This gives the assertion.
\endproof

For tautological bundles on the minimal resolution $X_{1}$, we need the following computations.
For $k \in \frac{1}{2}\Z$, we consider a $T^{2}$-equivariant sheaf $\mo_{X_{1}}(kC-\ell_{\infty})$ and a $T^2$-modules $L_{k}(t_1,t_2)=H^1(X_{1}, \mo_{X_{1}}(kC-\ell_{\infty}))$.
% and $L'_k(t_1,t_2)=H^1(X_{1}, \mo_X_{1}(kC-F))$

\begin{lem}
We have the following isomorphisms of $T^2$-modules.
\begin{eqnarray*}
L_{k} (t_1, t_2) 
\cong
\begin{cases}
{\displaystyle
  \bigoplus_{\substack{i,j \ge 0\\ i+j \le 2k-2\\ i+j \equiv 2k \text{ mod }2}}
  }
         t_1^{i+1}t_2^{j+1} & \text{if }k>\frac{1}{2}\\
{\displaystyle
  \bigoplus_{\substack{i,j \ge 0\\ i+j \le -2k-2\\ i+j \equiv 2k \text{ mod }2}}
  }
t_1^{-i}t_2^{-j} & \text{if }k<-\frac{1}{2}\\
0 & \text{otherwise}
\end{cases}
\end{eqnarray*}
\end{lem}
\proof
%(1) 
For $k=0, \pm \frac{1}{2}$, since $\chi(\mo_{X_{1}}(kC-\ell_{\infty}))=0$ we have $H^1(X_{1},\mo_{X_{1}}(kC-\ell_{\infty}))=0$ by Lemma \ref{vanish}. 
For $k>\frac{1}{2}$ we consider the exact sequence
$$
0 \to \mo_{X_{1}}((k-1)C-\ell_{\infty}) \to \mo_{X_{1}}(kC-\ell_{\infty}) \to \mo_{C}(kC)=\mo_{\PP^1}(-2k) \to 0.
$$
From the cohomology long exact sequence
$$
0 \to H^1(X_{1}, \mo_{X_{1}}((k-1)C-\ell_{\infty})) \to H^1(X_{1}, \mo_{X_{1}}(kC-\ell_{\infty})) \to H^1(\PP^1, \mo_{\PP^1}(-2k)) \to 0
$$ 
we get a decomposition $H^1(X_{1}, \mo_{X_{1}}(kC-\ell_{\infty})) = H^1(X_{1}, \mo_{X_{1}}((k-1)C-\ell_{\infty})) \oplus H^1(\PP^1, \mo_{\PP^1}(-2k))$.
Since we have the $T^2$-equivariant dualizing sheaf $t_1^{-1}t_2^{-1} \mo_{\PP^1}(-2)$, by the Serre duality we have 
$$
H^1(\PP^1, \mo_{\PP^1}(-2k))=(t_1^{-1}t_2^{-1}H^0(\PP^1, \mo_{\PP^1}(2k-2))^{\ast}
=\bigoplus_{i,j\ge0, i+j=2k-2}t_1^{i+1}t_2^{j+1}.
$$ 
Repeating this we get the assertion.

For $k<-\frac{1}{2}$, we consider the exact sequence
$$
0 \to \mo_{X_{1}}(kC-\ell_{\infty}) \to \mo_{X_{1}}((k+1)C-\ell_{\infty}) \to \mo_{C}((k+1)C)=\mo_{\PP^1}(-2k-2) \to 0.
$$
We get a decomposition 
$$
H^1(X_{1}, \mo_{X_{1}}(kC-\ell_{\infty})) = H^1(X_{1}, \mo_{X_{1}}((k+1)C-\ell_{\infty})) \oplus H^0(\PP^1, \mo_{\PP^1}(-2k-2)).
$$
Repeating this procedure we get the assertion.
%(2) It is similar to the proof of (1), or it can be obtained directly from (1) since $\mo_{X_{1}}(kC -F)=\mo_{X_{1}} \left( (k+\frac{1}{2}) C - \ell_{\infty} \right)$.
\endproof

\begin{prop}
\label{tautological1}
%We have the following.\\
%(1) 
For a $\tilde{T}$-fixed point $\left( E, \Phi \right) \in M_{X_{1}}(c)$ corresponding to a datum $(\vec{k}, \vec{Y}^1, \vec{Y}^2)$, the $\tilde{T}$-module $\mc V_{s}|_{(E,\Phi)}$ is isomorphic to 
\begin{eqnarray*}
\bigoplus_{\alpha=1}^{r} e_{\alpha} \left(
L_{k_{\alpha} + \frac{s}{2}} (t_1, t_2) \oplus
\bigoplus_{(i,j)\in Y^1_{\alpha}}t_1^{2(k_{\alpha}-i+1+\frac{s}{2})}\left(\frac{t_2}{t_1}\right)^{-j+1} \oplus 
\bigoplus_{(i,j)\in Y^2_{\alpha}}\left(\frac{t_1}{t_2}\right)^{-i+1} t_2^{2(k_{\alpha}-j+1+\frac{s}{2})} \right)
\end{eqnarray*}
for $s=0,1$.
\end{prop}
\proof
%(1) 
The torsion free sheaf $E$ is decomposed into $\bigoplus_{\alpha=1}^{r}I_{\alpha}(k_{\alpha}C)$.
The $T^{r}$-action on each component $H^{1}(X_{1}, I_{\alpha}\left( \left( k_{\alpha} + \frac{s}{2} \right) C-\ell_{\infty} \right) )$ is given by a multiplication of $e_{\alpha}$ for $\alpha=1, \ldots, r$. 
Hence it is enough to compute a $T^2$-module structure on $H^{1}(X_{1}, I_{\alpha}\left( \left( k_{\alpha} + \frac{s}{2} \right)C-\ell_{\infty} \right) )$ induced by the natural $T^2$-equivariant structure on $I_{\alpha} \left( \left( k_{\alpha} + \frac{s}{2} \right)C-\ell_{\infty} \right)$.
By the exact sequence
$$
0 \to I_{\alpha} \left( \left( k_{\alpha} + \frac{s}{2} \right)C-\ell_{\infty} \right) \to \mo_{X_{1}} \left( \left( k_{\alpha} + \frac{s}{2} \right)C-\ell_{\infty} \right) \to \mo_{Z_{\alpha}} \left( \left( k_{\alpha} + \frac{s}{2} \right)C-\ell_{\infty} \right) \to 0,
$$
we get a decomposition $H^1(X_{1},I_{\alpha} \left( \left( k_{\alpha} + \frac{s}{2} \right) C-\ell_{\infty} \right) )=L_{k_{\alpha}}(t_1,t_2) \oplus H^0(X_{1}, \mo_{Z_{\alpha}}(k_{\alpha}C-\ell_{\infty}))$.
We have $Z=Z_{\alpha}^1 \amalg Z_{\alpha}^2$, where $Z_{\alpha}^i$ is the sub-scheme supported at $P_i$ corresponding to $Y_{\alpha}^i$.
The multiplication of $yx_i^{2k}/x_0^{2k-1}$ gives an {\it equivariant} isomorphism $\mo_{Z_{\alpha}^i}(k_{\alpha}C-\ell_{\infty}) \cong t_i^{2k_{\alpha}}\mo_{Z_{\alpha}^i}$ for $i=1,2$.
Hence we have the desired isomorphism 
\begin{eqnarray*}
&&
H^0 \left( X_{1},\mo_{Z_{\alpha}} \left( \left(k_{\alpha} - \frac{s}{2} \right) C-\ell_{\infty} \right) \right) \\
&\cong&
\bigoplus_{(i,j)\in Y^1_{\alpha}}t_1^{2(k_{\alpha}-i+1+\frac{s}{2})}\left(\frac{t_2}{t_1}\right)^{-j+1} \oplus 
\bigoplus_{(i,j)\in Y^2_{\alpha}}\left(\frac{t_1}{t_2}\right)^{-i+1} t_2^{2(k_{\alpha}-j+1+\frac{s}{2})}.
\end{eqnarray*}
%(2) It is similar to the proof of (1).
\endproof

%%%%%%%%%%%%%%%%%%%%%%%%%%%%%%%%%%%%%%%%%%%%%%%%%%%%%%%%%%%%%%%%%%%%%%%%%%%

\subsection{$\tilde{T}$-module structures of tangent bundles on framed moduli}
\label{subsec:tild2}
We also compute the $\tilde{T}$-module structure of the tangent bundle of framed moduli $M_{X_{\kappa}}(c)$.
Let $Y_{\alpha}=\lbrace \lambda_{\alpha, 1}, \lambda_{\alpha,2}, \cdots, \rbrace$ be a Young diagram where $\lambda_{\alpha,i}$ is the height of the $i$-th column.
We set $\lambda_{\alpha, i}=0$ when $i$ is larger than the width of the diagram $Y_{\alpha}$.
Let $Y_{\alpha}^{T}=\lbrace \lambda_{\alpha, 1}', \lambda_{\alpha,2}', \cdots \rbrace $ be its transpose.
For a box $s=(i,j)$ in the $i$-th column and the $j$-th row, we define its arm-length $a_{Y_{\alpha}}(s)$ and leg-length $l_{Y_{\alpha}}(s)$ with respect to the diagram $Y_{\alpha}$ by $a_{Y_{\alpha}}(s)=\lambda_{\alpha, i}-j$ and $l_{Y_{\alpha}}(s)=\lambda_{\alpha, j}'-i$. 

We consider framed moduli $M(r,n)$ of torsion free sheaves on the plane $\PP^{2}$ with the rank $r$ and $c_{2}=n$.
We recall from \cite[Theorem 2.11]{NY1} that the fibre of $TM(r,n)$ over a fixed point corresponding to a datum $\vec{Y_{1}} = ( Y_{1}, \ldots, Y_{r} )$ consisting of Young diagrams is isomorphic to $\bigoplus_{\alpha,\beta=1}^{r} N_{\alpha,\beta}(t_1,t_2)$ as $\tilde{T}$-modules, where 
$$
N_{\alpha, \beta}(t_1,t_2)=e_{\beta}e_{\alpha}^{-1} \times \left\lbrace \bigoplus_{s\in Y_{\alpha}}\left(t_1^{-l_{Y_{\beta}}(s)}t_2^{a_{Y_{\alpha}}(s)+1}\right) \oplus 
\bigoplus_{t \in Y_{\beta}}\left(t_1^{l_{Y_{\alpha}}(t)+1}t_2^{-a_{Y_{\beta}}(t)}\right) \right\rbrace.
$$

\begin{prop}
\label{tangent0}
The fibre of $TM_{X_{0}}(c)$ over a fixed point corresponding to a datum $\vec{Y} = ( Y_{1}, \ldots, Y_{r} )$ is isomorphic to $\tilde{T}$-modules $\bigoplus_{\alpha,\beta=1}^{r} [N_{\alpha,\beta}(t_1,t_2)]_{0}$, where $ [N_{\alpha,\beta}(t_1,t_2)]_{0}$ is the degree $0$ parts of $N_{\alpha,\beta}(t_1,t_2)$ in $R(\tilde{T})$.
\end{prop}
\proof
Let $(E, \Phi)$ be a $\tilde{T}$-fixed point corresponding to $\vec{Y}$. 
Then the $\tilde{T}$-module structure of $T_{(E,\Phi)}M_{X_{0}}(c)=\Ext_{X_{0}}(E, E(-\ell_{\infty}))$ is computed similarly as in \cite[Theorem 2.11]{NY1}.
But in addition we must take $\Z_{2}$-invariant sections.
In this way we get the assertion.
\endproof

We also compute the $\tilde{T}$-module structure of the tangent bundle of $M_{X_{1}}(c)$.

\begin{prop}
\label{tangent1}
The fibre of the tangent bundle $TM_{X_{1}}(c)$ over a fixed point corresponding to a datum $(\vec{k}, \vec{Y}^{1}, \vec{Y}^{2})$ is isomorphic to 
$$
\bigoplus_{\alpha, \beta=1}^{r} \left(
e_{\beta}e_{\alpha}^{-1}L_{k_{\beta}-k_{\alpha}}(t_1,t_2) \oplus t_1^{2k_{\beta}-2k_{\alpha}} M^1_{\alpha, \beta}(t_1, t_2) \oplus t_2^{2k_{\beta}-2k_{\alpha}} M^2_{\alpha, \beta}(t_1, t_2) \right)$$
as a $\tilde{T}$-module, where  
$M^1_{\alpha, \beta}(t_1, t_2)$ (resp. $M^{2}_{\alpha,\beta}(t_1, t_2)$) is equal to $N_{\alpha, \beta}(t_1^2, t_2/t_1)$ (resp. $N_{\alpha, \beta}(t_1/t_2, t_2^2)$), with $(Y_{\alpha}, Y_{\beta})$ replaced by $(Y^{1}_{\alpha}, Y^{1}_{\beta})$ (resp. $(Y^{2}_{\alpha}, Y^{2}_{\beta})$).
\end{prop}
\proof
Let $\Ext^{\ast}_{X_{1}}$ be the alternating sum $\sum (-1)^{i}\Ext^{i}_{X_{1}}$ of $\tilde{T}$-modules.
By Lemma \ref{smooth}, we have 
\begin{eqnarray*}
T_{(E, \Phi)}M_{X_{1}}(c) 
&=& -\Ext^{\ast}_{X_{1}}(E,E(-\ell_{\infty}))\\
&=& \bigoplus_{\alpha, \beta =1}^{r} -\Ext^{\ast}_{X_{1}}\left(I_{\alpha}(k_{\alpha}C), I_{\beta}(k_{\beta}C-\ell_{\infty})\right).
\end{eqnarray*}
Each summand is multiplied by $e_{\beta}e_{\alpha}^{-1}$ for $T^r$-action.
In the rest of proof we compute the $T^2$-action on each summand.

By the exact sequence $0 \to I_{\alpha} \to \mo_{X_{1}} \to \mo_{Z_{\alpha}} \to 0,$
we get the following decomposition of $T_{(E, \Phi)}M_{X_{1}}(c)$:
\begin{eqnarray}
% \bigoplus_{\alpha,\beta=1}^{r} -\Ext^{\ast}_{X_{1}}\left( I_{\alpha}(k_{\alpha}C), I_{\beta}(k_{\beta}C-\ell_{\infty})\right)= \notag\\
\bigoplus_{\alpha,\beta=1}^{r}\left( -\Ext^{\ast}_{X_{1}}\left( \mo_{X_{1}}(k_{\alpha}C), \mo_{X_{1}}(k_{\beta}C-\ell_{\infty})\right) + \Ext^{\ast}_{X_{1}}\left( \mo_{X_{1}}(k_{\alpha}C), \mo_{Z_{\beta}}(k_{\beta}C-\ell_{\infty})\right) \right. \notag\\
\label{decomp2}
\left. +\Ext^{\ast}_{X_{1}}\left( \mo_{Z_{\alpha}}(k_{\alpha}C), \mo_{X_{1}}(k_{\beta}C-\ell_{\infty})\right) - \Ext^{\ast}_{X_{1}}\left( \mo_{Z_{\alpha}}(k_{\alpha}C), \mo_{Z_{\beta}}(k_{\beta}C-\ell_{\infty}) \right) \right) .
\end{eqnarray}
The first component in \eqref{decomp2} is isomorphic to $\bigoplus_{\alpha,\beta=1}^{r} L_{k_{\beta}-k_{\alpha}}(t_1, t_2)$.

For $\alpha= 1, \ldots, r$, we have $Z_{\alpha}=Z_{\alpha}^{1} \amalg Z_{\alpha}^{2}$, where $Z_{\alpha}^{i}$ are closed sub-schemes supported at $P_i$ corresponding to $Y_{\alpha}^{i}$ for $i=1,2$.
By an equivariant isomorphism $\mo_{Z_{\alpha}^{i}}(kC)\cong t_i^{2k} \mo_{Z_{\alpha}}$, the last three terms are isomorphic to 
\begin{eqnarray} 
\bigoplus_{i=1,2}\bigoplus_{\alpha,\beta=1}^{r}t_i^{2k_{\beta}-2k_{\alpha}} \left( \Ext^{\ast}_{X_{1}}\left( \mo_{X_{1}}, \mo_{Z_{\beta}^{i}}(-\ell_{\infty})\right) + \Ext^{\ast}_{X_{1}}\left( \mo_{Z_{\alpha}^{i}}, \mo_{X_{1}} \left( -\ell_{\infty} \right) \right) \right. \notag\\ 
\label{derived}
\left. - \Ext^{\ast}_{X_{1}} \left( \mo_{Z_{\alpha}^{i}}, \mo_{Z_{\beta}^{i}} \left( -\ell_{\infty}\right) \right) \right).
\end{eqnarray}

We consider these components as derived functors from the category of coherent sheaves supported at the origin of $\C^2$ via the coordinate $(y_1, \frac{x_1}{x_2})$ and $(\frac{x_2}{x_1}, y_2)$ around $P_1$ and $P_2$ respectively.
Then we have $\Ext^{\ast}_{X_{1}}=\Ext^{\ast}_{\PP^2}$.
If we write by $I_{\alpha}^{i}$ the ideal sheaf of $Z_{\alpha}^{i}$ in $\PP^2$, then \eqref{derived} is isomorphic to
$$
\bigoplus_{i=1,2}\bigoplus_{\alpha,\beta=1}^{r}t_i^{2k_{\beta}-2k_{\alpha}} \left( \Ext^{\ast}_{\PP^2}\left( \mo_{\PP^2}, \mo_{\PP^2}(-\ell_{\infty})\right) - \Ext^{\ast}_{\PP^2}\left( I_{\alpha}, I_{\beta}(-\ell_{\infty}) \right) \right).
$$
Since $\Ext_{\PP^2}^{\ast}\left(\mo_{\PP^2}, \mo_{\PP^2}(-\ell_{\infty})\right) = 0$ and $\Ext^{i}_{\PP^2}\left( I_{\alpha}, I_{\beta}(-\ell_{\infty}) \right)=0$ for $i=0,2$, it is isomorphic to 
$$
\bigoplus_{i=1,2}\bigoplus_{\alpha,\beta=1}^{r} t_i^{2k_{\beta}-2k_{\alpha}} \Ext^{1}_{\PP^2}\left( I_{\alpha}, I_{\beta}(-\ell_{\infty}) \right) = \bigoplus_{i=1,2}\bigoplus_{\alpha,\beta=1}^{r}t_i^{2k_{\beta}-2k_{\alpha}}  M^{i}_{\alpha,\beta}(t_1,t_2)
$$
as desired.
\endproof

%%%%%%%%%%%%%%%%%%%%%%%%%%%%%%%%%%%%%%%%%%%%%%%%%%%%%%%%%%%%%%%%%%%%%%%%%%%

\subsection{Comparison to $Z_{\PP^{2}}(\mbi \e, \mbi a, \mbi m, q)$}
\label{subsec:comp1}
We consider a $\tilde{T}$-equivariant bundle
$$
\mc F_{r}(\mc V_{0})=\left( \mc V_{0} \otimes \frac{e^{m_{1}}}{\sqrt{t_{1}t_{2}}} \right) \oplus \cdots \oplus \left( \mc V_{0} \otimes \frac{e^{m_{2r}}}{\sqrt{t_{1}t_{2}}} \right)
$$
on $M_{X_{\kappa}}(\alpha)$,
where $(e^{m_{1}}, \ldots, e^{m_{2r}})$ is an element in the last component $T^{2r}$ of $\tilde{T}$.
Here we consider a homomorphism $\tilde{T}' =\tilde{T} \to \tilde{T}$ defined by
$$
(t_{1}', t_{2}', e^{\mbi a'}, e^{\mbi m'}) \mapsto (t_1, t_{2}, e^{\mbi a}, e^{\mbi m})=((t_{1}')^2, (t_{2}')^2, e^{\mbi a'}, e^{\mbi m'}), 
$$
and use identification $t_{1}'=\sqrt{t_{1}}, t_{2}'=\sqrt{t_{2}}$ and $A^{\ast}_{\tilde{T}'}(\text{pt}) \otimes \mc S \cong \mc S$.
Nekrasov partition functions are defined by
\begin{eqnarray*}
Z_{X_{\kappa}}^{k}(\mbi \e, \mbi a, \mbi m, q)&=&\sum_{\substack{\alpha \in K(X_{\kappa})\\ \mbi r(\alpha)=\mbi r, k(\alpha)=k}} q^{n(\alpha)} \int_{M_{X_\kappa}(\alpha)} e(\mc F_{r}(\mc V_{0})).
\end{eqnarray*}
as in the introduction.

We consider the other Nekrasov partition function
$$
Z_{\PP^{2}} (\mbi \e, \mbi a, \mbi m, q)= \sum_{n=0}^{\infty} q^{n} \int_{M(r,n)} e(\mc F_{r}(\mc V)) = \sum_{\vec{Y}} q^{|\vec{Y}|} \frac{\iota_{\vec{Y}}^{\ast} e(\mc F_{r}(\mc V)) }{ e(TM(r,n))}\in \mc S[[q]]
$$
where $M(r,n)$ is framed moduli of torsion free sheaves on the plane $\PP^{2}$ with the rank $r$ and $c_{2}=n$.
We consider the equivariant Euler class $e(\mc F_{r}(\mc V))$ in $\mc S$ of a $\tilde{T}$-equivariant vector bundle 
\begin{eqnarray}
\label{fr}
\mc F_{r}(\mc V)= \left( \mc V \otimes \frac{e^{m_{1}}}{\sqrt{t_{1}t_{2}}} \right) \oplus \cdots \oplus \left( \mc V \otimes \frac{e^{m_{2r}}}{\sqrt{t_{1}t_{2}}} \right),
\end{eqnarray}
where $\mc V$ is the tautological bundle on $M(r,n)$.

We write $Z_{\PP^{2}} (\mbi \e, \mbi a, \mbi m, q) = \sum_{n=0}^{\infty} \alpha_{n} q^{n}$, where we also have description of $\alpha_{n}$ in terms of Young diagrams \cite{NY1}.
Then by Proposition \ref{tautological1} and Proposition \ref{tangent1}, we have
\begin{eqnarray}
\label{mult}
Z_{X_{1}}^{k} (\mbi \e, \mbi a, \mbi m, q) 
= 
\sum_{\vec{k} \in \mc K(\mbi w, k)} q^{\sum k_{\alpha}^{2}}  \ell_{\vec{k}}(\mbi \e, \mbi a, \mbi m) Z_{\PP^{2}}(\mbi \e^{0}, \mbi a^{0}, \mbi m, q)  Z_{\PP^{2}}(\mbi \e^{1}, \mbi a^{1}, \mbi m, q)
\end{eqnarray}
for $k \in \frac{1}{2} \Z$ as in \cite[(3.5)]{IMO}.
Here 
$
\mc K(\mbi w, k) = \left\lbrace \vec{k} \in \Z^{w_{0}} \oplus \left( 1/2 + \Z \right)^{w_{1}} \ \Big| \  \sum_{\alpha=1}^{r} k_{\alpha} = k \right\rbrace, 
$
$$
\mbi \e^{0} =( 2 \e_{1}, -\e_{1}+\e_{2} ) , \mbi \e^{1} = ( \e_{1} - \e_{2} , 2 \e_{2}), \mbi a^{0}= \mbi a + 2 \e_{1} \vec{k}, \mbi a^{1} = \mbi a + 2 \e_{2} \vec{k},
$$
and 
$$
\ell_{\vec{k}}(\mbi \e, \mbi a, \mbi m)
=
\frac{\prod_{f=1}^{2r} e \left( \bigoplus_{\alpha=1}^{r} e_{\alpha} L_{k_{\alpha}}(t_{1}, t_{2}) \frac{e^{m_{f}}}{\sqrt{t_{1}t_{2}}} \right)}
{\prod_{\alpha, \beta = 1}^{r} e \left( e_{\beta} e_{\alpha}^{-1} L_{k_{\beta} - k_{\alpha}}(t_{1}, t_{2}) \right)}.
$$

%%%%%%%%%%%%%%%%%%%%%%%%%%%%%%%%%%%%%%%%%%%%%%%%%%%%%%%%%%%%%%%%%%%%%%%%%%%

\subsection{Comparison to Ito-Maruyoshi-Okuda}
\label{subsec:comp2}
Ito-Maruyoshi-Okuda \cite{IMO} introduced a similar partition functions
$$
Z^{\C^{2}/\Z_{2}}_{N_{f}=2N, \text{inst}, c}(\vec{a}, \vec{I}; \mbi \mu; q; \e_{1}, \e_{2}), Z^{A_{1}\text{-resolved}}_{N_{F}=2N}(\vec{a}, \vec{I}; \mbi \mu; q;  \e_{1}, \e_{2}).
$$
We substitute $\e_{1}=-\e_{1}, \e_{2}=-\e_{2}$, $\vec{a}=\mbi a$, $\mu_{i} = m_{i} - \frac{\e_{1}+\e_{2}}{2}, \mu_{r+i} =- m_{r+i} + \frac{\e_{1}+\e_{2}}{2}$ for $i=1, \ldots, r$, $c=-k$, $N=r$, and 
$$
\vec{I} =( \overbrace{0, \ldots, 0}^{r_{0}}\overbrace{1, \ldots, 1}^{r_{1}}).
$$
Then we have
$$
Z_{X_{0}}^{k}(\mbi \e, \mbi a, \mbi m, q)=Z^{\C^{2}/\Z_{2}}_{N_{f}=2N, \text{inst}, c}(\vec{a}, \vec{I}; \mbi \mu; q; \e_{1}, \e_{2}), Z_{X_{1}}^{k}(\mbi \e, \mbi a, \mbi m, q)=Z^{A_{1}\text{-resolved}}_{N_{F}=2N}(\vec{a}, \vec{I}; \mbi \mu; q;  \e_{1}, \e_{2}).
$$
Furthermore, after this substitution their proposed relations \cite[(4.1)]{IMO} coincides with Theorem \ref{main}.

\begin{rem}
These computations in Appendix B can be justified by \cite{N3} and earlier results without using framed moduli spaces constructed in \cite{BPSS} (cf. \cite[line 4 - 6 in p.1179]{BPSS}).
\end{rem}
%%%%%%%%%%%%%%%%%%%%%%%%%%%%%%%%%%%%%%%%%%%%%%%%%%%%%%%%%%%%%%%

%\subsection{$\mu$ map}
%
%We consider the class $\mu(C)$ on $A^{\ast}( M_{\kappa}(c) )$ defined by $\mu(C) = c_{1}(\det \mc V_{0}^{\vee} \otimes \det \mc V_{1})$ (cf. \cite[\S 4.1]{B}).
%
%For a universal sheaf $\E$ on $X_{1} \times M_{X_{1}}(c)$, we put $\Delta(\E) = c_{2}(\E) - \frac{r-1}{2r} (c_{1}(\E))^{2} \in A^{\ast}_{\tilde{T}}(X_{1} \times M_{X_{1}}(c))$, and 
%$$
%\mu(C) = \int_{X_{1}} \Delta(\E) \cap [C \times M(n)] \in A^{\ast}_{\tilde{T}}(M_{X_{1}}(c)),
%$$
%where $\int_{X_{1}}$ is the push-forward by the projection $X_{1} \times M_{X_{1}}(c) \to M_{X_{1}}(c)$.
%
%\begin{prop}
%We have 
%$$
%\iota_{(k, \vec{Y}^{1}, \vec{Y}^{2})}^{\ast} \mu(C)=
%2 (|\vec{Y}^{1}| \e_{1} + |\vec{Y}^{2}| \e_{2} +k (a_{1} - a_{2}) + k^{2} (\e_{1} + \e_{2})). 
%$$
%\end{prop}

%%%%%%%%%%%%%%%%%%%%%%%%%%%%%%%%%%%%%%%%%%%%%%%%%%%%%%%%%%%%%%%%%%%%%%%%%%%
%%%%%%%%%%%%%%%%%%%%%%%%%%%%%%%%%%%%%%%%%%%%%%%%%%%%%%%%%%%%%%%%%%%%%%%%%%%
%%%%%%%%%%%%%%%%%%%%%%%%%%%%%%%%%%%%%%%%%%%%%%%%%%%%%%%%%%%%%%%%%%%%%%%%%%%

\noindent
              R. Ohkawa \\
              Waseda University, 3--4--1 Okubo, Shinjuku-ku, Tokyo 169--8555, Japan\\
              ohkawa.ryo@aoni.waseda.jp

\end{document}